\documentclass[final]{siamart1116}
\usepackage{graphics,graphicx,epstopdf}

\usepackage{longtable}
\usepackage[tight]{subfigure}
\usepackage{float}
\usepackage{makecell}

\usepackage[leqno]{amsmath}
\usepackage{amsxtra, amsfonts,amscd,amssymb,graphicx,epsf}
\usepackage{bbm, dsfont}

\usepackage{url}

\usepackage[algo2e, vlined,ruled]{algorithm2e}

\numberwithin{equation}{section}

\usepackage{xcolor}
\usepackage[normalem]{ulem}

\usepackage{amsmath}
\usepackage{float}
\usepackage{multirow}

\usepackage[algo2e,vlined]{algorithm2e}

\usepackage{blkarray}
\usepackage{cases}
\usepackage{graphicx}
\usepackage{epstopdf}
\usepackage{color}

\begin{document}

\title{ Adaptive Low-Nonnegative-Rank Approximation for State Aggregation of Markov Chains }

\author{
  Yaqi Duan\thanks{Department of Operations Research and Financial Engineering, Princeton University,
USA (\email{yaqid@princeton.edu}).}
  \and
  Mengdi Wang\thanks{Department of Operations Research and Financial Engineering, Princeton University,
USA (\email{mengdiw@princeton.edu}).}
  \and
  Zaiwen Wen\thanks{Beijing International Center for Mathematical
Research, Peking University, CHINA (\email{wenzw@pku.edu.cn}).  Research supported in part by NSFC grants 11831002 and 11421101, and by the National Basic Research Project under grant 2015CB856002.}
\and
	Yaxiang Yuan\thanks{State Key Laboratory of Scientific and Engineering Computing, Academy of Mathematics and Systems Science, Chinese Academy of Sciences,
		China (\email{yyx@lsec.cc.ac.cn}). }
}

\maketitle

\begin{abstract}
This paper develops a low-nonnegative-rank approximation method to identify the state aggregation structure of a finite-state Markov chain under an assumption that the state space can be mapped into a handful of meta-states. The number of meta-states is characterized by the nonnegative rank of the Markov transition matrix. Motivated by the success of the nuclear norm relaxation in low rank minimization problems, we propose an atomic regularizer as a convex surrogate for the nonnegative rank and formulate a convex optimization problem. Because the atomic regularizer itself is not computationally tractable, we instead solve a sequence of problems involving a nonnegative factorization of the Markov transition matrices by using the proximal alternating linearized minimization method.  Two methods for adjusting the rank of factorization are developed so that local minima are escaped. One is to append an additional column to the factorized matrices, which can be interpreted as an approximation of a negative subgradient step.  The other is to reduce redundant dimensions by means of linear combinations.  Overall, the proposed algorithm very likely converges to the global solution. The efficiency and statistical properties of our approach are illustrated on synthetic data. We also apply our state aggregation algorithm on a Manhattan transportation data set and make extensive comparisons with an existing method.
\end{abstract}


\begin{AM}    65K05, 90C06, 90C40 \end{AM}

\section{Introduction} \label{Section1}

Markov chain is a basic model of stochastic process. As a variant of Markov chain, Markov decision process lies at the core of dynamic programming and reinforcement learning \cite{DPandOptimalControl,AbstractDP,SuttonRL}, and has wide applications in
engineering systems, operations research, artificial intelligence and computer
games. {The present reinforcement learning algorithms may perform poorly if the state space is of huge ambient dimension.} Fortunately, empirical experiences tell us that stochastic systems in real world can usually be characterized by a handful of key features. Therefore, we consider developing a computational method to identify reduced-dimension representations of a Markov chain, so that decision problems can be solved efficiently with a compressed state space.

In this paper, we are particularly interested in situations where the stochastic system is driven by a latent Markov chain with a much smaller state space. The states in the latent Markov chain are {\it soft aggregations} of the states in the original system. To be more specific, suppose that the original stochastic system has a state space $\mathcal{S} = \{ s_1, s_2, \cdots, s_d \}$ with ambient dimension $d$, and the latent Markov chain has a state space $\tilde{\mathcal{S}} = \{\tilde{s}_1, \tilde{s}_2, \cdots, \tilde{s}_r\}$ with intrinsic dimension $r$, where $r \ll d$.
If the system is now at $s_i \in \mathcal{S}$, it first chooses a meta-state in $\tilde{\mathcal{S}}$ according to some multinomial distribution $[u_{i1}, u_{i2}, \cdots, u_{ir}]$, $i.e.$, $\tilde{s}_k$ is selected with probability $u_{ik}$.
Governed by a transition matrix $\tilde{P}\in\mathbb{R}^{r\times r}$, a one-step transition in $\tilde{\mathcal{S}}$ is then conducted, resulting in a new meta-state $\tilde{s}_l$ with probability $\tilde{p}_{kl}$. From there we finally come back to $\mathcal{S}$ according to another multinomial distribution $[v_{1l}, v_{2l}, \cdots, v_{dl}]$.
In this way, the transition probabilities of the original system can be expressed by
\begin{equation*} \label{introduction1} p_{ij} = \sum_{k,l=1}^r u_{ik} \tilde{p}_{kl} v_{jl}, \quad \forall s_i,s_j \in \mathcal{S}. \end{equation*}
In a more compact matrix form,
\begin{equation} \label{introduction2} P = U \tilde{P} V^T, \quad U,V \in \mathbb{R}^{d \times r}, \tilde{P} \in \mathbb{R}^{r \times r}. \end{equation}
Note that the decomposition is not unique.

Based on the state aggregation model above, our goal is then to recover $U$, $V$ and $\tilde{P}$ of a Markov chain (up to linear transformation). In many real applications, however, the exact dynamics of the system is never revealed to us, and we only have access to realized trajectories of state-to-state transitions.
Our proposed algorithm just requires an approximation of the true probability transition matrix $P$ as input.

\subsection{Literature Review}

Recently, there have been interests in the compression of Markov chains.
In \cite{2018arXiv180202920Z, 2018arXiv180400795L}, the authors propose methods for low-rank approximations of Markov chains and provide theoretical guarantees.
It is proved in \cite{2018arXiv180202920Z} that the spectral method performs
well in the recovery of hard state aggregation structures for ``lumpable complex
networks". In this special case, the group membership indicator vectors are
mutually orthogonal, so that the singular value decomposition (SVD) result coincides with the state aggregation
structure. However, for a more general soft state aggregation model
\eqref{introduction2}, the orthogonality gets easily violated, because the
probability matrices $U$ and $V$ in \eqref{introduction2} are nonnegative and
often dense. This makes spectral methods no longer applicable, and poses
additional difficulties in algorithm design and theoretical analysis.  A
non-convex estimator based on rank-constrained likelihood maximization  is
proposed in \cite{2018arXiv180400795L} where a difference
of convex function programming algorithm is used to handle the rank-constrained
non-smooth optimization problem. Statistical upper bounds and convergence
results are provided.

As is implied by \eqref{introduction2}, our target matrix $P$ has a low-rank and nonnegative factorization. Recovery of low-rank matrices is intensively studied in the past decade, and the matrix completion (MC) is a canonical problem closely related to our work.
MC aims to recover a low-rank matrix based on only a fraction of its entries.
It is often formulated as finding the matrix with smallest rank satisfying a collection of linear constraints. However, such a nonconvex formulation is NP-hard in general, and the convex surrogate using nuclear norm is a remedy. It is proved by \cite{candes2009exact,doi101137070697835} that under mild conditions, the nuclear norm relaxation identifies the exact solution to MC.
As is shown in \cite{Chandrasekaran2012}, the success of nuclear norm can be interpreted from the perspective of atomic norms. It inspires us to propose a convex heuristic for state aggregation using similar techniques introduced in \cite{Chandrasekaran2012}.

Our treatment of a problem with low-rank solution through matrix factorization has similar versions in the context of semidefinite programming \cite{journee2010low, Samuel2001A}, MC \cite{wen2012solving, 2009arXiv09105260K}, nuclear norm minimization \cite{mishra2013low}, tensor factorization \cite{2015arXiv150607540H}, etc.
The scheme presented in our paper is distinct from others since it involves a
nonnegative matrix factorization (NMF). The NP-hardness
\cite{doi101137070709967} and nonconvexity of NMF make it hard for global
optimization, and strategies have been proposed to find local solutions. The
algorithms in \cite{NIPS2000_1861, doi101162neco200719102756} use alternating
(projected) gradient descending minimization with carefully chosen step sizes,
while methods in \cite{kim2008nonnegative, kim2011fast} are developed under the
alternating nonnegative least square (ANLS) framework. 
In our state aggregation problem, we adopt another alternating direction method,
namely the proximal alternating linearized method (PALM), whose convergence to
critical points is proved in \cite{Bolte2014}.

\subsection{Scope of This Paper}

In this work, we develop a  low-nonnegative-rank approximation  method to identify an underlying state
aggregation structure based on trajectories of a Markov chain. After
constructing an empirical transition matrix, we compress the state space via an
optimization problem regularized by the nonnegative rank \cite{COHEN1993149},
which is similar to the usual notion of rank but has nonnegative requirements on
the factorization. Due to the combinatorial nature of nonnegative rank, {the
regularized optimization problem is nonconvex and hard to solve.} Therefore, we propose an atomic regularizer as a convex surrogate. The atomic regularizer is an extension of the atomic norm in \cite{Chandrasekaran2012}. It achieves empirical successes in maintaining low-nonnegative-rank structure, and substantially reduces recovery errors compared to the vanilla empirical transition matrix. We also note that, our atomic regularizer is quite general and not restricted to stochastic matrices. Hence, it can be of independent interest in NMF studies.

For the sake of an efficient algorithm, we reform the convexified problem into a
factorized optimization model in terms of $U,V \in \mathbb{R}_+^{d \times s}$, where $X =
UV^T$ is an NMF of the variable $X \in \mathbb{R}^{d \times d}$ and $s$ refers to the rank that needs to be adapted. Our
factorized problem is more complicated than the usual NMF due to the extra linear constraints
on $U$ and $V$ so that $X=UV^T$ is still a Markov transition matrix.
After investigating the first-order optimality properties, we prove that a local
minimizer of the factorized optimization model provides a solution to the convex
problem if and only if a certificate derived from the
convex model is satisfied.
Based on this observation, we devise a scheme that involves minimizing a sequence of factorized optimization models until the resulting local minimizer represents a global solution.
Two strategies for adjusting the rank $s$ are developed so that local minima are escaped.
One is to append an additional column to $U$ and $V$ simultaneously, which can be interpreted as an approximation of negative subgradient step. 
The other is to reduce the dimension of $U$ and $V$ together by means of linear
combinations so that the rank can be reduced to a proper range. The convergence
of our proposed algorithm is  very likely ensured since these two methods of rank adjustment gurantee monotone
improvement of the objective function.

We illustrate our new approach on synthetic data whose ground-truth solution is known.
In situations without sampling noise, the algorithm stops exactly at the intrinsic dimension of the system.
Experiments are also conducted on simulated trajectories, where we investigate how the ambient, intrinsic dimension and the length of trajectory affect the recovery of $P^*$. 
We finally apply our model to analyze a Manhattan transportation data set. Compared with an existing method in \cite{2017arXiv170507881Y}, our zoning results are more aligned to the geometric location without knowing it in advance.

\subsection{Outline}

The remainder of the paper is organized as follows.
In Section \ref{Section2}, we formulate state aggregation into an optimization problem regularized by nonnegative rank.
In Section \ref{Section3}, we introduce the atomic regularizer as a convex surrogate for the nonnegative rank, and derive necessary and sufficient conditions for the global optimality.
In Section \ref{Section4}, we design {a factorized optimization model} and look at its KKT conditions that govern the local solutions.
In Section \ref{Section5}, we develop strategies to escape from local minima and
establish an adaptive rank algorithm. 
 Numerical results are presented in Section \ref{Section6}.

\subsection{Notation}

For a column vector $\mathbf{x} \in \mathbb{R}^p$, we denote its $i$-th entry as $x_i$, the Euclidean norm as $\|\mathbf{x}\|_2$, and the $\ell_{\infty}$-norm as $\|{\bf x}\|_{\infty} := \max_{1\leq i \leq p}|x_i|$.
For a matrix $A \in \mathbb{R}^{p \times q}$, we denote $(A)_{ij}$ or $a_{ij}$
as the entry in the $i$-th row and $j$-th column of $A$. We denote $A_j \in \mathbb{R}^p$ as the $j$-th column of $A$, and $A^i \in \mathbb{R}^q$ as the transpose of the $i$-th row. The inner product of two matrices $A,B \in \mathbb{R}^{p \times q}$ is represented by $\langle A, B \rangle = \sum_{i,j}a_{ij}b_{ij}$.
If $A \in \mathbb{R}^{p \times q}$ is entrywisely nonnegative, we write that $A \geq 0$ or $A \in \mathbb{R}_+^{p \times q}$.
The nonnegative part of a matrix $A$ is denoted by $[A]_+$, $i.e.$,
$ \left( [A]_+ \right)_{ij} = \max\{a_{ij},0\} $.
In this paper, several matrix norms are involved, including the nuclear norm $\|
A \|_*$ (the sum of singular values), the spectral norm $\| A \|_2$ (the
largest singular value), the Frobenius norm $\| A \|_F =
\sqrt{\sum_{i,j}a_{ij}^2}$, and the $\ell_1$-norm $\| A \|_{\ell_1} = \sum_{i,j}|a_{ij}|$. We denote the $k$-th largest singular value of a matrix $A$ by $\sigma_k(A)$.
Moreover,  ${\bf 1}_p$ and ${\bf 0}_p$ are the $p$-dimensional all-one and
all-zero column vectors, respectively. We define ${\bf e}_i \in \mathbb{R}^p$ as
the vector whose $i$-th entry equals to one and other entries are zeros. ${\bf
0}_{p \times q}$ represents the $p$-by-$q$ zero matrix. For a vector ${\bf x}
\in \mathbb{R}^p$, ${\bf x}^{\bot}$ stands for the orthogonal complement of the
subspace spanned by ${\bf x}$, $i.e.$, ${\bf x}^{\bot} = \left\{ {\bf y} \in
\mathbb{R}^p \mid {\bf y}^T {\bf x} = 0 \right\}$. We denote by ${\bf
diag}\{{\bf x}\}$ or ${\bf diag}\{x_i\}_{i=1}^{p}$ the $p$-by-$p$ diagonal
matrix with ${\bf x} \in \mathbb{R}^p$ on its diagonal. For an index set
$\mathcal{I} = \{ j_1,j_2,\cdots,j_k \} \subseteq \{1,2,\cdots,q\}$, we define
a $k$-by-$k$ diagonal matrix ${\bf diag}\{ x_i \}_{i \in \mathcal{I}} := {\bf
diag}\{ x_{j_i} \}_{i=1}^k$, and $A_{\mathcal{I}} :=
[A_{j_1},A_{j_2},\cdots,A_{j_k}]$, the matrix composed of the
$j_1,j_2,\cdots,j_k$-th columns in $A$. The convex hull of a set $\mathcal{C}
\subseteq \mathcal{R}^p$ is written as ${\rm conv}(\mathcal{C})$.

\section{Problem Setup} \label{Section2}

Consider a discrete-time finite-state Markov chain with unknown transition matrix. Suppose that we only have access to observations of state-to-state trajectories. In this work, we are concerned about computational methods to identify state aggregation structures of the Markov chain. 

\subsection{Markov Chains}

Suppose that a discrete-time Markov chain has a state space $\mathcal{S} = \{ s_1, s_2,\cdots, s_d \}$ and is driven by an unknown probability transition matrix $P^* \in \mathbb{R}^{d \times d}$. Let $(i_0, i_1, i_2,\cdots)$ be a trajectory generated by the Markov chain. By definition,
\[ \mathbb{P}\left[ i_t = s_j \mid i_{t-1} = s_i, i_{t-2}, \cdots, i_0 \right] = p^*_{ij}, \quad t = 1,2,\cdots.  \]

In our problem, we are given an observed trajectory $(i_0, i_1, \cdots, i_n)$ of length ($n+1$). Based on this sample trajectory, we want to have an initial estimate of the transition probabilities.
We assume that the Markov chain is ergodic and has a unique stationary distribution $\xi^* \in \mathbb{R}^d$. One can estimate $\xi^*$ using an empirical stationary distribution $\hat{\xi}^{(n)}$ given by
\begin{equation} \label{FormulateXi}
    \hat{\xi}_j^{(n)} := \frac{1}{n} \sum_{t=1}^n \mathbbm{1}_{\{i_t = s_j\}}, \quad j = 1,2,\cdots,d,
\end{equation}
where $\mathbbm{1}_{\{\cdot\}}$ is an indicator function that takes value $1$ if the subscript event happens or takes $0$ otherwise.
Similarly, an empirical probability transition matrix $\hat{P}^{(n)}$ is formulated as
\begin{equation} \label{FormulateP}
    \hat{p}^{(n)}_{ij} := \left\{ \begin{aligned}
    & \frac{\sum_{t=1}^n \mathbbm{1}_{\{i_{t-1}=s_i,i_t=s_j\}}}{\sum_{t=1}^n\mathbbm{1}_{\{i_{t-1}=s_i\}}}, & \quad & \text{if $\sum_{t=1}^n\mathbbm{1}_{\{i_{t-1}=s_i\}} \geq 1$}, \\
    & \frac{1}{p}, & \quad & \text{if $\sum_{t=1}^n\mathbbm{1}_{\{i_{t-1}=s_i\}} = 0$},
    \end{aligned} \right.
\end{equation}
which is also the maximum likelihood estimator of $P^*$ \cite{1023072237025}.
Asymptotically, the ergodic theorem of Markov chain 
ensures that $\hat{\xi}^{(n)}$ and $\hat{P}^{(n)}$ are strongly consistent. Finite sampling error bounds are also available, see \cite{2018arXiv180202920Z}.
However, due to the curse of dimensionality, a large $d$ makes $\hat{P}^{(n)}$ a poor estimator. We next seek for a method to derive a better estimator from $\hat{P}^{(n)}$.

\subsection{State Aggregation Problem}

In this part, we leverage the state aggregation model \eqref{introduction2} to compress the state space of the Markov chain, following the idea of \cite{2018arXiv180202920Z}. Based on this additional structure, we expect to recover the ground-truth transition matrix $P^*$ with higher accuracy, compared with the initial estimator $\hat{P}^{(n)}$.
A more rigorous definition of the state aggregation structure is shown in Definition \ref{definition}.
\begin{definition}[See \cite{2018arXiv180202920Z}] \label{definition}
    A $d$-state Markov chain generated by a probability transition matrix $P^* \in \mathbb{R}^{d \times d}$ admits a state aggregation structure with intrinsic dimension $r$, if there exist matrices $U^* \in \mathcal{U}^{d \times r} := \left\{ U \in \mathbb{R}_+^{d \times r} \mid U {\bf 1}_r = {\bf 1}_d \right\}$, $V^* \in \mathcal{V}^{d \times r} := \left\{ V \in \mathbb{R}_+^{d \times r} \mid V^T {\bf 1}_d = {\bf 1}_r \right\}$ such that
    \begin{equation}\label{factorization} P^* = U^* \left( V^* \right)^T . \end{equation}
    The entries of $U^*$ and $V^*$ refer to aggregation and disaggregation probabilities, respectively.
\end{definition}

Definition \ref{definition} refines model \eqref{introduction2}.
We notice that in expression \eqref{introduction2}, the design of a latent Markov chain $\tilde{P}$ is flexible. Without loss of generality, one can simply take $\tilde{P}$ to be the identity matrix $I_r$. Therefore, we merge $\tilde{P}$ and $U$ (or $V^T$) together in \eqref{factorization}. It is important to note that $U^*$ and $V^*$ here are still not unique. One can only expect to recover the state aggregation structure up to linear transformation.

According to Definition \ref{definition}, any $d$-by-$d$ stochastic matrix trivially admits a state aggregation structure with intrinsic dimension $d$, whereas in practical applications it is desirable to have $r \ll d$.
The constraints $U^* \in \mathcal{U}^{d \times r}$ and $V^* \in \mathcal{V}^{d \times r}$ in Definition \ref{definition} are proposed so that the aggregation and disaggregation probabilities are meaningful, which leads to an NMF structure in \eqref{factorization}.
In the following, we will denote $\mathcal{U}^{1 \times q}$ by $\mathcal{U}^q$ and $\mathcal{V}^{p \times 1}$ by $\mathcal{V}^p$ for simplicity. If $p=q$, the elements in $\mathcal{U}^{p \times q}$ are stochastic matrices.

In our problem, we assume that the system is implicitly driven by an $r$-state
latent Markov chain as described in section \ref{Section1}
 (see also \cite{2018arXiv180202920Z}, Figure 1), where $r \ll d$ and is unknown. A state aggregation structure \eqref{factorization} with intrinsic dimension $r$ is then taken as an underlying assumption.
For the convenience of discussions, we introduce the notion of {\it nonnegative rank},
\[ {\bf rank}_+(A) := \min\left\{m \mid A = BC^T, B \in \mathbb{R}_+^{d \times m}, C \in \mathbb{R}_+^{d \times m} \right\}, \quad \forall A \in \mathbb{R}^{d \times d}. \]
It can be naturally derived from Definition \ref{definition} that the ground-truth transition matrix $P^*$ in our problem satisfies ${\bf rank}_+(P^*) \leq r$.
We will next propose an optimization problem to conduct state aggregation by controlling the nonnegative rank of the estimator.

Our goal here is to identify a better estimator $X \in \mathcal{U}^{d \times d}$ of $P^*$ based on the empirical $\hat{P}^{(n)}$.
To measure the discrepancy between $X$ and $\hat{P}^{(n)}$, we define
\begin{equation} \label{definition_g}
    g(X) := \frac{1}{2} \left\| \hat{\Xi}(\hat{P}^{(n)} - X) \right\|_F^2,
\end{equation}
where $\hat{\Xi} = {\bf diag}\{\hat{\xi}^{(n)}\}$ is a scaling matrix that assigns more weight to states that have occurred more frequently.
An optimization problem is then formulated as
\begin{equation}\label{NonnegativeRankProblem}
    {\rm minimize}_{X \in \mathbb{R}^{d \times d}} \quad g(X) + \chi_{\mathcal{E}}(X) + \lambda {\bf rank}_+(X),
\end{equation}
where $\chi_{\mathcal{E}}$ is a {\it characteristic function} defined as
\begin{equation} \label{definition_characteristic}
\chi_{\mathcal{E}}(X) := \left\{ \begin{aligned} & 0, & \quad & X \in \mathcal{E}, \\ & +\infty, & \quad & \text{otherwise}, \end{aligned} \right. \qquad \forall X \in \mathbb{R}^{d \times d},
\end{equation}
$\mathcal{E}:=\{ X \in \mathbb{R}^{d \times d} \mid X{\bf 1}_d = {\bf 1}_d \}$ is the set of row-normalized matrices,
and $\lambda$ refers to the regularization parameter. In \eqref{NonnegativeRankProblem}, $g(X)$ indicates the fidelity of data. The implicit constraints $\chi_{\mathcal{E}}(X) < +\infty$ and ${\bf rank}_+(X) < + \infty$ imply that $X{\bf 1}_d = {\bf 1}_d$ and $X \geq 0$, forcing $X$ to be a stochastic matrix. $\lambda {\bf rank}_+(X)$ is regarded as a regularization term that enforces low-nonnegative-rank property. When $\lambda$ gets larger, ${\bf rank}_+(X)$ tends to be smaller, so that the degree of aggregation increases.

Note that, a low-nonnegative-rank solution $\hat{X}$ to problem \eqref{NonnegativeRankProblem} helps recover a state aggregation structure \eqref{factorization}. In fact, if $\hat{X} \in \mathbb{R}^{d \times d}$ satisfies ${\bf rank}_+(\hat{X}) = s \leq d$, then there exists a factorization $\hat{X} = \tilde{U}\tilde{V}^T$ for some $\tilde{U}, \tilde{V} \in \mathbb{R}_+^{d \times s}$. Without loss of generality, we assume that $\tilde{V}$ does not have zero column. Taking $\hat{U} := \tilde{U}{\bf diag}\{ {\bf 1}_d^T\tilde{V} \}$ and $\hat{V} := \tilde{V}\left({\bf diag}\{ {\bf 1}_d^T\tilde{V} \}\right)^{-1}$, we have $\hat{V}^T{\bf 1}_d = {\bf 1}_s$ and $\hat{U}{\bf 1}_s = \tilde{U}(\tilde{V}^T{\bf 1}_d) = \hat{X}{\bf 1}_d = {\bf 1}_d$. Hence, $\hat{U} \in \mathcal{U}^{d \times s}$, $\hat{V} \in \mathcal{V}^{d \times s}$, and the factorization $\hat{X} = \hat{U}\hat{V}^T$ identifies a state aggregation structure of $\hat{X}$.

Due to the combinatorial nature of the nonnegative rank, it is hard to solve \eqref{NonnegativeRankProblem} directly. In the following sections, we will develop computation tools for solving \eqref{NonnegativeRankProblem} with theoretical guarantees.

\section{Convexified Formulation via Atomic Regularizers} \label{Section3}

In this section, we look for a convex surrogate for ${\bf rank}_+$ that can be
optimized efficiently and in practice yields desirable solutions with low
nonnegative rank. An atomic regularizer is proposed, inspired by the
concept of atomic norm relaxation in \cite{Chandrasekaran2012}. First-order optimality conditions are also developed for the convexified problem.

\subsection{Atomic Regularizer}

We propose a convex surrogate function of ${\bf rank}_+$ within the general framework of atomic norm in \cite{Chandrasekaran2012}.
The goal of \cite{Chandrasekaran2012} is to represent a vector ${\bf x} \in \mathbb{R}^p$ with a few elementary building blocks which we call {\it atoms}.
To be specific, let $\mathcal{A}$ be an {\it atomic set}, we want ${\bf x} \in \mathbb{R}^p$ to be formed as
\begin{equation} \label{atom}
    {\bf x} = \sum_{i = 1}^k c_i {\bf a}_i, \qquad {\bf a}_i \in \mathcal{A}, c_i > 0, \quad i = 1,2,\cdots,s.
\end{equation}
In MC, for instance, one can take the atomic set to be the collection of unit-spectral-norm rank-one matrices, which we denote by $\mathcal{A}_*$. The affine rank of a matrix is then interpreted as the smallest number of atoms in $\mathcal{A}_*$ whose conic combination can represent the matrix.

The main idea of \cite{Chandrasekaran2012} is to recover a structure \eqref{atom} with a small $k$ by using a convex program that minimizes the following gauge function
\begin{equation} \label{AtomicNorm}
     \| {\bf x} \|_{\mathcal{A}} := \{ t > 0 \mid {\bf x} \in t \cdot {\rm conv}(\mathcal{A}) \}, \quad \forall {\bf x} \in \mathbb{R}^p.
\end{equation}
If $\mathcal{A}$ is centrally symmetric around the origin, then $\| \cdot \|_{\mathcal{A}}$ is indeed a norm and is called an {\it atomic norm}. In MC, the atomic set $\mathcal{A}_*$ induces the nuclear norm $\| \cdot \|_*$ \cite{Chandrasekaran2012}, which we usually take as a convex surrogate for affine rank. The nuclear norm relaxation consistently yields promising results both theoretically and experimentally \cite{candes2009exact,doi101137070697835}. It motivates us to propose a convex relaxation of ${\bf rank}_+$ in an analogous way.

It is important to note that, nuclear norm is not suitable for our problem, since it does not guarantee nonnegativity while \eqref{factorization} actually requires an NMF.
Therefore, we build the following atomic set using nonnegative atoms in $\mathcal{A}_*$,
\begin{equation} \label{definition_atomicset}\mathcal{A}_+ := \left\{ {\bf u}{\bf v}^T \left| \ {\bf u},{\bf v} \in \mathbb{R}_+^d, \|{\bf u}\|_2 = 1, \|{\bf v}\|_2 = 1 \right. \right\}. \end{equation}
In this way, the state aggregation problem \eqref{NonnegativeRankProblem} can be interpreted as approximating $\hat{P}^{(n)}$ with as few atoms in $\mathcal{A}_+$ as possible.
$\mathcal{A}_+$ is not centrally symmetric. However, we can still derive from $\mathcal{A}_+$ a convex relaxation of ${\bf rank}_+$ that is analogous to \eqref{AtomicNorm},
\begin{equation} \label{definition_atomicnorm} \begin{aligned}
     & \Omega(X) & := & \inf\left\{ t > 0 \mid X \in t \cdot {\rm conv}(\mathcal{A}_+) \right\} \\
     & & = & \inf \left\{ \left. \sum_{j=1}^s \|U_j\|_2\|V_j\|_2 \right| \text{ $X = UV^T$ with $U,V \in \mathbb{R}_+^{d \times s}$} \right\}, \quad \forall X \in \mathbb{R}^{d \times d}.
\end{aligned} \end{equation}
Since $\Omega$ is no longer a norm, we only call it an {\it atomic regularizer}.
For an arbitrary $X \in \mathbb{R}^{d \times d}$, if there exists a factorization $X = UV^T$ with $U,V \in \mathbb{R}_+^{d \times s}$ that achieves the infimum in the definition of $\Omega$, $i.e.$, $\Omega(X) = \sum_{j=1}^s \|U_j\|_2\|V_j\|_2$, then we say that it is an {\it optimal factorization (with respect to $\Omega$)}.

Replacing ${\bf rank}_+(X)$ in \eqref{NonnegativeRankProblem} by the atomic regularizer $\Omega(X)$ yields
\begin{equation}\label{reformM_new}
{\rm minimize}_{X \in \mathbb{R}^{d \times d}} \qquad f_{\lambda}(X) := g(X) + \chi_{\mathcal{E}}(X) + \lambda \Omega(X),
\end{equation}
which is now a convex optimization problem.


\subsection{Global Optimality Conditions} \label{GlobalOptimalityConditions}

Since \eqref{reformM_new} is convex over $\mathbb{R}^{d \times d}$, its global solutions can be characterized by first-order information. The objective function $f_\lambda$ consists of three parts, among which the discrepancy term $g(X)$ is differentiable and the remaining two terms $\chi_{\mathcal{E}}(X)$ and $\Omega(X)$ have subdifferentials with explicit form. Next we derive expressions of $\partial \chi_{\mathcal{E}}(X)$ and $\partial \Omega(X)$, as well as the sufficient and necessary conditions for global optimality.

The following two lemmas are concerned with $\partial \chi_{\mathcal{E}}(X)$ and $\partial \Omega(X)$, whose proofs can be found in Appendix \ref{lemma_subgradient_chi:proof} and \ref{lemma_subgradient_Omega:proof}.
\begin{lemma} \label{lemma_subgradient_chi}
    The subdifferential of characteristic function $\chi_{\mathcal{E}}(X)$ is
    \begin{equation} \label{subgradient_chi}
        \partial \chi_{\mathcal{E}}(X) = \left\{ \mu {\bf 1}_d^T \mid \mu \in \mathbb{R}^d \right\}.
    \end{equation}
\end{lemma}

\begin{lemma} \label{lemma_subgradient_Omega}
    The subdifferential of the atomic regularizer $\Omega(X)$ is
    \begin{equation}\label{subgradient_Omega} \partial \Omega(X) = \left\{ W \left| \ \Omega^{\circ}(W) \leq 1,  \langle W, X \rangle = \Omega(X) \right.  \right\}, \end{equation}
    where $\Omega^{\circ}(W)$ is the support function given by
    \begin{equation}\label{DefinitionOmegaCirc} \Omega^{\circ}(W) := \sup_{Z: \Omega(Z) \leq 1} \langle W, Z \rangle. \end{equation}
\end{lemma}

    Based on Lemma \ref{lemma_subgradient_chi} and \ref{lemma_subgradient_Omega}, we now establish the sufficient and necessary conditions for global optimality in Theorem \ref{GlobalKKT}. The proof can be found in Appendix \ref{GlobalKKT:proof}.

    \begin{theorem}[Sufficient and necessary conditions for global optimality of \eqref{reformM_new}]\label{GlobalKKT} ~\\
    Suppose that $\hat{X} \in \mathbb{R}^{d \times d}$ is factorized as $\hat{X}
    = \hat{U}\hat{V}^T$, where $\hat{U} \in \mathcal{U}^{d \times s}, \hat{V}
    \in \mathcal{V}^{d \times s}$ for some $s$, and $\hat{U}$ does not have zero column. Then
    $\hat{X}$ is globally optimal for problem \eqref{reformM_new} with an optimal factorization $\hat{X} = \hat{U}\hat{V}^T$
    if and only if there exists $\mu \in \mathbb{R}^d$ such that
        \begin{subnumcases}{}
            & ${\bf u}^T\left( \mu{\bf 1}_d^T - \nabla g(\hat{X}) \right){\bf v} \leq \lambda, \quad \forall {\bf u},{\bf v} \in \mathbb{R}^d_+ \ s.t. \ \|{\bf u}\|_2=\|{\bf v}\|_2=1$, \label{Theorem1(1)} \\
            & $\left[\mu {\bf 1}_s^T - \nabla g(\hat{X})\hat{V} \right]_+  = \lambda \hat{U}{\bf diag} \left\{ \frac{\|\hat{V}_j\|_2}{\|\hat{U}_j\|_2} \right\}_{j=1}^s$,\label{Theorem1(2)} \\
            & $\left[ {\bf 1}_d \mu^T \hat{U} - \left(\nabla g(\hat{X})\right)^T \hat{U} \right]_+ = \lambda \hat{V} {\bf diag} \left\{ \frac{\|\hat{U}_j\|_2}{\|\hat{V}_j\|_2} \right\}_{j=1}^s$.\label{Theorem1(3)}
        \end{subnumcases}
    The vector $\mu$ is unique if existing.
    \end{theorem}


    In \eqref{Theorem1(1)}, ${\bf u}^T( \mu{\bf 1}_d^T - \nabla g(\hat{X}) ){\bf v}$ approximates the decrease of $g(X)+\chi_{\mathcal{E}}(X)$ in direction ${\bf u}{\bf v}^T$ at $\hat{X}$. Intuitively, the inequality means that for a small perturbation along ${\bf u}{\bf v}^T$, this decrease is dominated by the increase in $\lambda\Omega(X)$. 
    When $({\bf u},{\bf v})=(\hat{U}_j,\hat{V}_j)$ for $j=1,2,\cdots,s$, \eqref{Theorem1(2)} and \eqref{Theorem1(3)} suggest that the inequality \eqref{Theorem1(1)} holds as equality. 
    Conditions \eqref{Theorem1(1)} - \eqref{Theorem1(3)} provide a certification of global optimality.
    In order to determine whether $\hat{X} = \hat{U}\hat{V}^T$ is a global solution to \eqref{reformM_new}, one can first calculate a vector $\mu$ according to \eqref{Theorem1(2)} and \eqref{Theorem1(3)}, and next verify \eqref{Theorem1(1)} with the resulting $\mu$.

\section{Factorized Optimization Model} \label{Section4}

    In the previous section, the convex formulation \eqref{reformM_new} is
    proposed for estimating the state aggregation structure from empirical data. However, due to
    the implicit form of the atomic regularizer $\Omega(X)$, it is impractical
    to directly minimize the function $f_{\lambda}(X)$ over the matrix $X$.
    {Letting $X = UV^T$, we reformulate \eqref{reformM_new} into a factorized optimization model. Then we investigate its equivalence to the convex model and its KKT conditions.}

    \subsection{Factorized Optimization Model and its Equivalence to Convexified Model}

    Although the nonnegative rank of the optimal solution $\hat{X}$ to \eqref{reformM_new} is usually unknown, we can pick an arbitrary $s$ and have an initial guess that ${\bf rank}_+(\hat{X}) \leq s$. We recast \eqref{reformM_new} by nonnegative matrix factorization $X = UV^T$ where $(U, V) \in \mathbb{R}_+^{d \times s} \times \mathbb{R}_+^{d \times s}$. The search space can be further narrowed down to $\mathcal{U}^{d \times s} \times \mathcal{V}^{d \times s}$, since each $(U,V) \in \mathbb{R}_+^{d \times s} \times \mathbb{R}_+^{d \times s}$ with $UV^T \in \mathcal{U}^{d \times d}$ can be mapped into
    \[ \left( U{\bf diag}\{ V^T{\bf 1}_d \}, V\left({\bf diag}\{ V^T{\bf 1}_d \}\right)^{-1} \right) \in \mathcal{U}^{d \times s} \times \mathcal{V}^{d \times s}\]
    under the assumption that $V$ does not have zero column. This transform does not change the product $X = UV^T$.

    We formulate a factorized optimization model as follows
        \begin{equation} \label{fixrM}
         {\rm minimize}_{U \in \mathcal{U}^{d \times s}, V \in \mathcal{V}^{d \times s}} \quad  F_{\lambda}(U,V) = g(UV^T) + \lambda \sum_{j=1}^s \|U_j\|_2 \|V_j\|_2,
        \end{equation}
    where $s$ refers to the rank of the model and is a parameter to be adjusted.
    Intuitively, the problem \eqref{fixrM} not only optimizes $f_{\lambda}(X)$ but also identifies an optimal factorization of $X$ (with respect to $\Omega$). The equivalence between \eqref{fixrM} and the convexified problem \eqref{reformM_new} is provided in Lemma \ref{LemmaEquivalence}. See Appendix \ref{LemmaEquivalence:proof} for the proof.
    \begin{lemma}\label{LemmaEquivalence}
        When $s$ is sufficiently large, \eqref{reformM_new} and \eqref{fixrM} are equivalent.
    \end{lemma}

    Lemma \ref{LemmaEquivalence} indicates that one can solve \eqref{reformM_new} by optimizing a factorized model \eqref{fixrM} with sufficiently large $s$.
     Two strategies for adjusting the rank $s$ are developed in subsection
     \ref{Criteria}.
 

    \subsection{The KKT Conditions of Factorized Optimization Model}

    Compared to \eqref{reformM_new}, the factorized optimization model \eqref{fixrM}
    searches a space of much smaller dimension and provides a tractable way to evaluate the regularization term. However, it is no longer convex.
    In order to bridge the gap between a local minimum of \eqref{fixrM} and a global solution to \eqref{reformM_new}, we develop KKT conditions for \eqref{fixrM} as a counterpart of Theorem \ref{GlobalKKT}. 
    The results are summerized in Theorem \ref{LocalKKT}, whose proof can be found in Appendix \ref{LocalKKT:proof}.

    \begin{theorem}[KKT conditions of \eqref{fixrM}] \label{LocalKKT}
        Suppose that $\left(\hat{U},\hat{V}\right)$ is a local solution to \eqref{fixrM}, and take $\hat{X} = \hat{U}\hat{V}^T$, $\mathcal{I} = \{j \mid \hat{U}_j \neq 0 \}$. Then, there exists $\mu \in \mathbb{R}^d$ such that
            \begin{equation}\label{Theorem2}
             \left\{ \begin{aligned}
            & \left[\mu {\bf 1}_{|\mathcal{I}|}^T - \nabla g(\hat{X})\hat{V}_{\mathcal{I}} \right]_+ = \lambda \hat{U}_{\mathcal{I}}{\bf diag} \left\{ \frac{\|\hat{V}_j\|_2}{\|\hat{U}_j\|_2} \right\}_{j \in \mathcal{I}}, \\
            & \left[ {\bf 1}_d \mu^T \hat{U}_{\mathcal{I}} - \left(\nabla g(\hat{X})\right)^T \hat{U}_{\mathcal{I}} \right]_+ = \lambda \hat{V}_{\mathcal{I}} {\bf diag} \left\{ \frac{\|\hat{U}_j\|_2}{\|\hat{V}_j\|_2} \right\}_{j\in\mathcal{I}}.
            \end{aligned} \right.
        \end{equation}
        The vector $\mu$ is unique if existing and refers to the Lagrangian multiplier.
        \end{theorem}

     According to \eqref{Theorem2}, the existence of a Lagrangian multiplier $\mu$ can help determine whether $(\hat{U},\hat{V})$ solves \eqref{fixrM} locally.
     In addition, given a local solution $(\hat{U},\hat{V})$ to \eqref{fixrM}, one readily
     notices that it satisfies all but one condition in Theorem \ref{GlobalKKT}.
     Comparison between Theorems \ref{GlobalKKT} and \ref{LocalKKT} therefore provides the following certificate that tells us how to verify the global optimality of $\hat{X} = \hat{U}\hat{V}^T$.

     \begin{theorem}[Global optimality certificate]\label{TheoremLocaltoGlobal}
        Suppose that $(\hat{U},\hat{V})$ is a local solution to \eqref{fixrM}. Then $\hat{X} = \hat{U}\hat{V}^T$ is globally optimal for \eqref{reformM_new} if and only if
        \begin{equation} \label{LocaltoGlobal}
            {\bf u}^T \left[ \mu {\bf 1}_d^T - \nabla g(\hat{X}) \right] {\bf v} \leq \lambda , \qquad \forall {\bf u}, {\bf v} \geq 0, \quad \|{\bf u}\|_2 = \|{\bf v}\|_2 = 1,
        \end{equation}
        where $\mu$ is the Lagrangian multiplier defined in Theorem \ref{LocalKKT}.
     \end{theorem}

\section{An Adaptive Rank Estimation Algorithm}\label{Section5}


In this section, we propose an adaptive method that solves a sequence of \eqref{fixrM} while adjusting values of $s$ until reaching the global optimality condition. Details of the algorithm and convergence properties are provided below.

\subsection{A Subroutine: Proximal Alternating Linearized Minimization}

    We first deal with the factorized optimization model for a fixed $s$. In order to apply PALM \cite{Bolte2014}, we rewrite \eqref{fixrM} into
    \begin{equation} \label{PALMfunction}
        {\rm minimize}_{U,V \in \mathbb{R}^{d \times s}} \quad \tilde{F}_{\lambda}(U,V) := F_{\lambda}(U,V) + \chi_{\mathcal{U}^{d \times s}}(U) + \chi_{\mathcal{V}^{d \times s}}(V),
    \end{equation}
    where $\chi_{\mathcal{U}^{d \times s}}$ and $\chi_{\mathcal{V}^{d \times s}}$ are characteristic functions analogous to $\chi_{\mathcal{E}}$. PALM is a Gauss-Seidel-like method that alternately minimizes the linearized objective function with respect to $U$ and $V$. In the following, we first develop a method to ensure the differentiability of $F_{\lambda}$ in each iteration, and next show that the subproblems of PALM can be solved efficiently. We provide two types of step sizes. One is derived from \cite{Bolte2014}, which guarantees the convergence of PALM in theory. The other combines the Barzilai-Borwein (BB) step sizes with non-monotone line search. We use it in numerical experiments.

    Under the PALM framework, we generate a sequence $\{(U^k,V^k)\}_{k\in\mathbb{N}}$ to solve \eqref{PALMfunction}. In order that PALM is {well-defined} in each iteration, the objective function $F_{\lambda}$ is supposed to be differentiable at each $(U^k,V^k)$. In other words, it is required that $U^k$ does not have zero columns. To this end, before the $(k+1)$-th iteration, we remove the columns in $U^k$ of which the Euclidean norms are smaller than $\varepsilon_0$, and normalize $U^k$ so that each row sums to one. Here, $\varepsilon_0$ is a user-defined constant concerning the computation precision. We choose $\varepsilon_0 = 10^{-14}$ by default. The corresponding columns in $V^k$ are also removed. Due to the constraint $U{\bf 1}_s = {\bf 1}_d$, the resulting matrices are not empty, and we denote them by $\tilde{U}^k$ and $\tilde{V}^k$.

    After preprocessing, we apply the following modified PALM scheme,
    \begin{eqnarray}
    	& \begin{aligned} & U^{k+1}&\in &\arg\min_U \left\{ \left\langle \nabla_U F_{\lambda}(\tilde{U}^k,\tilde{V}^k), U-\tilde{U}^k \right\rangle + \chi_{\mathcal{U}^{d \times s}}(U) + \frac{1}{2 c_k}\| U-\tilde{U}^k \|_F^2 \right\} \\
    	&&=& {\bf proj}_{\mathcal{U}^{d \times s}}\left( \tilde{U}^k - c_k \nabla_U F_{\lambda}(\tilde{U}^k,\tilde{V}^k)\right), \end{aligned} \label{scheme1}\\
    	& \begin{aligned} & V^{k+1}& \in & \arg\min_V \left\{ \left\langle \nabla_V F_{\lambda}(U^{k+1},\tilde{V}^k),  V-\tilde{V}^k \right\rangle + \chi_{\mathcal{V}^{d \times s}}(V) + \frac{1}{2 d_k}\| V-\tilde{V}^k \|_F^2 \right\} \\
    	&&= &{\bf proj}_{\mathcal{V}^{d \times s}}\left( \tilde{V}^k - d_k \nabla_V F_{\lambda}(U^{k+1},\tilde{V}^k)\right),\end{aligned} \label{scheme2}
    \end{eqnarray}
    where $c_k$ and $d_k$ are some suitable step sizes, and
    \[ {\bf proj}_{\mathcal{C}}(A) := \arg\min_{ X \in \mathcal{C}} \| X - A \|_F, \quad \forall A \in \mathbb{R}^{d \times s}, \]
    for $\mathcal{C} = \mathcal{U}^{d \times s}$ or $\mathcal{V}^{d \times s}$.
    One can derive the closed-form solutions of ${\bf proj}_{\mathcal{U}^{d \times s}}$ and ${\bf proj}_{\mathcal{V}^{d \times s}}$ from Lemma \ref{projsol}.
    \begin{lemma}[Theorem 2.2 in \cite{ProjOperator}]\label{projsol}  For a vector ${\bf y} \in \mathbb{R}^p$, let $y^{(1)} \geq y^{(2)} \geq \cdots \geq y^{(p)}$ be a rearrangement of the entries in ${\bf y}$,
    and take \[ l = \max \left\{j \ \left| \ \sum\nolimits_{k=1}^{j-1}(y^{(k)} - y^{(j)})< 1, 1 \leq j \leq p \right.\right\},
        \text{ \ and \ }  \eta = \frac{1}{l}\left(1 - \sum_{k=1}^{l} y^{(k)}\right). \]
    Then, the projection of ${\bf y}$ onto $ \mathcal{V}^p = \left\{ {\bf x} \in \mathbb{R}_+^p \mid {\bf x}^T{\bf 1}_p = 1 \right\} $
    is given by
    \[ {\bf proj}_{\mathcal{V}^p}({\bf y}) = \left[ {\bf y} + \eta {\bf 1}_p \right]_+. \]
    \end{lemma}
    Lemma \ref{projsol} also provides a numerically efficient way to calculate \eqref{scheme1} and \eqref{scheme2}.

    We now consider the step sizes $c_k$ and $d_k$ suggested by \cite{Bolte2014}, which rely on the Lipschitz smoothness of $F_{\lambda}$. The partial gradients of $F_{\lambda}$ are
    \[ \begin{aligned}
     & \nabla_U F_{\lambda}(U,V) = -\hat{\Xi}^2(\hat{P}^{(n)} - UV^T)V + \lambda U {\bf diag}\left\{\frac{\|V_j\|_2}{\|U_j\|_2} \right\}_{j=1}^s, \\
     & \nabla_V F_{\lambda}(U,V) = -(\hat{P}^{(n)}-UV^T)^T\hat{\Xi}^2U + \lambda V {\bf diag}\left\{\frac{\|U_j\|_2}{\|V_j\|_2} \right\}_{j=1}^s,
     \end{aligned}\]
    if existing.
    Within the feasible set, $\nabla_U F_{\lambda}(U,V)$ is Lipschitz continuous with constant
    $L_1(V) = \|\hat{\Xi}\|_F^2\|V^TV\|_F + \lambda \varepsilon_0^{-1}$
    with respect to $U$ for a fixed $V$, since
    { \[ \begin{aligned} & \left\| \nabla_U F(U',V) - \nabla_U F(U,V) \right\|_F \leq   \left\| \hat{\Xi}^2(U' - U)V^TV \right\|_F + \lambda\varepsilon_0^{-1}\|U'-U\|_F \\ \leq & \left(\|\hat{\Xi}\|_F^2\|V^TV\|_F+\lambda\varepsilon_0^{-1}\right)\|U'-U\|_F, \qquad \text{$\forall$ feasible $U'$, $U$}. \end{aligned}\]}
    Similarly, $\nabla_U F(U,V)$ is also Lipschitz continuous with constant $L_2(U) = \|U^T\hat{\Xi}^2U \|_F$ \\ $ + \lambda\sqrt{d}\|U\|_F$ with respect to $V$ for a fixed $U$. By taking
    \begin{equation}\label{cd}
    c_k = \left(\gamma_1 L_1(V^k)\right)^{-1},\quad d_k = \left(\gamma_2 L_2(U^k)\right)^{-1}
    \end{equation}
    with $\gamma_1 > 1$, $\gamma_2 >1$, one can prove that the modified PALM converges theoretically.

    In fact, as \eqref{scheme1} and \eqref{scheme2} iterate, the rank $s$ monotonically decreases and converges to a positive integer. Therefore, when analyzing asymptotic behaviors, we only need to deal with situations where $s$ is fixed and $\|U_j^k\| \geq \varepsilon_0$ for $j=1,2,\cdots,s$. Under these conditions, the results in \cite{Bolte2014} can be extended to our proposed scheme. We summarize the convergence properties in Theorem \ref{GlobalConverge}, and defer the proof to Appendix \ref{GlobalConverge:proof}.
    \begin{theorem}[Corollary of Lemma 3 and Theorem 1 in \cite{Bolte2014}]\label{GlobalConverge}
        Let $\{(U^k,V^k)\}_{k\in \mathbb{N}}$ be a sequence generated by \eqref{scheme1} and \eqref{scheme2} iteratively and $\|U_j^k\|_2 \geq
        \varepsilon_0$ for $j = 1,2,\cdots,s$. 
        \begin{enumerate}
        	\item {\rm (Lemma 3 (i))} The sequence $\{ F_{\lambda}(U^k,V^k) \}_{k \in \mathbb{N}}$ is non-increasing.
            \item {\rm (Theorem 1 (i))}The sequence $\{ (U^k,V^k) \}_{k \in \mathbb{N}} $ has finite length, which means,
                \begin{equation*} 
                \sum_{k=1}^{\infty} \left(\|U^{k+1}-U^k\|_F + \|V^{k+1}-V^k\|_F\right) < \infty. \end{equation*}
            \item {\rm (Theorem 1 (ii))}The sequence $\{ (U^k,V^k) \}_{k \in \mathbb{N}} $ converges to a stationary point of problem \eqref{PALMfunction}.
        \end{enumerate}
    \end{theorem}

Theorem \ref{GlobalConverge} shows that the step sizes in \eqref{cd} guarantee monotone improvement of function value and the global convergence of PALM to a stationary point.

{In terms of convergence rate, we notice that $\tilde{F}_{\lambda}$ is semi-algebraic according to the proof of Theorem \ref{GlobalConverge}. Remark 6(ii) in \cite{Bolte2014} suggests that a semi-algebraic function satisfies the {\L}ojasiewicz inequality for some $\theta \in [0,1)$. In other words, there exist $C_0>0$ and $k_0 \in \mathbb{N}$ such that
\[ \left( \tilde{F}_{\lambda}(U^k,V^k) - \tilde{F}_{\lambda}(\hat{U},\hat{V}) \right)^{\theta} \leq C_0 \inf \left\{ \| G \|_F \left| \ G \in \partial \tilde{F}_{\lambda}(U^k,V^k) \right. \right\}, \quad \forall k \geq k_0, \]
where $(\hat{U},\hat{V})$ is the limit point of $\{(U^k,V^k)\}_{k \in \mathbb{N}}$. If $\theta \leq 1/2$, $\{(U^k,V^k)\}_{k \in \mathbb{N}}$ converges linearly to $(\hat{U},\hat{V})$. Otherwise, we only have theoretical guarantees for a geometric convergence rate.}

Aside from \eqref{cd}, we also adjust BB step sizes to our problem and adopt the non-monotone line search technique to enhance numerical performances. Letting
\[ c_{k+1,1} = \frac{\left|\langle S_U^{k}, Y_U^{k}\rangle\right|}{\| Y_U^k \|_F^2}, \quad
   c_{k+1,2} = \frac{\|S_U^{k}\|_F^2}{\left|\langle S_U^{k}, Y_U^{k} \rangle\right|}, \]
where $S_U^k = U^{k+1} - \tilde{U}^{k}$, $Y_U^k = \nabla_U f(U^{k+1},\tilde{V}^{k}) - \nabla_U f(\tilde{U}^{k},\tilde{V}^{k})$, and
\[ d_{k+1,1} = \frac{\left|\langle S_V^{k}, Y_V^{k}\rangle\right|}{\| Y_V^k \|_F^2}, \quad
   d_{k+1,2} = \frac{\|S_V^{k}\|_F^2}{\left|\langle S_V^{k}, Y_V^{k} \rangle\right|}, \]
where $S_V^k = V^{k+1} - \tilde{V}^{k}$, $Y_V^k = \nabla_V f(U^{k+1},V^{k+1}) - \nabla_V f(U^{k+1},\tilde{V}^{k})$, we take \begin{equation} \label{BB}
    c_k = c_{k,1} \delta^p \text{ or } c_k = c_{k,2} \delta^p, \quad d_k = d_{k,1} \delta^q \text{ or } d_k = d_{k,2} \delta^q
\end{equation}
for some $\delta \in (0,1)$, $p, q \in \mathbb{N}$. Here, $p$ and $q$ are respectively the smallest integers such that
\[ \begin{aligned}
    & f(U^{k+1}, \tilde{V}^k) \leq \frac{1}{5}\sum_{t=k-4}^k f(\tilde{U}^t,\tilde{V}^t) - \eta c_k \| \nabla_U f(\tilde{U}^k,\tilde{V}^k) \|_F^2,  \\
    & f(U^{k+1}, V^{k+1}) \leq \frac{1}{5}\sum_{t=k-4}^k f(U^{t+1},\tilde{V}^t) - \eta d_k \| \nabla_V f(U^{k+1},\tilde{V}^k) \|_F^2.
\end{aligned} \]

In numerical experiments, we say that the subroutine PALM converges locally if the decrease in function value satisfies \begin{equation} \label{ConvergeLocally} \left( \bar{f}^k - f^k \right)/f^k < 10^{-3}, \qquad \bar{f}^k = \frac{1}{30}\sum_{t=k-30}^{k-1}f(U^t,V^t). \end{equation} The BB step sizes in \eqref{BB} appear to be efficient in most of our test problems. For instance, when solving a factorized optimization model with $d = 1000$, $s = 10$ and a randomly generated initial point $(U^0,V^0)$, 
our scheme converges within $200$ steps in most scenarios. See Section \ref{Subsection6_1} and \ref{Subsection6_2} for further details of the numerical settings.

\subsection{Successive Refinements to Escape from Local Solutions} \label{Criteria}
	

    When the subroutine PALM converges to a stationary point, one needs to further refine the solution unless it represents a global minimum of the convexified problem \eqref{fixrM}. 
    We propose numerical criteria for the global optimality and develop two methods to escape from local minima without increasing the function value.

    \subsubsection{Criteria to Determine Global Optimality}


    Recall that Theorem \ref{TheoremLocaltoGlobal} provides a certificate \eqref{LocaltoGlobal} to determine the global optimality of a local solution $(\hat{U},\hat{V})$. We now consider how to verify \eqref{LocaltoGlobal} numerically.

     {\bf Exact stopping rule:}
     We adopt the gradient projection method to solve the following optimization problem
     \begin{equation} \label{OmegaDual} \begin{aligned}
        \sigma := \quad & {\rm maximize}_{{\bf u}, {\bf v}} & \quad & {\bf u}^T \left( \mu {\bf 1}_d^T - \nabla g(\hat{X}) \right) {\bf v} \\
        & s.t. & \quad & \|{\bf u}\|_2 = 1, \|{\bf v}\|_2 = 1, \\
        &      &       & {\bf u} \geq 0, {\bf v} \geq 0.
     \end{aligned} \end{equation}
     The algorithm for \eqref{OmegaDual} usually converges in a few steps, and the iterations can be easily calculated. It is safe to say that compared with the PALM process, the time spent on \eqref{OmegaDual} is negligible. As \eqref{LocaltoGlobal} suggests, the solution $\hat{X} = \hat{U}\hat{V}^T$ is considered to be globally optimal if $\sigma \leq (1+\varepsilon_{\rm Exa})\lambda$, where $\varepsilon_{\rm Exa} > 0$ is a user-defined threshold concerning the precision of the solution.
     
     In some scenarios, it requires a very large rank $s$ to reach a reasonably small $\varepsilon_{\rm Exa}$, and the objective function $F_{\lambda}$ is slightly improved as $s$ grows. It is a better choice to stop at a smaller $s$ so that we can compress the state space in a more compact way, even if $\hat{X}=\hat{U}\hat{V}^T$ does not solve \eqref{reformM_new} precisely. To this end, we propose an early stopping rule as follows.
     
     \vspace{0.1cm}

     {\bf Early stopping rule:}
     We define a function
     \[\varphi({\bf v}) = \left\| \left[\left( \mu {\bf 1}_d^T - \nabla g(\hat{X}) \right){\bf v}\right]_+ \right\|_2\] over the set $\Pi := \{ {\bf v} \in \mathbb{R}_+^{d} \mid \|{\bf v}\|_2 = 1 \}$. The condition \eqref{LocaltoGlobal} means that the superlevel set $L_{\lambda} = \{ {\bf v} \in \Pi \mid \varphi({\bf v}) > \lambda \}$ is empty.
     Assume that $(\hat{U},\hat{V})$ is a limit point generated by the modified PALM. If $\hat{U}\hat{V}^T$ is close to a global minimum of \eqref{reformM_new}, $L_{\lambda}$ is supposed to have small measure due to the continuity of $\varphi$. We represent the measure of $L_{\lambda}$ by probabilities, and make an estimate using a group of $i.i.d.$ test vectors $\{\bar{\bf v}_k\}_{k = 1}^N$ uniformly distributed on $\Pi$. In our experiments, $N = 5000$. If each $\bar{\bf v}_k$ satisfies $\varphi(\bar{\bf v}_k) \leq \lambda$, we say that $(\hat{U},\hat{V})$ identifies a global solution.
     This criterion is a relaxation of condition \eqref{LocaltoGlobal} and often results in an early stop in the scheme. If starting from an initial point $(U^0,V^0)$ with a small $s$, we can usually obtain a low-rank solution. Since $\hat{X} = \hat{U}\hat{V}^T$ is only required to lie in the neighborhood of the genuine minimum, this approximate solution is expected to be more robust to sampling noises and work better in real-world problems.

     These two stopping rules are applicable to different situations. When investigating the properties of the state aggregation model \eqref{reformM_new}, we prefer the exact stopping rule so that the problem can be solved more accurately. When dealing with real-world applications, the early stopping rule helps identify approximate solutions with much lower nonnegative rank, and compresses the state space more efficiently.


       \subsubsection{Appending a New Column}

       When condition \eqref{LocaltoGlobal} fails to hold, we need to refine the local solution $(\hat{U},\hat{V})$ so that PALM can resume from a new initial point.
       Either by the exact or early stopping rule, one can find a pair of vectors $\bar{\bf u}, \bar{\bf v} \in \mathbb{R}^d$ such that
       \begin{equation}\label{LocaltoGlobal_Counter} \bar{\bf u}^T \left( \mu {\bf 1}_d^T - \nabla g(\hat{X}) \right) \bar{\bf v} > \lambda, \qquad \bar{\bf u}, \bar{\bf v} \geq 0, \quad \|\bar{\bf u}\|_2 = \|\bar{\bf v}\|_2 = 1. \end{equation}
       Intuitively, $\bar{\bf u}\bar{\bf v}^T$ approximates 
       the negative subgradient directions of $f_{\lambda}$ at $\hat{X} = \hat{U}\hat{V}^T$.
       In fact, any subgradient $Z \in \partial f_{\lambda}(\hat{X})$ can be expressed as $Z = \nabla g(\hat{X}) - \mu {\bf 1}_d^T + \lambda W$ for some $W \in \partial \Omega(\hat{X})$. According to Lemma \ref{subgradient_Omega}, $\Omega^{\circ}(W) \leq 1$, and it implies that
       \[ \left\langle Z , \bar{\bf u}\bar{\bf v}^T \right\rangle = \lambda \bar{\bf u}^T W \bar{\bf v} - \bar{\bf u}^T \left( \mu {\bf 1}_d^T - \nabla g(\hat{X}) \right) \bar{\bf v} \leq \lambda - \bar{\bf u}^T \left( \mu {\bf 1}_d^T - \nabla g(\hat{X}) \right) \bar{\bf v} < 0. \]
       Inspired by this, we consider appending scaled $\bar{\bf u}$ and $\bar{\bf v}$ as new columns to $\hat{U}$ and $\hat{V}$ respectively, and propose \eqref{Append} to construct a better solution.
       See Theorem \ref{NewColumn} for details, whose proof can be found in Appendix \ref{NewColumn:proof}.

    \begin{theorem}[Escaping local minima]\label{NewColumn}
        Suppose that $(\hat{U},\hat{V})$ is a local solution to $\eqref{fixrM}$ where $\hat{U}$ does not have zero column, and $\hat{X}=\hat{U}\hat{V}^T$ is not globally optimal for \eqref{reformM_new}. According to Theorem \ref{TheoremLocaltoGlobal}, there exist $\bar{\bf u},\bar{\bf v} \in \mathbb{R}^d$ satisfying \eqref{LocaltoGlobal_Counter}.
        Take
        \begin{equation}\label{Append}
        \bar{U} := \left[ \begin{array}{cc} {\bf diag}\left\{ {\bf 1}_d - \kappa \bar{\bf u} \right\} \hat{U} & \kappa \bar{\bf u} \end{array} \right], \quad
           \bar{V} := \left[ \begin{array}{cc} \hat{V} & \left(\bar{\bf v}^T {\bf 1}_d\right)^{-1}\bar{\bf v} \end{array} \right] \end{equation}
        for some sufficiently small $\kappa>0$. Then, $(\bar{U}, \bar{V})$ is feasible for problem \eqref{fixrM} and
        \[ {F}_{\lambda}(\bar{U},\bar{V}) < {F}_{\lambda}(\hat{U},\hat{V}). \]
    \end{theorem}
        Theorem \ref{NewColumn} suggests that as long as $(\hat{U},\hat{V})$ is not globally optimal, one can always reduce the global objective value of \eqref{reformM_new} by appending columns. This will guarantee monotone improvement of the solutions.
        
        In numerical experiments, we adopt backtracking line search to determine the appropriate $\kappa$ in \eqref{Append}. We take $\kappa = 0.5^p \|\bar{\bf u}\|_{\infty}^{-1}$, where $p\in \mathbb{N}$ is the smallest integer such that $F_{\lambda}(\bar{U},\bar{V}) < (1-10^{-5})F_{\lambda}(\hat{U},\hat{V})$. If $\kappa < 10^{-8}$, we say that the global optimality is achieved and terminate the algorithm.

        \subsubsection{Removing Redundant Dimensions}

        Suppose that $(\hat{U},\hat{V}) \in \mathcal{U}^{d \times s} \times
        \mathcal{V}^{d \times s}$ is a local solution to \eqref{fixrM}. In the case where $\{ \hat{U}_j\hat{V}_j^T \}_{j=1}^s$ are linearly dependent, the linear combination $\hat{X} = \hat{U}\hat{V}^T = \sum_{j=1}^s \hat{U}_j\hat{V}_j^T$ can be expressed using less than $s$ nonnegative rank-one matrices. In other words, $\hat{X}$ actually admits a state aggregation structure with less than $s$ meta-states.
        To this end, we hope to identify a factorization $\hat{X} =\tilde{U}\tilde{V}^T$ with rank smaller than $s$, and $(\tilde{U},\tilde{V})$ preserves the objective value.

        Recall that $F_{\lambda}(U,V) = g(UV^T) + \lambda \sum_j \| U_j \|_2 \| V_j \|_2$. Since the product $\tilde{U}\tilde{V}^T=\hat{U}\hat{V}^T$, the first term $g(UV^T)$ is unchanged after the adjustment. As for the second term $\sum_{j=1}^{s}\| U_j \|_2 \| V_j \|_2$, the following corollary of Theorem \ref{LocalKKT} suggests that it remains the same under some linear combinations of the rank-one matrices. See Appendix \ref{corollary:proof} for the proof.
        \begin{corollary}[of Theorem \ref{LocalKKT}] \label{corollary}
            Suppose that $(\hat{U},\hat{V}) \in \mathcal{U}^{d \times s} \times
            \mathcal{V}^{d \times s}$ is a local solution to \eqref{fixrM} and
            there exists $\alpha \in \mathbb{R}^s$ such that $\sum_{j=1}^s
            \alpha_j \hat{U}_j\hat{V}_j^T = {\bf 0}_{d \times d}$. Then,
            \[ \sum_{j=1}^s \alpha_j \| \hat{U}_j \|_2 \| \hat{V}_j \|_2 = 0. \]
        \end{corollary}

       Corollary \ref{corollary} inspires us to get rid of the ``redundant dimensions" in $(\hat{U},\hat{V})$ by deleting carefully selected columns and rescaling the rest.
       To be specific, suppose that the linear equation $ \sum_{j=1}^s \alpha_j \hat{U}_j\hat{V}_j^T = {\bf 0}_{d \times d} $ has $s'$ linearly independent solutions $\alpha^1, \alpha^2, \cdots, \alpha^{s'} \in \mathbb{R}^s$ where $1 \leq s' \leq s-1$. One can then construct a vector $\theta \in \mathbb{R}^{s'}$ such that $\sum_{k=1}^{s'} \theta_k \alpha^k \leq {\bf 1}_s$ and the equality holds on at least $s'$ positions. It can be done, for instance, by solving a corresponding linear program with the simplex method.
       By taking \begin{equation}\label{Compress} \tilde{U} = \hat{U} {\bf diag}\left\{ {\bf 1}_s - \sum_{k=1}^{s'}\theta_k\alpha^k \right\}, \quad \tilde{V} = \hat{V}, \end{equation}
       and removing all the zero columns in $\tilde{U}$ and the corresponding columns in $\tilde{V}$, both $\tilde{U}$ and $\tilde{V}$ have at most $s-s'$ columns, which is smaller than $s$. $(\tilde{U},\tilde{V})$ is guaranteed to be feasible to \eqref{fixrM} and Corollary \ref{corollary} implies that $F_{\lambda}(\tilde{U},\tilde{V}) = F_{\lambda}(\hat{U},\hat{V})$.
       This modification may result in a non-stationary point, from where PALM resumes for a smaller $s$.
       
       By removing redundant dimensions in local solutions, we have $s \leq d^2+1$ throughout the computing process. It also guarantees that the algorithm will terminate at a rank no larger than $d^2$.

       In numerical experiments, if 
       \begin{equation}\label{Delta}
            \Delta = \sum_{j=1}^s\sum_{k=1}^{s'} \theta_k\alpha^k_j \hat{U}_j \hat{V}_j^T
       \end{equation}
       satisfies
       $ \left\| \hat{\Xi}\Delta \right\|_F \leq \varepsilon $,
       we consider $\Delta$ as a negligible component of $\hat{X}$ and apply \eqref{Compress} to reduce the rank of factorization. Here, $\varepsilon$ refers to the linear dependence threshold that will be determined later.

\subsection{An Adaptive Rank Algorithm}

        Based on the preliminaries, we now propose an adaptive rank algorithm for the convexified problem \eqref{reformM_new}.
        We apply PALM as a subroutine to solve the factorized optimization model \eqref{fixrM} for fixed $s$ until converging to a stationary point. If the local solution has redundant dimensions as described in Corollary \ref{corollary}, we compress it by \eqref{Compress} and continue PALM. Otherwise, we adopt either exact or early stopping criterion to check condition \eqref{LocaltoGlobal}. If there exists $(\bar{\bf u},\bar{\bf v})$ violating \eqref{LocaltoGlobal}, we append columns to the local solution according to Theorem \ref{NewColumn}, and resume PALM with a new initial point. The full algorithm is summarized as Algorithm \ref{MetaAlg}.

        \begin{algorithm}
        \caption{An Adaptive Rank Algorithm for \eqref{reformM_new}}\label{MetaAlg}
        \KwIn{\begin{minipage}[t]{0.8\linewidth} Initial point $(U^0,V^0) \in \mathcal{U}^{d \times {s_0}} \times \mathcal{V}^{d \times {s_0}}$. \\ Termination criterion: exact or early stopping rule.\end{minipage}}
        \While{``the termination criterion" is not met}{
            \While{the local convergence criterion \eqref{ConvergeLocally} is not met}{
                For each $j=1,2,\cdots, s$ such that $\|U_j^k\|_2 \leq \varepsilon_0$, remove the $j$-th column from both $U^k$ and $V^k$. Update $s$. \\
                Apply \eqref{scheme1} and \eqref{scheme2} to solve \eqref{fixrM} and obtain $(U^{k+1},V^{k+1})$. $k \leftarrow k+1$.
            }
            \If{$\Delta$ in \eqref{Delta} satisfies $\|\hat{\Xi}\Delta\|_F \leq \varepsilon$}{
                Compress $(U^k,V^k)$ according to \eqref{Compress}. Update $s$. }
            \Else{
                Check the global optimality of $(U^k,V^k)$ according to the termination criterion. \\
                    \If{``the termination criterion" is not met}{
                      Determine a step size $\kappa>0$ in \eqref{Append} with backtracking line search. \\ \If{$\kappa\geq10^{-8}$}{Append a column to $U^k$ and $V^k$ according to \eqref{Append}. $s \leftarrow s+1$.}
                     \Else{
                        Terminate the algorithm.
                        }
                    }
            }
        }
        \end{algorithm}

%

\section{Numerical Experiments} \label{Section6}

We present some numerical results to illustrate the efficiency of our proposed state aggregation method. 
When applied to synthetic transition matrices, our scheme stops exactly at the ground-truth intrinsic dimension $r$ when the regularization constant $\lambda$ is appropriately chosen. We also generate Markovian trajectories with varying dimension $d$, $r$ and sampling size $n$, and investigate how these parameters influence the recovery error. Finally, we use our approach to analyze a real dataset of Manhattan transportation network, and conduct extensive comparison with an existing method in \cite{2017arXiv170507881Y}. 
The algorithm is implemented in MATLAB and all experiments are performed on a computer with an Intel i5 1.9GHz and 8GB of RAM.

\subsection{Evaluations of Solution Quality}
We adopt several metrics to evaluate the quality of our solutions.
The following \textit{local errors} $relLE1$ and $relLE2$ concern whether the factorized optimization problem \eqref{fixrM} is solved properly. According to the KKT conditions in Theorem \ref{LocalKKT}, given that the local scheme converges to $(\hat{U},\hat{V})$, we define
{ \[
    relLE1 := \frac{\left\| \left[ \mu {\bf 1}_s^T - \nabla g(\hat{X})\hat{V} \right]_+ - \lambda \hat{U} {\bf diag} \left\{ \frac{\|\hat{V}_j\|_2}{\|\hat{U}_j\|_2} \right\}_{j=1}^s \right\|_{\ell_1}}{\left\| \lambda \hat{U} {\bf diag} \left\{ \frac{\|\hat{V}_j\|_2}{\|\hat{U}_j\|_2} \right\}_{j=1}^s \right\|_{\ell_1}} \]}
and
{ \[
    relLE2 := \frac{\left\| \left[ {\bf 1}_d \mu^T \hat{U} - \left( \nabla g (\hat{X}) \right)^T \hat{U} \right]_+ - \lambda \hat{V} {\bf diag}\left\{ \frac{\|\hat{U}_j\|_2}{\|\hat{V}_j\|_2} \right\}_{j=1}^s \right\|_{\ell_1}}{\left\| \lambda \hat{V} {\bf diag}\left\{ \frac{\|\hat{U}_j\|_2}{\|\hat{V}_j\|_2} \right\}_{j=1}^s \right\|_{\ell_1}}
    \] }
to represent the precision of $(\hat{U},\hat{V})$. Here, the Lagrangian multiplier $\mu \in \mathbb{R}^d$ is estimated by
\[ \mu_i = \frac{\sum_{j=1}^s \left( \lambda \hat{u}_{ij} \frac{\|\hat{V}_j\|_2}{\|\hat{U}_j\|_2} + \left( \nabla g(\hat{X})\hat{V} \right)_{ij} \right) \mathbbm{1}_{\{\hat{u}_{ij} \neq 0\}} }{ \sum_{j=1}^s \mathbbm{1}_{ \{ \hat{u}_{ij} \neq 0 \} } }, \qquad i = 1,2,\cdots,d. \]

As for the global optimality condition \eqref{LocaltoGlobal}, the \textit{global error} is defined as
\[ GE := \Omega^{\circ} (\hat{W}) - 1, \qquad \hat{W} = \lambda^{-1} ( \mu {\bf 1}_d^T - \nabla g(\hat{X}) ), \]
where $ \Omega^{\circ}(\hat{W}) $ is calculated by the gradient projection method.

Aside from $GE$, the duality gap of problem \eqref{reformM_new} serves as another technique to certify global optimality. The Lagrangian dual formulation of \eqref{reformM_new} is given by
\begin{equation}\label{LagrangianDual}
    \begin{aligned}
         \max_{M \in \mathbb{R}^{d \times d}} & \quad  -g^*(M) \\
         s.t. &  \quad \Omega^{\circ}(-M) \leq \lambda,
    \end{aligned}
\end{equation}
where $g^*(M) = \sup_{X \in \mathcal{E}} \left\{ \langle M, X \rangle - g(X) \right\}$. After a routine calculation, we have $g^*(M) = \frac{1}{2} \| \hat{\Xi}^{-1}M + \hat{\Xi}\hat{P}^{(n)} \|_F^2 - \frac{1}{2}\| \hat{\Xi} \hat{P}^{(n)} \|_F^2 - \frac{1}{2d} \| \hat{\Xi}^{-1} M {\bf 1}_d \|_2^2 $. The duality gap can be estimated by $ g(\hat{X}) + \lambda \sum_{j=1}^s\|\hat{U}_j\|_2\|\hat{V}_j\|_2 + g^*(M) $ for some dual candidate $M$ such that $\Omega^{\circ}(-M) \leq \lambda$.
It is a good choice to take $M = -\frac{\lambda}{\Omega^{\circ}(\hat{W})}\hat{W}$. We therefore define the \textit{relative duality gap},
\[ relDG := \frac{  g(\hat{X}) + \lambda \sum_{j=1}^s\|\hat{U}_j\|_2\|\hat{V}_j\|_2 + g^*(M)  }{g(\hat{X}) + \lambda \sum_{j=1}^s\|\hat{U}_j\|_2\|\hat{V}_j\|_2}. \]

In situations where the ground truth is known in advance, we can also calculate the \textit{relative recovery error} of a solution $\hat{X}$,
\[ relRE : = \frac{1}{2} \left\| \hat{\Xi} (\hat{X} - P^*) \right\|_F^2 \left/ \frac{1}{2} \left\| \hat{\Xi} P^* \right\|_F^2. \right. \]
After comparing $relRE$ with the following \textit{relative sampling error}
\[ relSE : = \frac{1}{2} \left\| \hat{\Xi} (\hat{P}^{(n)} - P^*) \right\|_F^2 \left/ \frac{1}{2} \left\| \hat{\Xi} P^* \right\|_F^2, \right. \]
we can tell if solving problem \eqref{reformM_new} is efficient at revealing the system dynamics.

\subsection{Experiments with Exact Low-Nonnegative-Rank Matrices} \label{Subsection6_1}

In this part, we apply Algorithm \ref{MetaAlg} to problem \eqref{reformM_new} where $\hat{P}^{(n)}$ and $\hat{\xi}^{(n)}$ are replaced by $P^*$ and $\xi^*$. We want to check if the algorithm can successfully recover the underlying state aggregation structure when the transition matrix is already seen.

We carry out experiments with $d = 1000,2000,5000$, and set the inner dimension $r = 5$. A test Markov chain is created randomly by the following procedure. We first generate two random matrices $U^*,V^* \in \mathbb{R}^{d \times r}$. The rows of $U^*$ and columns of $V^*$ are independent and uniformly distributed on simplexes $\mathcal{U}^r$ and $\mathcal{V}^d$, respectively. We assemble a Markov transition matrix $P^* = U^*(V^*)^T$ and calculate its stationary distribution $\xi^*$.

When implementing Algorithm \ref{MetaAlg}, we adopt the exact stopping rule for higher precision. Since the algorithm appends columns one by one, whereas reduces redundant dimensions all at once, it is a better choice to start with a large initial dimension $s_0$ and let $s$ reduce to a proper range after some computations. To this end, we choose $s_0 = 300$ which is large enough so that $s_0$ does not affect the terminal dimension $\hat{s}$.

The numerical results are shown in Table \ref{exact1000}, \ref{exact2000} and
\ref{exact5000}. We learn that, if the regularization parameter $\lambda$ is sufficiently small, Algorithm \ref{MetaAlg} converges to a solution $(\hat{U},\hat{V}) \in \mathbb{R}^{d \times \hat{s}} \times \mathbb{R}^{d \times \hat{s}}$ with $\hat{s}=r$, the nonnegative rank of $P^*$. It implies that, in the convexified problem \eqref{reformM_new}, the atomic regularizer performs well in identifying low-nonnegative-rank structures. We also notice that, when $\lambda = 10^{-9}$, $relRE$ are small ($<(5\%)^2$), and the solutions to \eqref{reformM_new} are close to the ground-truth transition matrices. It is safe to say that our proposed method yields reliable results in recovering the state aggregation structures of the ground-truth $P^*$.

\begin{table}
    \centering
    {
    \begin{tabular}{c|cccc|cccc}
        \hline
        $\lambda$ & $\hat{s}$ & obj value & $relRE$ & {\footnotesize time (${\rm s}$)} & $relLE1$ & $relLE2$ & $GE$ & $relDG$ \\
        \hline
        $1$e${\text{-}06}$ & $1$ & $1.14$e${\text{-}06}$ & $1.31$e${\text{-}01}$ & $13.8$ & $0$ & $6.7$e${\text{-}16}$ & $\text{-}1.0$e${\text{-}15}$ & $2.8$e${\text{-}15}$ \\
        $1$e${\text{-}07}$ & $5$ & $1.60$e${\text{-}07}$ & $2.15$e${\text{-}02}$ & $38.2$ & $2.7$e${\text{-}03}$ & $5.6$e${\text{-}03}$ & $1.6$e${\text{-}03}$ & $5.8$e${\text{-}04}$ \\
        $1$e${\text{-}08}$ & $5$ & $1.80$e${\text{-}08}$ & $3.91$e${\text{-}04}$ & $77.2$ & $2.3$e${\text{-}03}$ & $2.4$e${\text{-}03}$ & $1.7$e${\text{-}03}$ & $2.6$e${\text{-}03}$ \\
        $1$e${\text{-}09}$ & $5$ & $1.84$e${\text{-}09}$ & $6.76$e${\text{-}06}$ & $99.4$ & $5.8$e${\text{-}02}$ & $3.5$e${\text{-}02}$ & $3.9$e${\text{-}02}$ & $4.5$e${\text{-}02}$\\
        \hline
    \end{tabular}
    }
    \caption{The recovery of an exact low-nonnegative-rank matrix $P^*$ with $d = 1000$, $r = 5$, $\sigma_1(\Xi^*P^*) = 1.20 \times 10^{-3}$, $\sigma_r(\Xi^*P^*) = 1.93 \times 10^{-4}$. Algorithm \ref{MetaAlg} starts with $s_0 = 300$. The linear dependence threshold $\varepsilon = 5 \times 10^{-5}$.} \label{exact1000}
    {
    \begin{tabular}{c|cccc|cccc}
        \hline
        $\lambda$ & $\hat{s}$ & obj value & $relRE$ & {\footnotesize time (${\rm s}$)} & $relLE1$ & $relLE2$ & $GE$ & $relDG$ \\
        \hline
        $1$e${\text{-}06}$ & $1$ & $1.04$e${\text{-}06}$ & $1.86$e${\text{-}01}$ & $36.8$ & $2.0$e${\text{-}18}$ & $7.6$e${\text{-}16}$ & $\text{-}1.0$e${\text{-}15}$ & $2.4$e${\text{-}15}$ \\
        $1$e${\text{-}07}$ & $1$ & $1.27$e${\text{-}07}$ & $1.08$e${\text{-}01}$ & $37.7$ & $7.4$e${\text{-}18}$ & $3.5$e${\text{-}06}$ & $2.8$e${\text{-}11}$ & $7.5$e${\text{-}09}$ \\
        $1$e${\text{-}08}$ & $5$ & $1.70$e${\text{-}08}$ & $4.65$e${\text{-}03}$ & $78.7$ & $1.4$e${\text{-}03}$ & $2.3$e${\text{-}03}$ & $1.8$e${\text{-}03}$ & $2.7$e${\text{-}03}$ \\
        $1$e${\text{-}09}$ & $5$ & $1.80$e${\text{-}09}$ & $7.35$e${\text{-}05}$ & $140.7$ & $1.7$e${\text{-}02}$ & $2.8$e${\text{-}02}$ & $3.8$e${\text{-}02}$ & $4.7$e${\text{-}02}$ \\
        \hline
    \end{tabular}
    }
    \caption{The recovery of an exact low-nonnegative-rank matrix $P^*$ with $d = 2000$, $r = 5$, $\sigma_1(\Xi^*P^*) = 5.97 \times 10^{-4}$, $\sigma_r(\Xi^*P^*) = 9.81 \times 10^{-5}$. Algorithm \ref{MetaAlg} starts with $s_0 = 300$. The linear dependence threshold $\varepsilon = 5 \times 10^{-5}$.} \label{exact2000}
    {
    \begin{tabular}{c|cccc|cccc}
        \hline
        $\lambda$ & $\hat{s}$ & obj value & $relRE$ & {\footnotesize time (${\rm s}$)} & $relLE1$ & $relLE2$ & $GE$ & $relDG$ \\
        \hline
        $1$e${\text{-}06}$ & $1$ & $1.01$e${\text{-}06}$ & $2.33$e${\text{-}01}$ & $75.0$ & $0$ & $5.0$e${\text{-}14}$ & $\text{-}5.0$e${\text{-}14}$ & $1.6$e${\text{-}14}$ \\
        $1$e${\text{-}07}$ & $1$ & $1.06$e${\text{-}07}$ & $1.67$e${\text{-}01}$ & $90.1$ & $3.0$e${\text{-}18}$ & $3.4$e${\text{-}14}$ & $\text{-}5.9$e${\text{-}14}$ & $1.1$e${\text{-}14}$\\
        $1$e${\text{-}08}$ & $5$ & $1.39$e${\text{-}08}$ & $8.61$e${\text{-}02}$ & $121.8$ & $3.7$e${\text{-}04}$ & $1.0$e${\text{-}03}$ & $1.5$e${\text{-}03}$ & $1.2$e${\text{-}03}$ \\
        $1$e${\text{-}09}$ & $5$ & $1.75$e${\text{-}09}$ & $2.11$e${\text{-}03}$ & $181.6$ & $1.3$e${\text{-}02}$ & $2.0$e${\text{-}02}$ & $1.3$e${\text{-}02}$ & $5.3$e${\text{-}03}$ \\
        \hline
    \end{tabular}
    }
    \caption{The recovery of an exact low-nonnegative-rank matrix $P^*$ with $d = 5000$, $r = 5$, $\sigma_1(\Xi^*P^*) = 2.38 \times 10^{-4}$, $\sigma_r(\Xi^*P^*) = 3.89 \times 10^{-5}$. Algorithm \ref{MetaAlg} starts with $s_0 = 300$. The linear dependence threshold $\varepsilon = 5 \times 10^{-5}$.} \label{exact5000}
\end{table}

\subsection{Experiments with Simulated Data} \label{Subsection6_2}

In this part, we investigate statistical properties of solutions of Algorithm \ref{MetaAlg} when the input is simulated transition trajectories.
Given a fixed trajectory, we solve problem \eqref{reformM_new} with different $\lambda$, and obtain a series of $\hat{X}$ and $relRE$. We plot $relRE$ against $\lambda$ and identify the most appropriate $\lambda^*$ that yields the smallest error $relRE^*$.
In order to generate a path of solutions corresponding to different values of $\lambda$, we employ the \textit{warm-restart} approach. To be specific, we initialize the algorithm with an obtained solution to a problem with slightly larger $\lambda$. 
In each group of our experiments, we investigate how $relRE^*$ is influenced by one parameter among $d$, $r/d$ and $n/d^2$. We let the parameter of interest take different values and fix the other two.

We provide some benchmark simulations for the state aggregation problem. A Markov chain $P^*$ with $d$ states and $r$ inherent meta-states is created in the same way as in Subsection \ref{Subsection6_1}. After choosing an initial state $i_0$ under invariant distribution $\xi^*$, we randomly generate $i_t, t = 1,\cdots,n$, step by step, and form a trajectory $(i_0, i_1, \cdots, i_n)$ of length $n+1$.

For the sake of higher precision, we adopt the exact stopping rule when applying Algorithm \ref{MetaAlg}. For the same reasons mentioned before, we also take the initial dimension $s_0 = 300$. $relLE1$, $relLE2$, $GE$ and $relDG$ are used to measure the optimization errors of Algorithm \ref{MetaAlg}. If these quantities are relatively small, we say that the optimization errors are dominated by the statistical error $relRE$. In this case, it is secure to study the statistical properties of model \eqref{reformM_new} with the computational results. Also, for each group of parameters $d$, $r/d$ and $n/d^2$, we run Algorithm \ref{MetaAlg} on $5$ independent trajectories so as to reduce the random errors.

\vspace{0.3cm}

\noindent {\bf Dimension $d$}

\vspace{0.1cm}

\begin{figure}
    \centering
    \begin{minipage}{0.495\linewidth}
        \includegraphics[width = \textwidth]{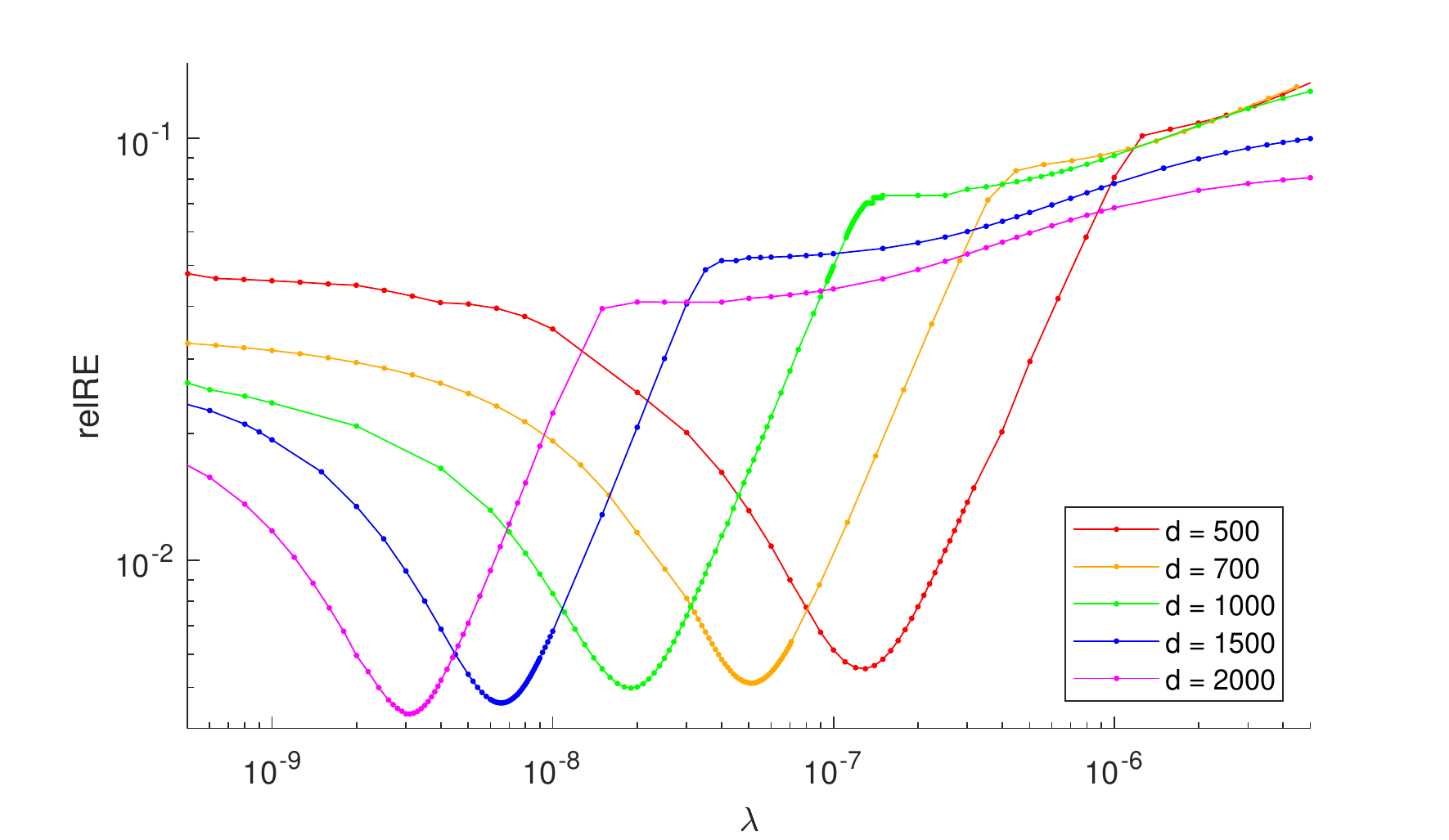}
        \caption{$relRE$-$\lambda$ curves with $r/d = 0.01$, $n/d^2 = 10$, $d = 500$, $700$, $1000$, $1500$, or $2000$. \protect\\ ~ \protect\\  }
        \label{figure_d_1}
    \end{minipage}
    \begin{minipage}{0.495\linewidth}
        \includegraphics[width = \textwidth]{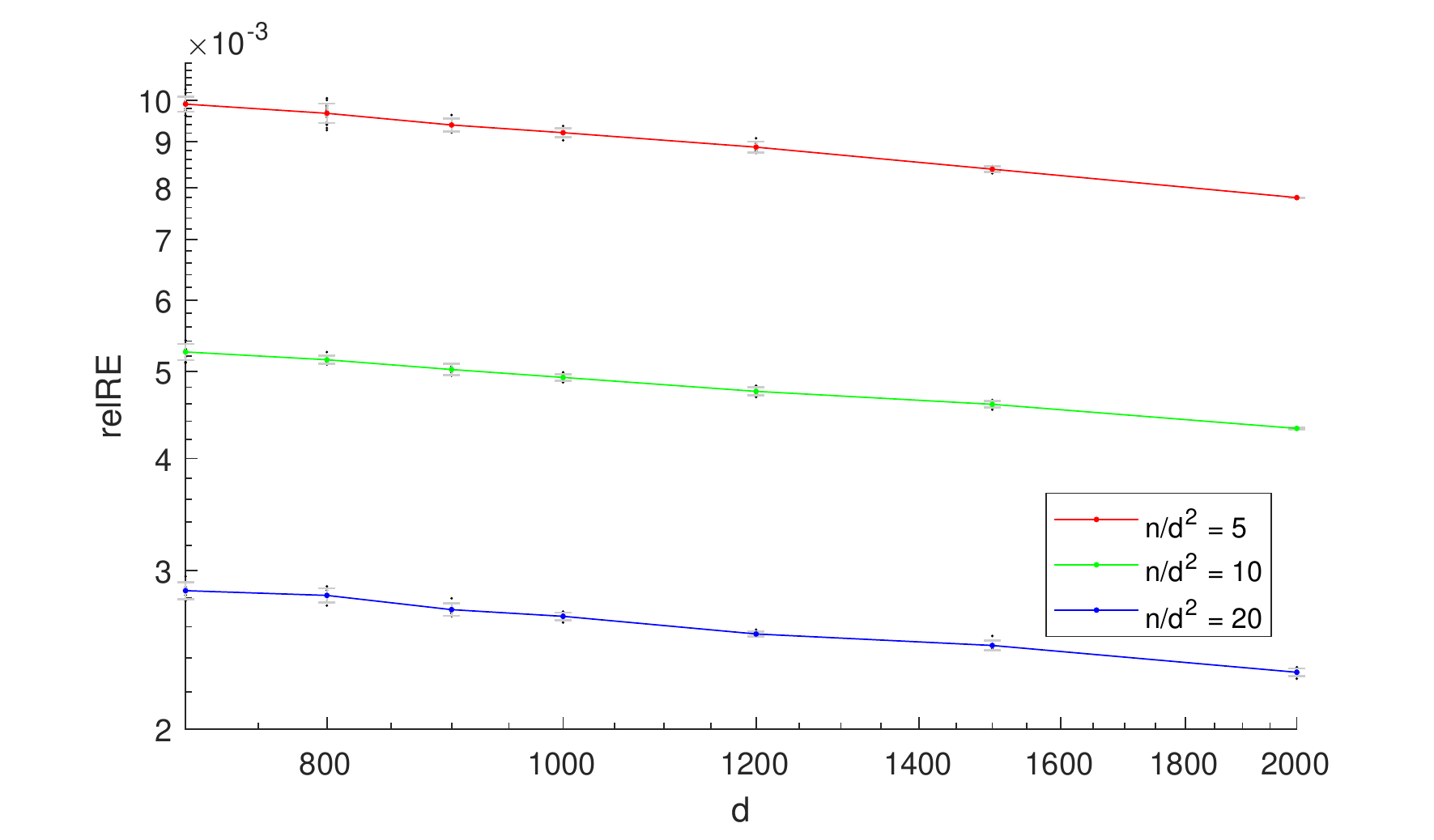}
        \caption{$relRE^*$-$d$ curves with $r/d = 0.01$, $n/d^2 = 5, 10$ or $20$. Each black point comes from one independent experiment. The grey bars stand for standard deviations (after logarithm). }
        \label{figure_d_2}
    \end{minipage}
\end{figure}

Markov chains with $d = 500, 700, 1000, 1500,2000$ and $r/d = 0.01$ are created. For each model, trajectories with sampling size $n/d^2 = 5, 10$ and $20$ are generated independently. A plot of the regularization paths is shown in Figure \ref{figure_d_1}. When the sparsity degree $r/d$ and valid sampling size $n/d^2$ are fixed, the shape of the regularization paths are almost identical, regardless of the dimension $d$. In Figure \ref{figure_d_2}, we plot the optimal relative recovery error $relRE^*$ against $d$. We can learn from the results that $relRE^*$ and $d$ have an approximately linear relationship in the log-log plot. The regression coefficients are shown in Table \ref{table_d_regression}. A summary of computational results is reported in Table \ref{table_d}. With the growth of $d$, $relSE$ increases whereas $relRE^*$ decreases. Therefore, when $r/d$ and $n/d^2$ are fixed, one can reconstruct the model more accurately if $d$ is larger.

\begin{table}
  \centering
  \begin{tabular}{cc|ccc}
    \hline
    $r/d$ & $ n/d^2 $ & {\footnotesize Constant $A_1$} & \begin{tabular}{c} {\footnotesize Regression} \\ {\footnotesize Coefficient $B_1$} \end{tabular} & {\footnotesize \begin{tabular}{c} Squared Correlation \\ Coefficient \end{tabular}} \\
    \hline
    $0.01$ & $5$ & $-3.115$ & $- 0.228$ & $0.997$ \\
    $0.01$ & $10$ & $-4.030$ & $-0.186$ & $0.998$ \\
    $0.01$ & $20$ & $-4.536$ & $-0.201$ & $0.996$ \\
    \hline
  \end{tabular}
  \caption{Coefficients in the linear regression $\ln relRE^* = A_1 + B_1 \ln d$.}\label{table_d_regression}
\end{table}

\begin{table}
    \centering
    {
    \begin{tabular}{c|cc|cc|cc}
        \hline
        $(d,\lambda^*)$ & \multicolumn{2}{c|}{($500$,$1.6\text{e-}07$)} & \multicolumn{2}{c|}{($700$,$6.2\text{e-}08$)} & \multicolumn{2}{c}{($1000$,$2.4\text{e-}08$)} \\
        \hline
        $\varepsilon$ & $1.5\text{e-}04$ & $5\text{e-}05$ & $1.5\text{e-}04$ & $5\text{e-}05$ & $1.5\text{e-}04$ & $5\text{e-}05$ \\
        $\hat{s}$ & $5$ & $13$ & $7$ & $18$ & $10$ & $22$ \\
        obj value & $4.70\text{e-}07$ & $4.70\text{e-}07$ & $2.10\text{e-}07$ & $2.09\text{e-}07$ & $9.16\text{e-}08$ & $9.15\text{e-}08$ \\
        $relRE^*$ & $5.70\text{e-}03$ & $5.96\text{e-}03$ & $5.85\text{e-}03$ & $5.83\text{e-}03$ & $5.45\text{e-}03$ & $5.61\text{e-}03$  \\
        time (${\rm s}$) & $19.5$ & $19.5$ & $28.9$ & $79.2$ & $49.4$ & $58.5$ \\
        \hline
        $relLE1$ & $1.6\text{e-}03$ & $1.8\text{e-}03$ & $2.1\text{e-}04$ & $8.2\text{e-}04$ & $1.7\text{e-}03$ & $5.6\text{e-}03$ \\
        $relLE2$ & $2.3\text{e-}03$ & $1.4\text{e-}03$ & $5.6\text{e-}04$ & $1.4\text{e-}03$ & $4.5\text{e-}03$ & $6.9\text{e-}03$ \\
        $GE$ & $1.3\text{e-}01$ & $4.9\text{e-}02$ & $1.6\text{e-}01$ & $5.6\text{e-}02$ & $2.0\text{e-}01$ & $8.8\text{e-}02$ \\
        $relDG$ & $7.1\text{e-}02$ & $2.7\text{e-}02$ & $7.9\text{e-}02$ & $2.8\text{e-}02$ & $8.7\text{e-}02$ & $3.9\text{e-}02$ \\
        \hline
        $relSE$ & \multicolumn{2}{c|}{$6.30\text{e-}02$} & \multicolumn{2}{c|}{$7.03\text{e-}02$} & \multicolumn{2}{c}{$7.73\text{e-}02$} \\
        $\sigma_1(\hat{\Xi}\hat{P}^{(n)})$ & \multicolumn{2}{c|}{$2.38\text{e-}03$} & \multicolumn{2}{c|}{$1.63\text{e-}03$} & \multicolumn{2}{c}{$1.10\text{e-}03$} \\
        $\sigma_r(\hat{\Xi}\hat{P}^{(n)})$ & \multicolumn{2}{c|}{$3.72\text{e-}04$} & \multicolumn{2}{c|}{$1.85\text{e-}04$} & \multicolumn{2}{c}{$8.61\text{e-}05$} \\
        $\sigma_{r+1}(\hat{\Xi}\hat{P}^{(n)})$ & \multicolumn{2}{c|}{$6.39\text{e-}05$} & \multicolumn{2}{c|}{$3.81\text{e-}05$} & \multicolumn{2}{c}{$2.16\text{e-}05$} \\
        \hline
    \end{tabular}

    \vspace{0.5cm}

    \begin{tabular}{c|cc|cc|cc}
        \hline
        $(d,\lambda^*)$ & \multicolumn{2}{c|}{($1200$,$1.5\text{e-}08$)} & \multicolumn{2}{c|}{($1500$,$7.5\text{e-}09$)} & \multicolumn{2}{c}{($2000$,$3.5\text{e-}09$)} \\
        \hline
        $\varepsilon$ & $1.5\text{e-}04$ & $5\text{e-}05$ & $1.5\text{e-}04$ & $5\text{e-}05$ & $1.5\text{e-}04$ & $5\text{e-}05$ \\
        $\hat{s}$ & $13$ & $27$ & $15$ & $33$ & $20$ & $42$ \\
        obj value & $6.09\text{e-}08$ & $6.09\text{e-}08$ & $3.53\text{e-}08$ & $3.55\text{e-}08$ & $1.86\text{e-}08$ & $1.86\text{e-}08$ \\
        $relRE^*$ & $5.43\text{e-}03$ & $5.77\text{e-}03$ & $4.51\text{e-}03$ & $5.21\text{e-}03$ & $4.31\text{e-}03$ & $4.38\text{e-}03$ \\
        time (${\rm s}$) & $90.0$ & $153.0$ & $256.5$ & $183.8$ & $257.0$ & $370.0$ \\
        \hline
        $relLE1$ & $5.7\text{e-}04$ & $3.2\text{e-}03$ & $2.6\text{e-}03$ & $2.1\text{e-}03$ & $8.8\text{e-}04$ & $2.9\text{e-}03$ \\
        $relLE2$ & $1.0\text{e-}03$ & $9.5\text{e-}03$ & $5.1\text{e-}03$ & $7.0\text{e-}03$ & $3.4\text{e-}03$ & $1.3\text{e-}02$ \\
        $GE$ & $2.0\text{e-}01$ & $5.2\text{e-}02$ & $2.1\text{e-}01$ & $9.9\text{e-}02$ & $2.6\text{e-}01$ & $1.6\text{e-}01$ \\
        $relDG$ & $8.5\text{e-}02$ & $2.2\text{e-}02$ & $8.0\text{e-}02$ & $3.8\text{e-}02$ & $9.5\text{e-}02$ & $5.8\text{e-}02$ \\
        \hline
        $relSE$ & \multicolumn{2}{c|}{$7.98\text{e-}02$} & \multicolumn{2}{c|}{$8.38\text{e-}02$} & \multicolumn{2}{c}{$8.76\text{e-}02$} \\
        $\sigma_1(\hat{\Xi}\hat{P}^{(n)})$ & \multicolumn{2}{c|}{$9.05\text{e-}04$} & \multicolumn{2}{c|}{$7.10\text{e-}04$} & \multicolumn{2}{c}{$5.24\text{e-}04$} \\
        $\sigma_r(\hat{\Xi}\hat{P}^{(n)})$ & \multicolumn{2}{c|}{$6.18\text{e-}05$} & \multicolumn{2}{c|}{$4.03\text{e-}05$} & \multicolumn{2}{c}{$2.23\text{e-}05$} \\
        $\sigma_{r+1}(\hat{\Xi}\hat{P}^{(n)})$ & \multicolumn{2}{c|}{$1.63\text{e-}05$} & \multicolumn{2}{c|}{$1.15\text{e-}05$} & \multicolumn{2}{c}{$7.36\text{e-}06$} \\
        \hline
    \end{tabular}
    }
    \caption{Numerical results with $r/d = 0.01$, $n/d^2 = 10$. The initial dimension $s_0 = 0.1d$.} \label{table_d}
\end{table}

\vspace{0.3cm}

\noindent {\bf Sparsity Degree $r/d$}

\vspace{0.1cm}

\begin{figure}
    \centering
    \begin{minipage}{0.495\linewidth}
        \includegraphics[width = \textwidth]{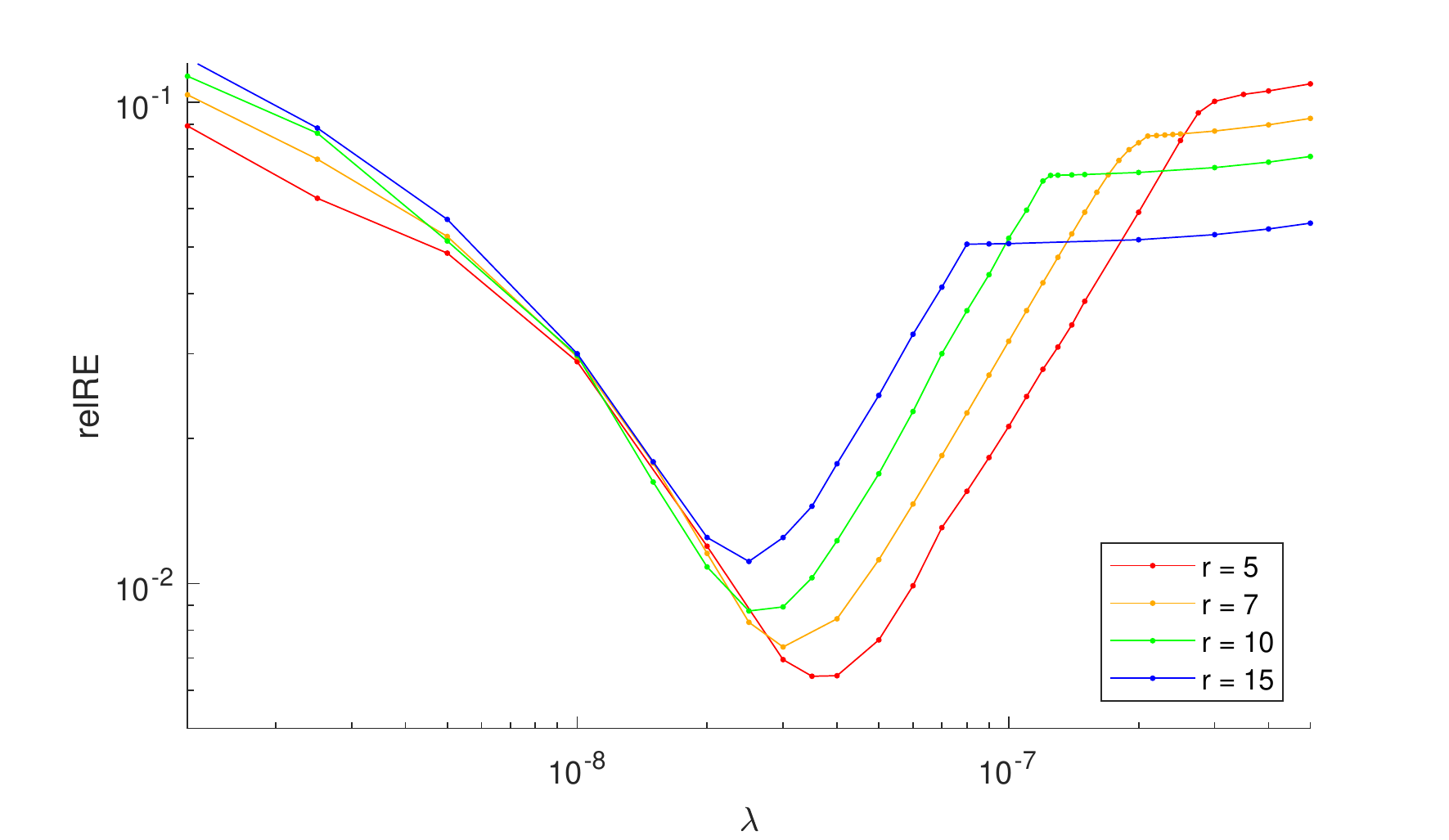}
        \caption{$relRE$-$\lambda$ curves with $d = 1000$, $n/d^2 = 5$, $r = 5$, $7$, $10$, or $1500$. \protect\\ ~ \protect\\ }
        \label{figure_r_1}
    \end{minipage}
    \begin{minipage}{0.495\linewidth}
        \includegraphics[width = \textwidth]{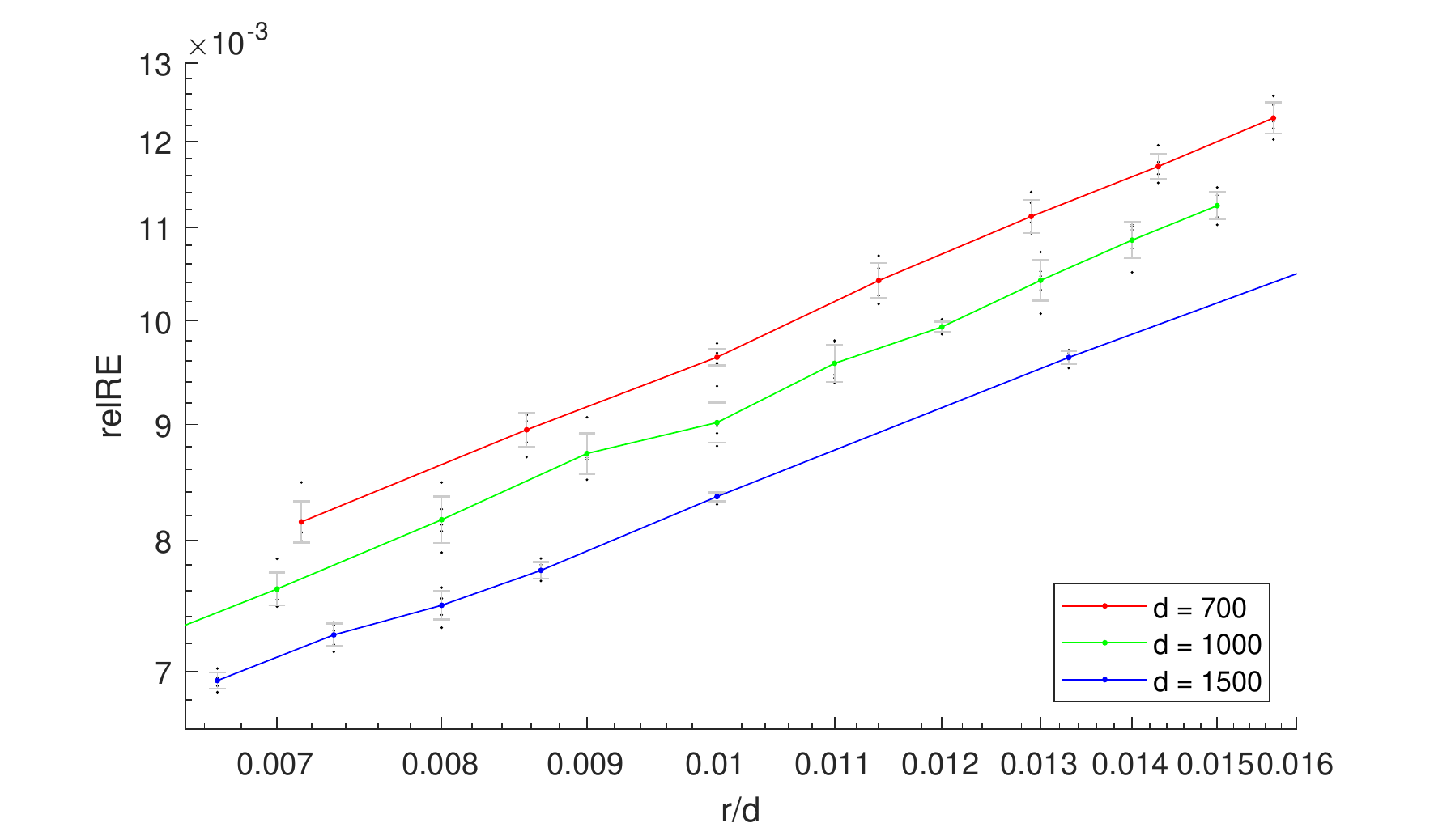}
        \caption{$relRE^*$-$r/d$ curves with $n/d^2 = 5$, $d = 700$, $1000$ or $1500$. Each black point comes from one independent experiment. The grey bars stand for standard deviations (after logarithm). }
        \label{figure_r_2}
    \end{minipage}
\end{figure}

Markov chains with $d = 1000$ and $r = 5, 7, 10, 15$ are created. The sampling size $n/d^2$ is taken to be $5$ or $10$. A plot of the regularization paths is shown in Figure \ref{figure_r_1}. We learn that, when $d$ and $n/d^2$ are fixed, one can recover the model more accurately if $r/d$ is smaller. In Figure \ref{figure_r_2}, we plot $relRE^*$ against $r/d$. There is also a linear relationship between $relRE^*$ and $d$ in the log-log plot. The linear regression coefficients are shown in Table \ref{table_r_regression}. More detailed computational results are presented in Table \ref{table_r}. When $r \ll d$, as $r/d$ gets smaller, $relSE$ decreases slightly, and $relRE^*$ drops even faster.

\begin{table}
  \centering
  \begin{tabular}{cc|ccc}
    \hline
    $ d $ & $n/d^2$ & {\footnotesize Constant $A_2$} & \begin{tabular}{c} {\footnotesize Regression} \\ {\footnotesize Coefficient $B_2$} \end{tabular} & {\footnotesize \begin{tabular}{c} Squared Correlation \\ Coefficient \end{tabular}} \\
    \hline
    $700$ & $5$ & $-0.962$ & $0.525$ & $0.9994$ \\
    $1000$ & $5$ & $-1.026$ & $0.507$ & $0.9990$ \\
    $1500$ & $5$ & $-1.122$ & $0.477$ & $0.9998$ \\
    \hline
  \end{tabular}
  \caption{Coefficients in the linear regression $\ln relRE^* = A_2 + B_2 \ln (r/d)$.}\label{table_r_regression}
\end{table}

\begin{table}
    \centering
    {
    \begin{tabular}{c|cc|cc}
        \hline
        $(r,\lambda^*)$ & \multicolumn{2}{c|}{$(5,3.2\text{e-}08)$} & \multicolumn{2}{c}{$(7,2.8\text{e-}08)$} \\
        \hline
        $\varepsilon$ & $1.5\text{e-}04$ & $5\text{e-}05$ & $1.5\text{e-}04$ & $5\text{e-}05$ \\
        $\hat{s}$ & $5$ & $14$ & $8$ & $19$ \\
        obj value & $1.04\text{e-}07$ & $1.04\text{e-}07$ & $9.83\text{e-}08$ & $9.82\text{e-}08$ \\
        $relRE^*$ & $3.77\text{e-}03$ & $3.86\text{e-}03$ & $4.62\text{e-}03$ & $4.74\text{e-}03$ \\
        time (${\rm s}$) & $63.3$ & $68.3$ & $85.8$ & $77.1$ \\
        \hline
        $relLE1$ & $5.0\text{e-}03$ & $1.8\text{e-}03$ & $1.2\text{e-}03$ & $4.8\text{e-}03$ \\
        $relLE2$ & $8.6\text{e-}03$ & $1.6\text{e-}03$ & $2.3\text{e-}03$ & $2.2\text{e-}03$ \\
        $GE$ & $1.9\text{e-}01$ & $7.6\text{e-}02$ & $1.2\text{e-}01$ & $5.5\text{e-}02$ \\
        $relDG$ & $9.3\text{e-}02$ & $3.8\text{e-}02$ & $5.6\text{e-}02$ & $2.6\text{e-}02$ \\
        \hline
        $relSE$ & \multicolumn{2}{c|}{$6.31\text{e-}02$} & \multicolumn{2}{c}{$6.99\text{e-}02$} \\
        $\sigma_1(\hat{\Xi}\hat{P}^{(n)})$ & \multicolumn{2}{c|}{$1.20\text{e-}03$} & \multicolumn{2}{c}{$1.14\text{e-}03$} \\
        $\sigma_r(\hat{\Xi}\hat{P}^{(n)})$ & \multicolumn{2}{c|}{$1.73\text{e-}04$} & \multicolumn{2}{c}{$1.28\text{e-}04$} \\
        $\sigma_{r+1}(\hat{\Xi}\hat{P}^{(n)})$ & \multicolumn{2}{c|}{$2.33\text{e-}05$} & \multicolumn{2}{c}{$2.27\text{e-}05$} \\
        \hline
    \end{tabular}

    \vspace{0.5cm}

    \begin{tabular}{c|cc|cc}
        \hline
        $(r,\lambda^*)$ & \multicolumn{2}{c|}{$(10,2.4\text{e-}08)$} & \multicolumn{2}{c}{$(15,2.1\text{e-}08)$} \\
        \hline
        $\varepsilon$ & $1.5\text{e-}04$ & $5\text{e-}05$ & $1.5\text{e-}04$ & $5\text{e-}05$ \\
        $\hat{s}$ & $10$ & $22$ & $16$ & $39$ \\
        obj value & $9.16\text{e-}08$ & $9.15\text{e-}08$ & $8.42\text{e-}08$ & $8.58\text{e-}08$ \\
        $relRE^*$ & $5.45\text{e-}03$ & $5.61\text{e-}03$ & $6.20\text{e-}03$ & $6.67\text{e-}03$ \\
        time (${\rm s}$) & $49.4$ & $58.5$ & $65.4$ & $136.4$ \\
        \hline
        $relLE1$ & $1.7\text{e-}03$ & $5.6\text{e-}03$ & $2.3\text{e-}04$ & $5.6\text{e-}04$ \\
        $relLE2$ & $4.5\text{e-}03$ & $6.9\text{e-}03$ & $1.1\text{e-}03$ & $1.7\text{e-}03$ \\
        $GE$ & $2.0\text{e-}01$ & $8.8\text{e-}02$ & $1.5\text{e-}01$ & $3.3\text{e-}02$ \\
        $relDG$ & $8.7\text{e-}02$ & $3.9\text{e-}02$ & $6.3\text{e-}02$ & $1.4\text{e-}02$ \\
        \hline
        $relSE$ & \multicolumn{2}{c|}{$7.73\text{e-}02$} & \multicolumn{2}{c}{$8.36\text{e-}02$} \\
        $\sigma_1(\hat{\Xi}\hat{P}^{(n)})$ & \multicolumn{2}{c|}{$1.10\text{e-}03$} & \multicolumn{2}{c}{$1.06\text{e-}03$} \\
        $\sigma_r(\hat{\Xi}\hat{P}^{(n)})$ & \multicolumn{2}{c|}{$8.61\text{e-}05$} & \multicolumn{2}{c}{$5.74\text{e-}05$} \\
        $\sigma_{r+1}(\hat{\Xi}\hat{P}^{(n)})$ & \multicolumn{2}{c|}{$2.16\text{e-}05$} & \multicolumn{2}{c}{$2.09\text{e-}05$} \\
        \hline
    \end{tabular}}
    \caption{Numerical results with $d = 1000$, $n/d^2 = 10$. The initial dimension $s_0 = 0.1d$.} \label{table_r}
\end{table}

\vspace{0.3cm}

\noindent {\bf Sample size $n/d^2$}

\vspace{0.1cm}

\begin{figure}
    \centering
    \begin{minipage}{0.495\linewidth}
        \includegraphics[width = \textwidth]{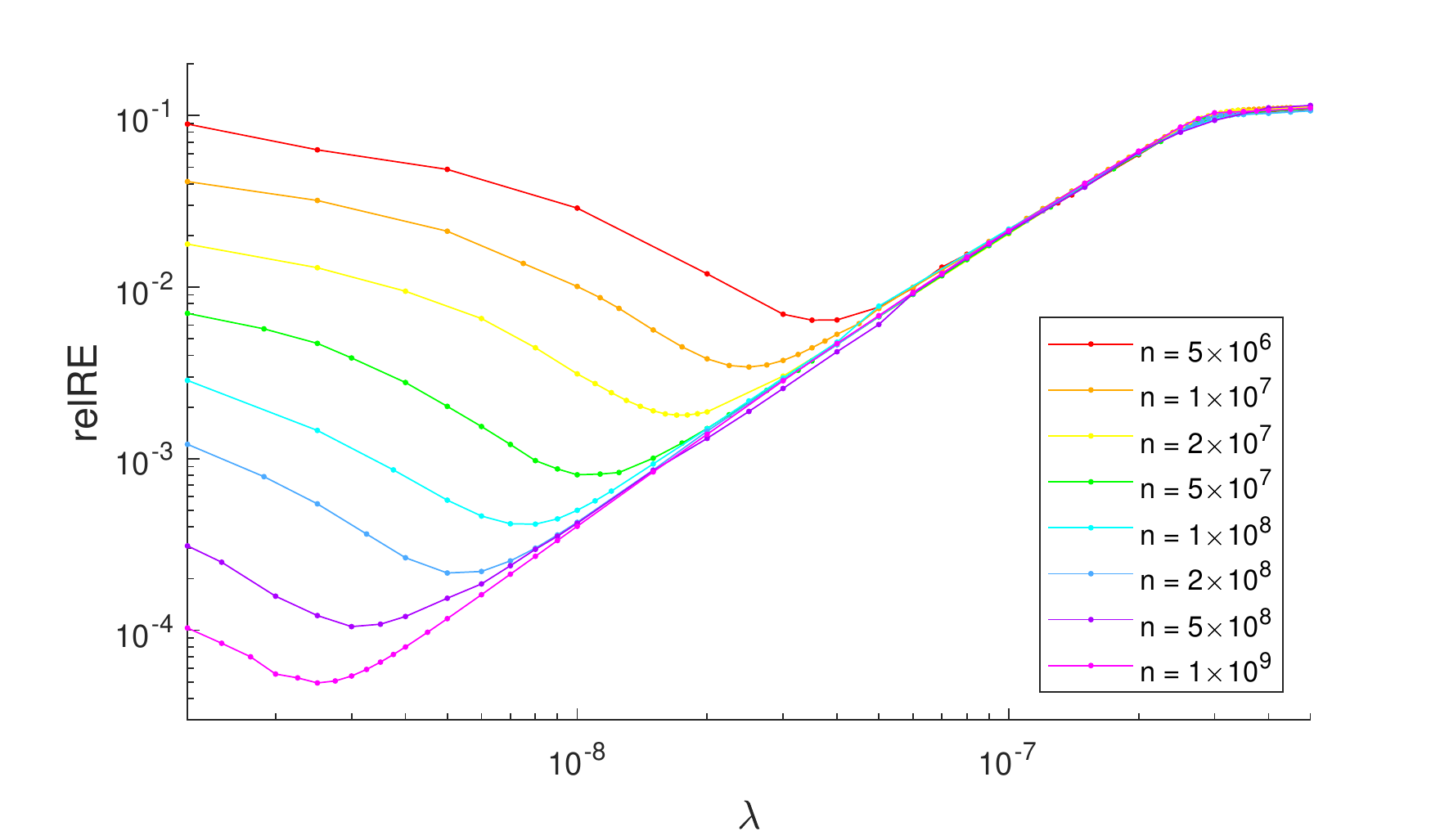}
        \caption{$relRE$-$\lambda$ curves with $d = 1000$, $r/d = 0.005$, $n/d^2 = 5$, $10$, $20$, $50$, $100$, $200$, $500$ or $1500$. \protect\\ ~ \protect\\ }
        \label{figure_n_1}
    \end{minipage}
    \begin{minipage}{0.495\linewidth}
        \includegraphics[width = \textwidth]{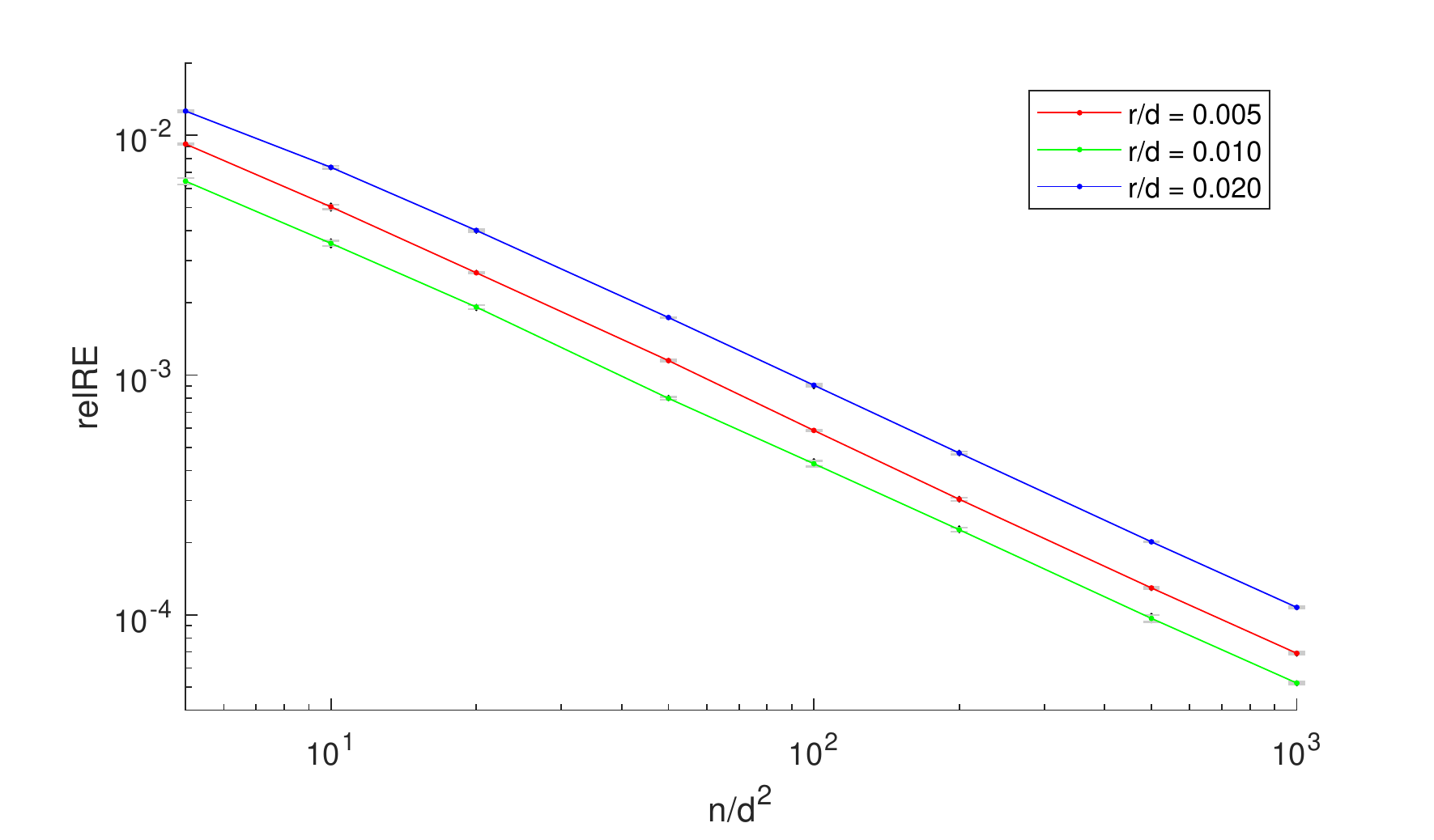}
        \caption{$relRE$-$n/d^2$ curves with $d = 1000$, $r/d = 0.005$, $0.010$ or $0.020$. Each black point comes from one independent experiment. The grey bars stand for standard deviations (after logarithm). }
        \label{figure_n_2}
    \end{minipage}
\end{figure}

A Markov chain with $d = 1000$ and $r = 5$ is created. We generate trajectories with $n/d^2 = 5, 10, 20, 50, 100, 200, 500, 1000$. A plot of regularization paths is shown in Figure \ref{figure_n_1}. As the sampling size $n$ grows, the empirical transition matrix $\hat{P}^{(n)}$ gets closer and closer to the ground-truth $P^*$, therefore, $\lambda^*$ is smaller and $relRE^*$ reduces to zero. In Figure \ref{figure_n_2}, we see the linear relationships between $relRE^*$ and $n$. The linear regression coefficients are shown in Table \ref{table_n_regression}. By the central limit theorem of Markov chain, $relSE$ is of order $n^{-1}$. However, the regression coefficient $B_3$ is slightly larger than $-1$. This results from the bias introduced by the regularization term. Numerical details are summarized in Table \ref{table_n}.

\begin{table}
  \centering
  \begin{tabular}{cc|ccc}
    \hline
    $d$ & $ r/d $ & {\footnotesize Constant $A_3$} & \begin{tabular}{c} {\footnotesize Regression} \\ {\footnotesize Coefficient $B_3$} \end{tabular} & {\footnotesize \begin{tabular}{c} Squared Correlation \\ Coefficient \end{tabular}} \\
    \hline
    $1000$ & $0.005$ & $-1.372$ & $- 0.930$ & $0.99987$ \\
    $1000$ & $0.010$ & $-1.539$ & $-0.916$ & $0.99990$ \\
    $1000$ & $0.015$ & $-1.231$ & $-0.910$ & $0.99950$ \\
    \hline
  \end{tabular}
  \caption{Coefficients in the linear regression $\ln relRE^* = A_3 + B_3 \ln (n/d^2)$.}\label{table_n_regression}
\end{table}

\begin{table}
    \centering
    {
    \begin{tabular}{c|cc|cc|cc}
        \hline
        $(n,\lambda^*)$ & \multicolumn{2}{c|}{$(5\text{e+}06,4.1\text{e-}08)$} & \multicolumn{2}{c|}{$(1\text{e+}07,3.2\text{e-}08)$} & \multicolumn{2}{c}{$(2\text{e+}07,2.4\text{e-}08)$} \\
        \hline
        $\varepsilon$ & $1.5\text{e-}04$ & $5\text{e-}05$ & $1.5\text{e-}04$ & $5\text{e-}05$ & $1.5\text{e-}04$ & $5\text{e-}05$ \\
        $\hat{s}$ & $6$ & $17$ & $5$ & $14$ & $5$ & $14$ \\
        obj value & $1.67\text{e-}07$ & $1.69\text{e-}07$ & $1.04\text{e-}07$ & $1.04\text{e-}07$ & $6.67\text{e-}08$ & $6.67\text{e-}08$ \\
        $relRE^*$ & $5.77\text{e-}03$ & $6.40\text{e-}03$ & $3.77\text{e-}03$ & $3.86\text{e-}03$ & $2.33\text{e-}03$ & $2.34\text{e-}03$ \\
        time (${\rm s}$) & $87.3$ & $108.8$ & $63.3$ & $68.3$ & $32.5$ & $38.4$ \\
        \hline
        $relLE1$ & $4.0\text{e-}03$ & $2.0\text{e-}03$ & $5.0\text{e-}03$ & $1.8\text{e-}03$ & $1.6\text{e-}03$ & $1.9\text{e-}03$ \\
        $relLE2$ & $5.3\text{e-}03$ & $1.4\text{e-}03$ & $8.6\text{e-}03$ & $1.6\text{e-}03$ & $2.8\text{e-}03$ & $2.4\text{e-}03$ \\
        $GE$ & $2.3\text{e-}01$ & $8.8\text{e-}02$ & $1.9\text{e-}01$ & $7.6\text{e-}02$ & $1.4\text{e-}01$ & $3.9\text{e-}02$ \\
        $relDG$ & $9.3\text{e-}02$ & $3.7\text{e-}02$ & $9.3\text{e-}02$ & $3.8\text{e-}02$ & $8.3\text{e-}02$ & $2.4\text{e-}02$ \\
        \hline
        $relSE$ & \multicolumn{2}{c|}{$1.26\text{e-}01$} & \multicolumn{2}{c|}{$6.31\text{e-}02$} & \multicolumn{2}{c}{$3.15\text{e-}02$} \\
        $\sigma_1(\hat{\Xi}\hat{P}^{(n)})$ & \multicolumn{2}{c|}{$1.20\text{e-}03$} & \multicolumn{2}{c|}{$1.20\text{e-}03$} & \multicolumn{2}{c}{$1.20\text{e-}03$} \\
        $\sigma_r(\hat{\Xi}\hat{P}^{(n)})$ & \multicolumn{2}{c|}{$1.74\text{e-}04$} & \multicolumn{2}{c|}{$1.73\text{e-}04$} & \multicolumn{2}{c}{$1.87\text{e-}04$} \\
        $\sigma_{r+1}(\hat{\Xi}\hat{P}^{(n)})$ & \multicolumn{2}{c|}{$3.29\text{e-}05$} & \multicolumn{2}{c|}{$2.33\text{e-}05$} & \multicolumn{2}{c}{$1.70\text{e-}05$} \\
        \hline
    \end{tabular}

    \vspace{0.5cm}

    \begin{tabular}{c|cc|cc|cc}
        \hline
        $(n,\lambda^*)$ & \multicolumn{2}{c|}{$(5\text{e+}07,1.7\text{e-}08)$} & \multicolumn{2}{c|}{$(1\text{e+}08,1.2\text{e-}08)$} & \multicolumn{2}{c}{$(2\text{e+}08,9.0\text{e-}09)$} \\
        \hline
        $\varepsilon$ & $1.5\text{e-}04$ & $5\text{e-}05$ & $1.5\text{e-}04$ & $5\text{e-}05$ & $1.5\text{e-}04$ & $5\text{e-}05$ \\
        $\hat{s}$ & $5$ & $17$ & $5$ & $18$ & $5$ & $12$ \\
        obj value & $3.97\text{e-}08$ & $3.97\text{e-}08$ & $2.61\text{e-}08$ & $2.62\text{e-}08$ & $1.85\text{e-}08$ & $1.85\text{e-}08$ \\
        $relRE^*$ & $1.18\text{e-}03$ & $1.20\text{e-}03$ & $6.37\text{e-}04$ & $6.83\text{e-}04$ & $3.91\text{e-}04$ & $3.94\text{e-}04$ \\
        time (${\rm s}$) & $42.0$ & $47.9$ & $42.7$ & $28.6$ & $43.0$ & $45.8$ \\
        \hline
        $relLE1$ & $2.4\text{e-}03$ & $2.5\text{e-}03$ & $3.2\text{e-}03$ & $3.8\text{e-}03$ & $7.4\text{e-}03$ & $2.0\text{e-}02$ \\
        $relLE2$ & $3.9\text{e-}03$ & $3.7\text{e-}03$ & $4.8\text{e-}03$ & $2.6\text{e-}03$ & $1.2\text{e-}02$ & $5.5\text{e-}03$ \\
        $GE$ & $6.5\text{e-}02$ & $1.2\text{e-}02$ & $7.0\text{e-}02$ & $1.5\text{e-}02$ & $2.7\text{e-}02$ & $3.7\text{e-}02$ \\
        $relDG$ & $4.6\text{e-}02$ & $8.4\text{e-}03$ & $5.2\text{e-}02$ & $1.2\text{e-}02$ & $2.1\text{e-}02$ & $2.0\text{e-}02$ \\
        \hline
        $relSE$ & \multicolumn{2}{c|}{$1.27\text{e-}02$} & \multicolumn{2}{c|}{$6.37\text{e-}03$} & \multicolumn{2}{c}{$3.23\text{e-}03$} \\
        $\sigma_1(\hat{\Xi}\hat{P}^{(n)})$ & \multicolumn{2}{c|}{$1.20\text{e-}03$} & \multicolumn{2}{c|}{$1.20\text{e-}03$} & \multicolumn{2}{c}{$1.20\text{e-}03$} \\
        $\sigma_r(\hat{\Xi}\hat{P}^{(n)})$ & \multicolumn{2}{c|}{$1.72\text{e-}04$} & \multicolumn{2}{c|}{$1.71\text{e-}04$} & \multicolumn{2}{c}{$1.72\text{e-}04$} \\
        $\sigma_{r+1}(\hat{\Xi}\hat{P}^{(n)})$ & \multicolumn{2}{c|}{$1.04\text{e-}05$} & \multicolumn{2}{c|}{$7.33\text{e-}06$} & \multicolumn{2}{c}{$5.25\text{e-}06$} \\
        \hline
    \end{tabular}}
    \caption{Numerical results with $d = 1000$, $r = 5$. The initial dimension $s_0 = 0.1d$.} \label{table_n}
\end{table}

\subsection{Experiments with Manhattan Taxi Data}

We use the state aggregation model to partition Manhattan transportation network into different regions.
Our experiment is based on a real dataset of $1.1 \times 10^7$ NYC Yellow cab trips in January 2016 \cite{NYCyellowcabJan2016}. Each record includes passenger pick-up and drop-off information (coordinates, time, etc.) of one trip.
We want to construct a stochastic model of the traffic flow.
The movements of taxis are nearly memoryless. Therefore, we admit that the stochastic system satisfies Markov property, and approximate it by a finite-state Markov chain. We divide the map into a fine grid and merge the locations in the same cell into one state. Each trip is a sampled one-step state transition between cells.
We formulate an empirical Markov transition matrix $\hat{P}$ and an empirical stationary distribution $\hat{\xi}$ based on the dataset, and seek for a method to simplify the stochastic system.

We further assume that the transportation system is driven by a Markov chain with fewer states that are invisible and are aggregations of the states in the original Markov chain. In order to identify the state aggregation structure, we apply model \eqref{reformM_new} to the estimated $\hat{P}$ and $\hat{\xi}$. The solution $\hat{U}$ and $\hat{V}$ in \eqref{reformM_new} help embed the states in the original Markov chain into an $s$-dimensional space. We then cluster the states according to their coordinates, which yields a partition of the transportation network. In this problem, it is desirable to have a small rank $s$, since clustering algorithms are very likely to fail in a high-dimensional space. For this reason, we choose the early stopping rule and start from $s_0=1$ when applying Algorithm \ref{MetaAlg}.

We now introduce some specific settings in our numerical experiments.
When preprocessing data, we first delete the records beyond area $80.07^{\circ} {\rm W} \sim 60.92^{\circ} {\rm W}$, $30.69^{\circ} {\rm N} \sim 50.85^{\circ} {\rm N}$,
and divide the rectangular into $0.001^{\circ} \times 0.001^{\circ}$ grid.
We count the number of times a state appears as a pick-up place, and discard the ones with frequency of occurrence lower than $10^{-4}$. We end up with $2017$ valid states denoted by $\mathcal{S} = \{1,2,\cdots,2017\}$, and $7.5 \times 10^{6}$ remaining records $\left(s_t^{\rm pickup}, s_t^{\rm dropoff}\right) \in \mathcal{S}^2$, $t = 1,2,3, \cdots, T$.
The empirical $\hat{P}$ and $\hat{\xi}$ are formulated by
\[ \hat{P}_{ij} = \frac{\sum_{t=1}^T \mathbbm{1}[s_t^{\rm pickup}=i, s_t^{\rm dropoff}=j]}{\sum_{t=1}^T\mathbbm{1}[s_t^{\rm pickup}=i]}, \quad i,j \in \mathcal{S}, \]
\[ \hat{\xi}_i = \frac{\sum_{t=1}^T \mathbbm{1}[s_t^{\rm pickup}=i]}{T}, \quad i \in \mathcal{S}. \]
We use MATLAB function \texttt{kmeans} to implement $k$-means++ algorithm, and specify $500$ replicates to help find a lower local minimum.

In order to illustrate the reliability of the state aggregation model, we compare our results with the batch partition procedure in \cite{2017arXiv170507881Y}, which provides another method to embed the states in a reduced-dimensional space. The network partition algorithm in \cite{2017arXiv170507881Y} begins with optimal rank-$s$ approximation of matrix $\hat{\Xi}\hat{P}$ given by SVD, which we refer to as $\check{U} \check{D} \check{V}^T$.
Here, $s$ is the number of parts that the map will be partitioned into.
It has been proved in \cite{Weinan2008Optimal} that, for a lumpable Markov chain, one can obtain the optimal partition of state space with the clustering results of $\hat{\Xi}^{-1}\check{U}$ or $\hat{\Xi}^{-1}\check{V}$.
In our numerical experiments, we used MATLAB function \texttt{svds} to implement approximate SVD and \texttt{kmeans} with $500$ replicates to cluster the states.

\begin{figure}
    \centering
    \includegraphics[width=0.225\textwidth]{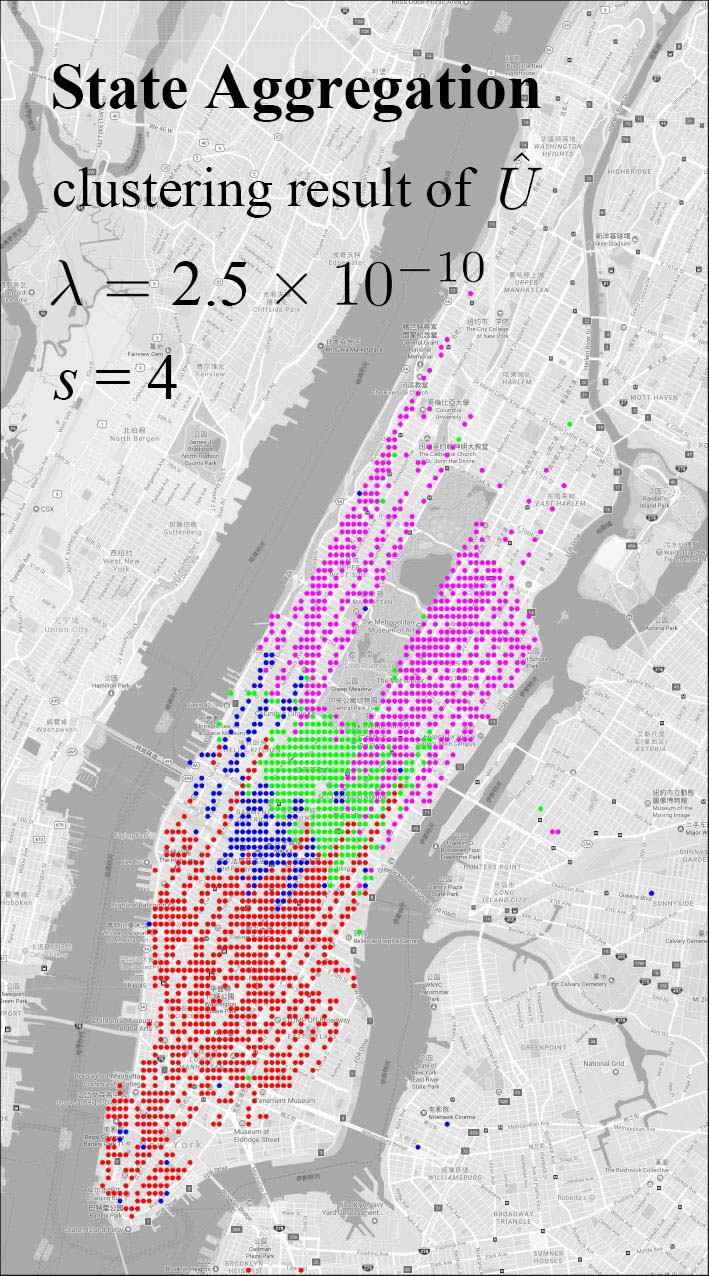} \ \includegraphics[width=0.225\textwidth]{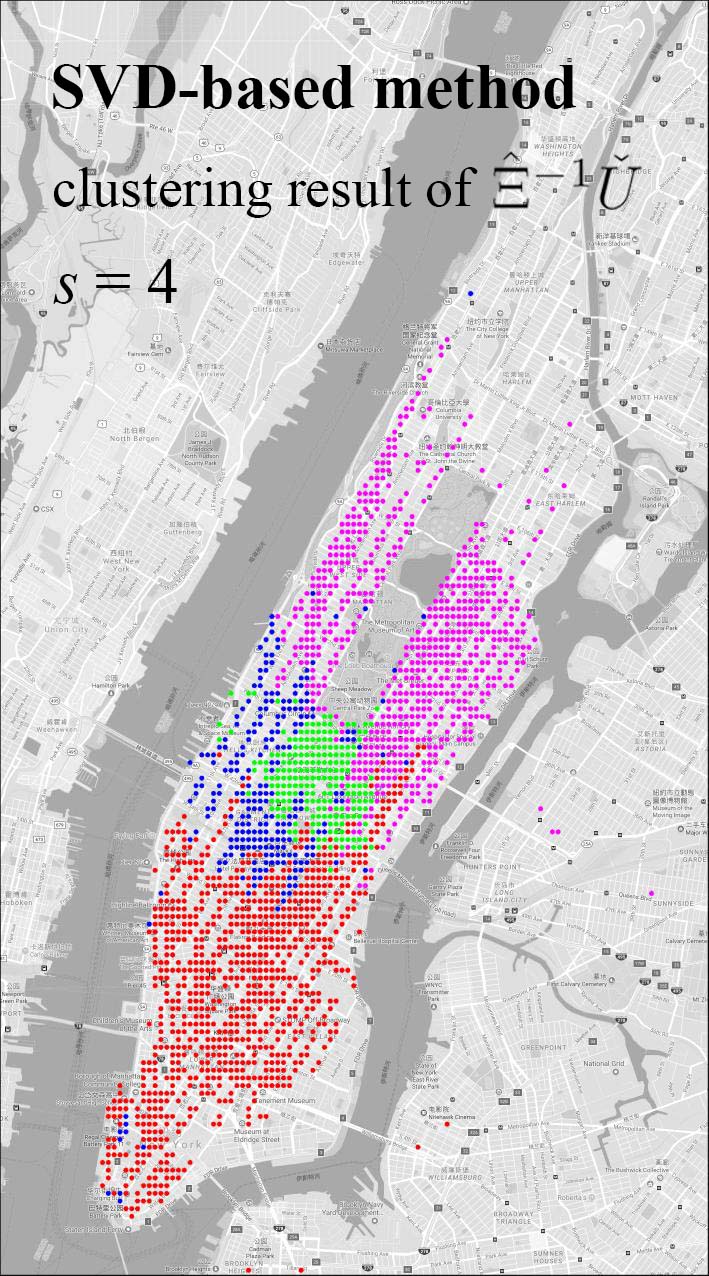} \hspace{0.5cm}
    \includegraphics[width=0.225\textwidth]{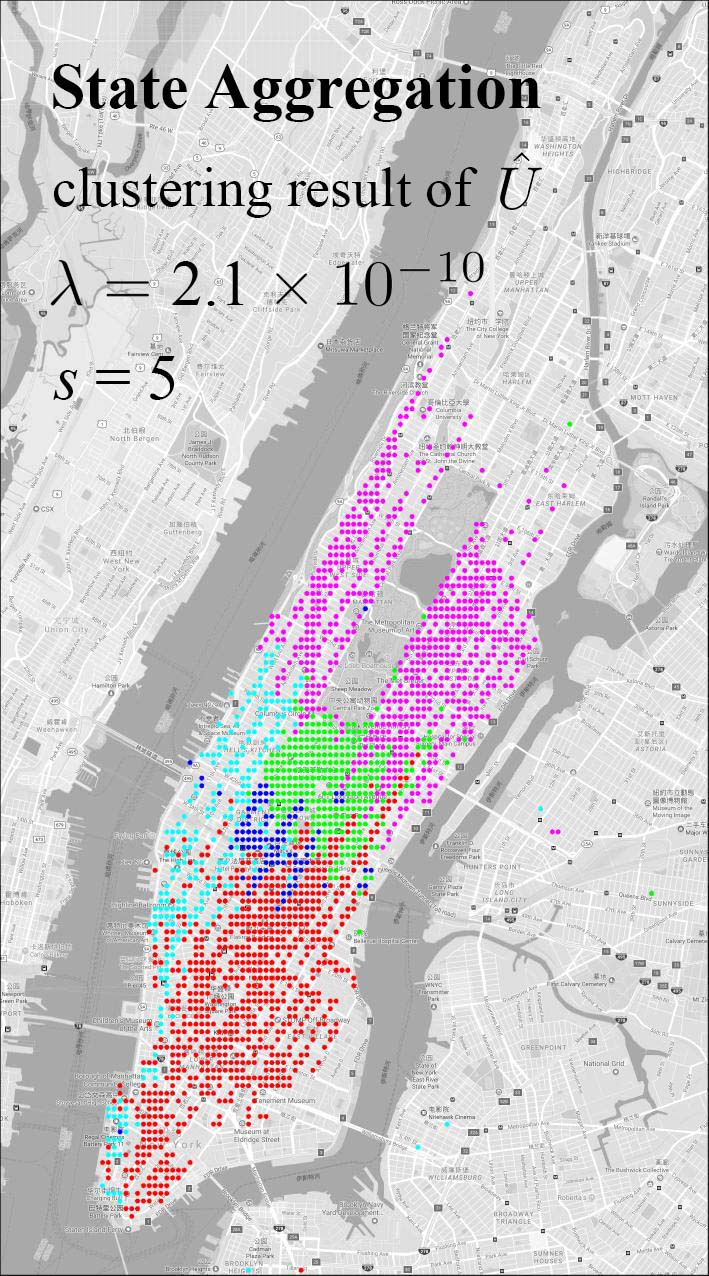} \ \includegraphics[width=0.225\textwidth]{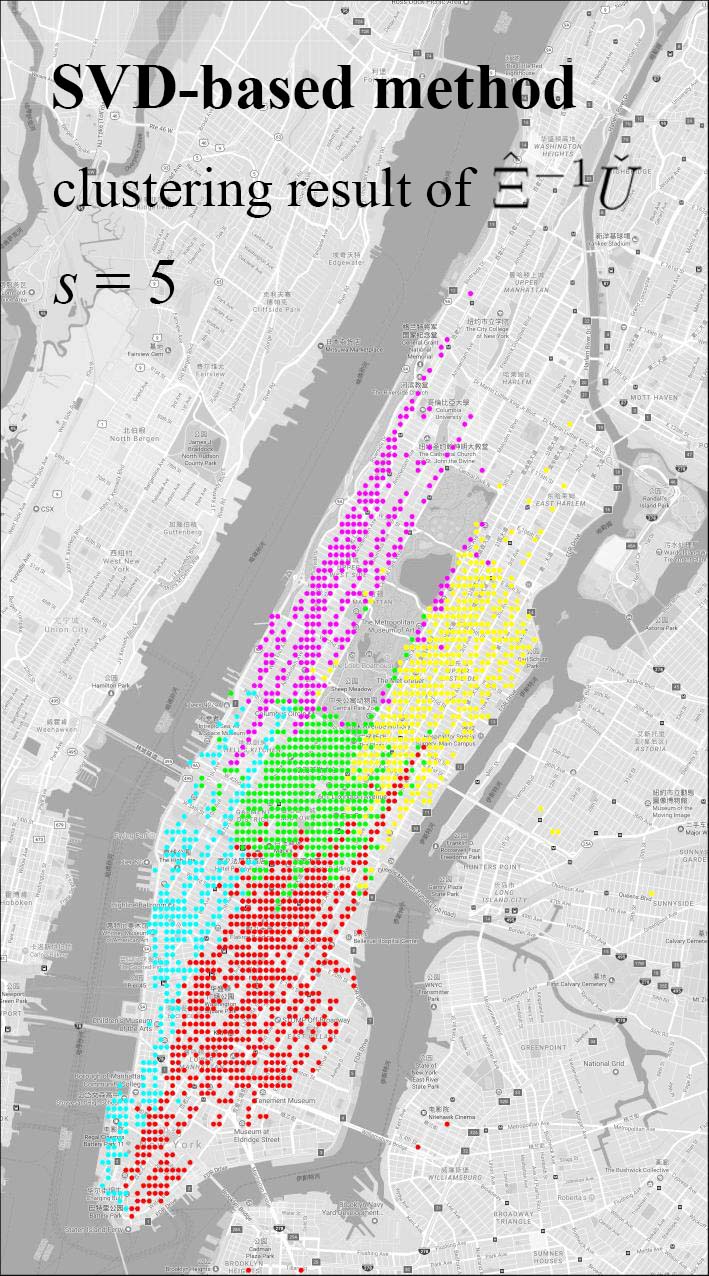} \\
    \vspace{0.7cm}
    \includegraphics[width=0.225\textwidth]{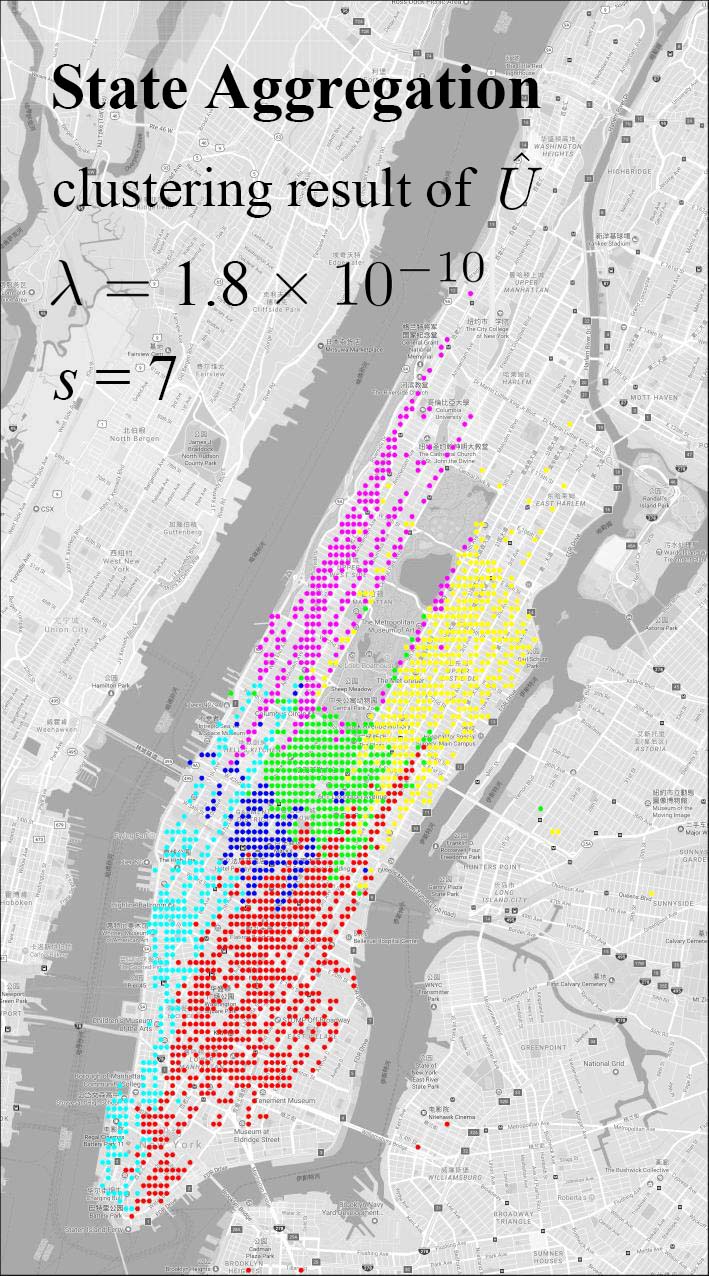} \ \includegraphics[width=0.225\textwidth]{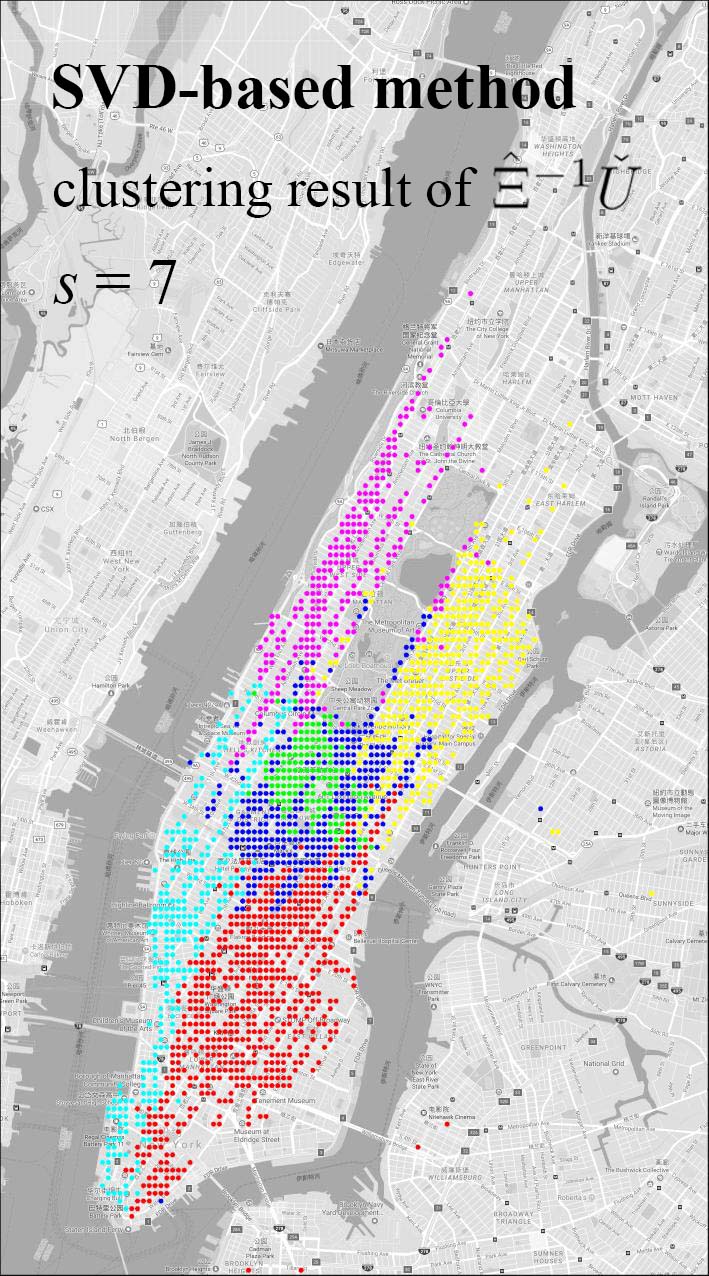} \hspace{0.5cm}
    \includegraphics[width=0.225\textwidth]{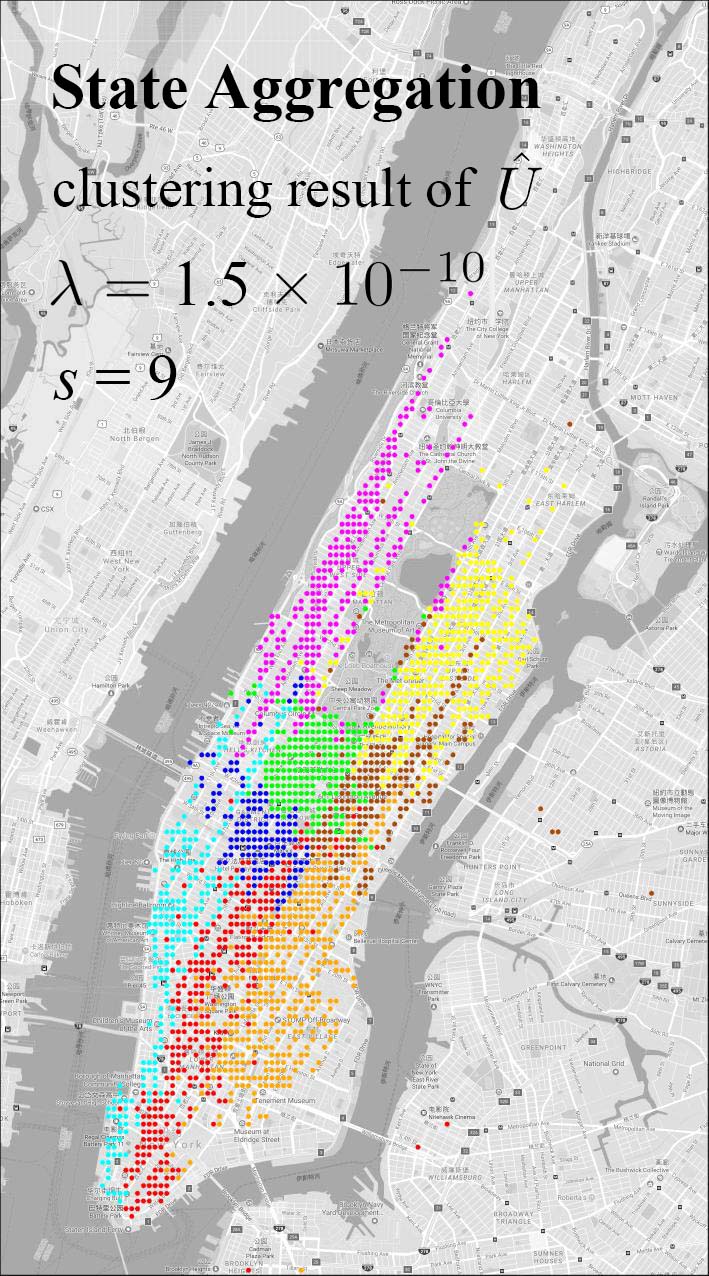} \ \includegraphics[width=0.225\textwidth]{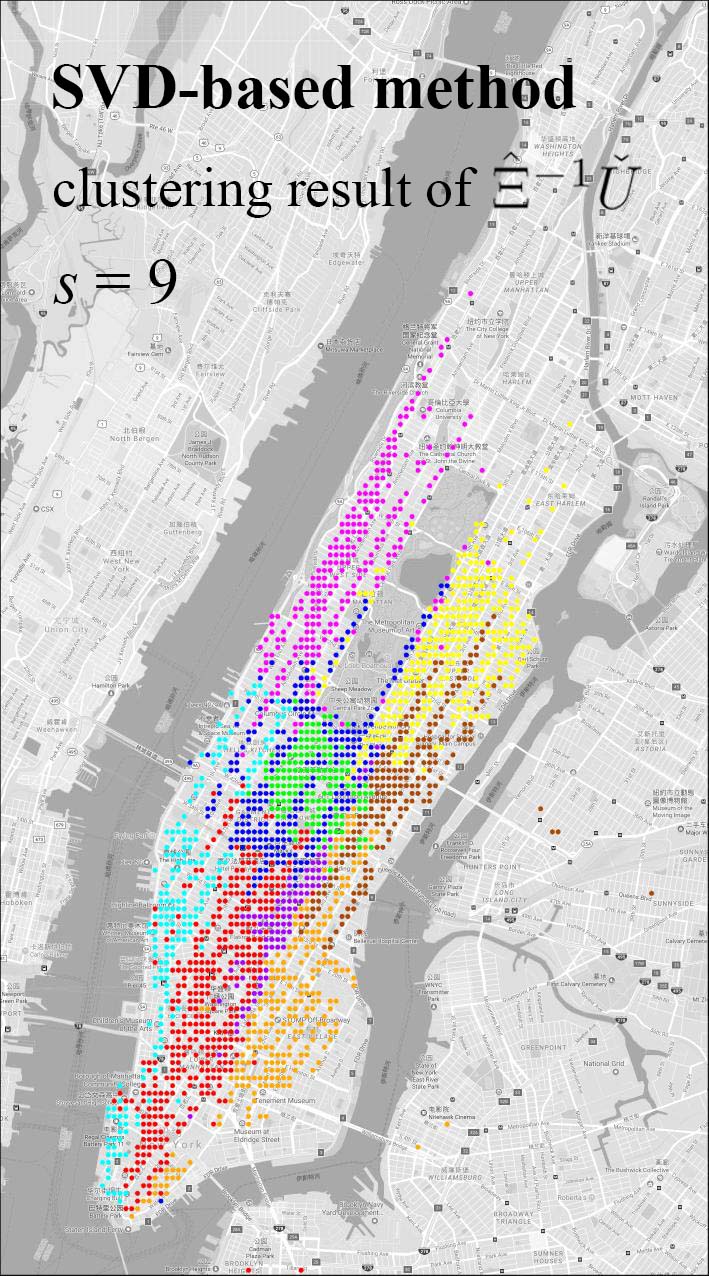} \\
    \vspace{0.7cm}
    \includegraphics[width=0.225\textwidth]{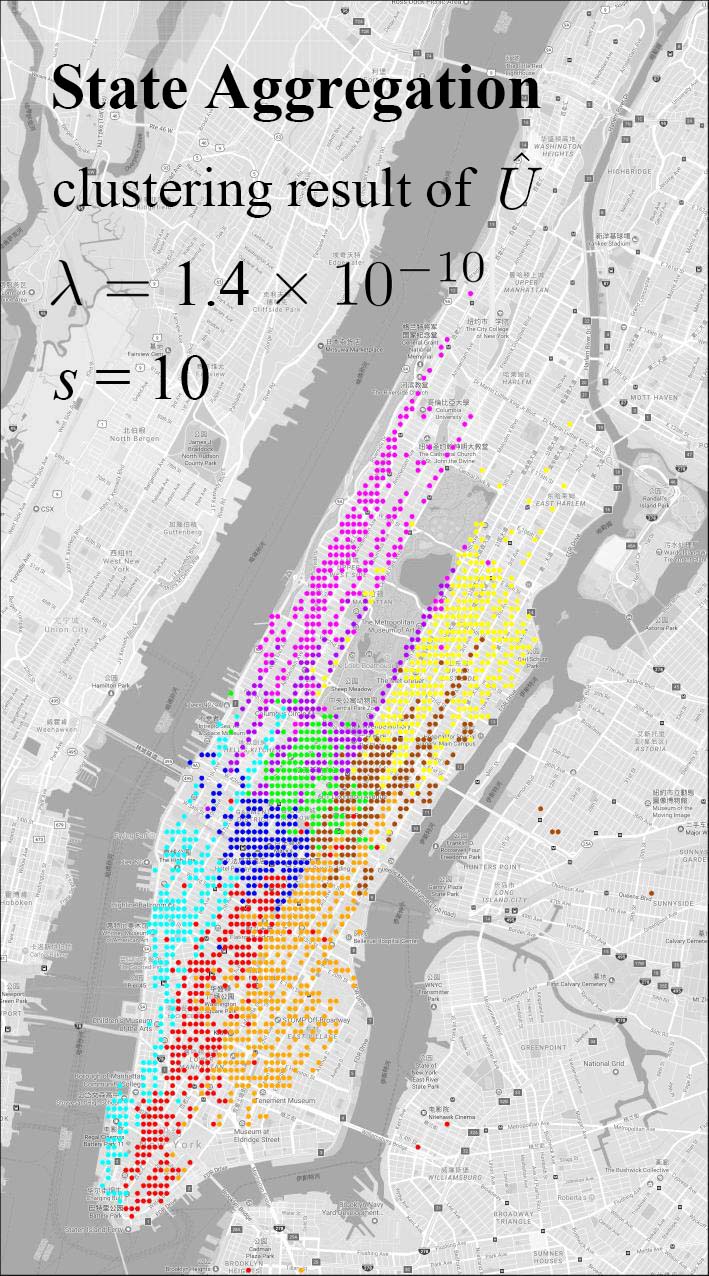} \ \includegraphics[width=0.225\textwidth]{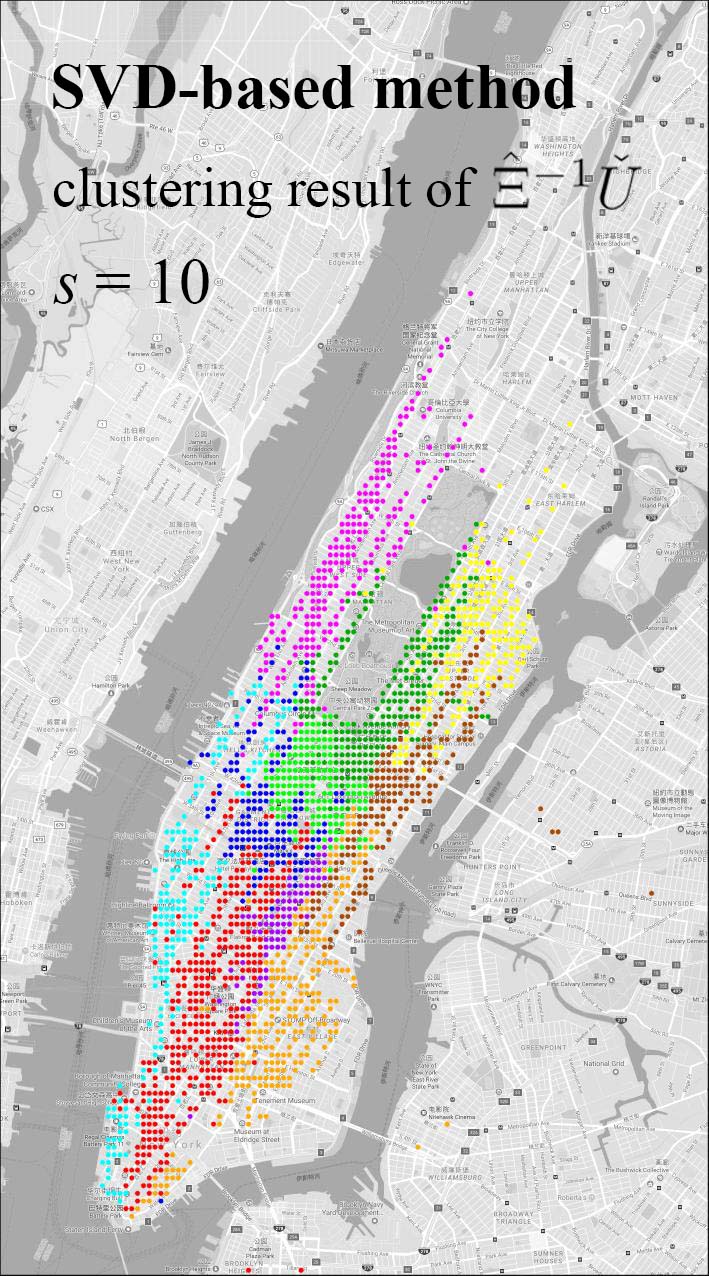} \hspace{0.5cm}
    \includegraphics[width=0.225\textwidth]{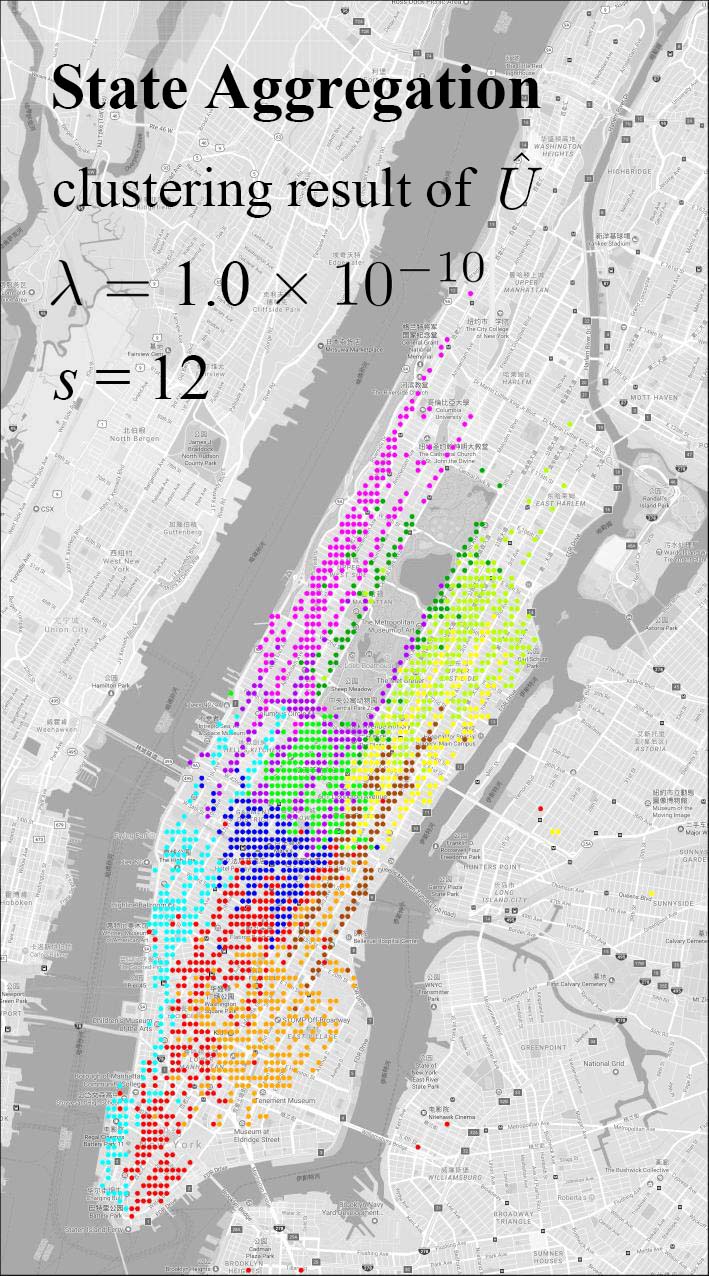} \ \includegraphics[width=0.225\textwidth]{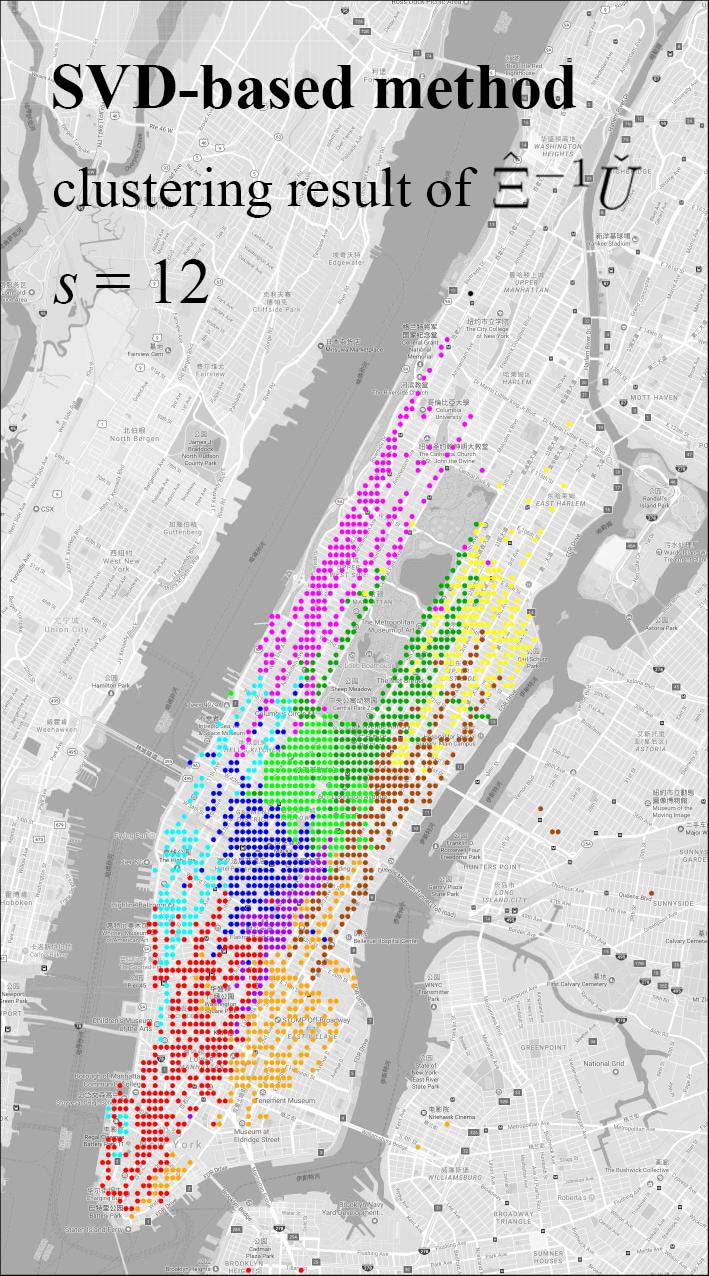} \\
    \caption{The partition of Manhattan transportation network. Each point in the map represents a valid state of the Markov chain. The figures in one pair have exactly the same number of regions, where the left one is produced by the state aggregation model and the right one is provided by the SVD-based method in \cite{2017arXiv170507881Y}. In some figures, there are less than $s$ regions appearing on the map, because some points are plotted beyond the boundaries.}
    \label{Manhattan}
\end{figure}

The zoning results of $\hat{U}$ in \eqref{reformM_new} and
$\hat{\Xi}^{-1}\check{U}$ in \cite{2017arXiv170507881Y} are shown in Figure
\ref{Manhattan}. Overall, the results of our approach are more aligned to
geometric location, while the partition results of the SVD-based method are more
scattered. The figures also show that the appropriate regularization parameter $\lambda$ ranges from $1.5 \times 10^{-10}$ to $1.8 \times 10^{-10}$. Within this interval, Manhattan island is divided into $6 \sim 9$ coherent regions. {Note that the regularization parameter $\lambda$ is small. One possible reason is that the maximal components of $\hat{P}$ is already of the order $10^{-3}$.} The partition results coincide with the division of lower, midtown and upper Manhattan and provide abundant information of how the taxi trips are distributed.

\section{Conclusion}
We propose a convex programming formulation for the state aggregation problem of Markov chains. An atomic regularizer is introduced to control the nonnegative rank of the solutions.  The convex formulation is solved by minimizing a sequence of fixed-rank nonnegative factorization models using the proximal alternating linearization method.  The first-order optimality conditions of the convex and factorized problems are both investigated. They enable us to establish a criterion of whether a stationary point of the factorization model is globally optimal to the convex one. We further develop strategies to adjust the rank of the factorization. By increasing the rank, we can ensure a monotone decrease of the objective function and escape from a local minimum.  In situations where the matrix factorization has ``redundant ranks", we can compress the matrices so that the state space is reduced into a proper size.  Numerical experiments show that our  method always converges to a global solution at a rank much smaller than the ambient dimension. Investigation of statistical properties and a real-world application on the Manhattan transportation data set are also provided.

The performance of our method can be further improved in several aspects, including designing a better regularization than the atomic regularizer and developing a faster method for solving the factorized optimization model. Two particularly important topics of investigation are  (i) theoretical analysis on the convergence to the globally optimal solution, and (ii) numerical study for real and practical Markov decision processes.

\appendix

\section{Proofs}

\subsection{Proof of Lemma \ref{lemma_subgradient_chi}} \label{lemma_subgradient_chi:proof}
\begin{proof}
    The subgradient of characteristic function $\chi_{\mathcal{E}}(X)$ is the normal cone of $\mathcal{E}$, $i.e.$
    \[ \begin{aligned} & \partial \chi_{\mathcal{E}}(X) & = & \left\{G \in \mathbb{R}^{d \times d} \mid \langle G, Y-X \rangle \leq 0, \ \forall Y \in \mathcal{E} \right\} \\
    & & = & \left\{G \in \mathbb{R}^{d \times d} \mid \langle G, Z \rangle \leq 0, \ \forall \text{$Z$ $s.t.$ $Z{\bf 1}_d = {\bf 0}_d$} \right\}. \end{aligned} \]
    The $i$-th row of $G$ satisfies ${\bf z}^TG^i \leq 0$, for all ${\bf z} \in
    {\bf 1}_d^{\bot} $. By definition, we also have $-{\bf z} \in {\bf
    1}_d^{\bot} $ so that ${\bf z}^TG^i \geq 0$. Hence, ${\bf z}^TG^i = 0$ for
    all ${\bf z} \in {\bf 1}_d^{\bot}$. It indicates that $G^i = \mu_i {\bf
    1}_d$ for some scalar $\mu_i \in \mathbb{R}$. Combining the rows $G^i$
    together gives \eqref{subgradient_chi}.
\end{proof}

\subsection{Proof of Lemma \ref{lemma_subgradient_Omega}} \label{lemma_subgradient_Omega:proof}
\begin{proof}
    Although similar results have been provided in Lemma 13 in \cite{2015arXiv150607540H}, we briefly write the proof here in order to keep the paper self-contained.
    By definition, a matrix $W \in \partial \Omega(X)$ if and only if
        \[ \Omega(Z) \geq \Omega(X) + \langle W, Z - X \rangle, \quad \forall Z \in \mathbb{R}^{d \times d}, \]
    which is equivalent to
        \begin{equation}\label{Theorem1_3} \begin{aligned}
            & \langle W, X \rangle - \Omega(X) & \geq &
            \sup_{Z \in \mathbb{R}^{d \times d}}\left\{ \langle W, Z \rangle - \Omega(Z) \right\} \\
            & & = & \sup_{k \in \mathbb{R}_+} \left\{ k \left( \sup_{Z: \Omega(Z) = 1} \langle W, Z \rangle - 1 \right) \right\}.
        \end{aligned} \end{equation}
    If $\Omega^{\circ}(W) > 1$, $i.e.$ $\sup_{Z: \Omega(Z) = 1} \langle W, Z
    \rangle > 1$, one can take $k \rightarrow + \infty$, and the right hand side
    of \eqref{Theorem1_3} goes to infinity. It contradicts the finiteness of the
    left hand side. Therefore, one must have $\Omega^{\circ}(W) \leq 1$. In this case, we set $k = 0$ in order to reach the supremum in \eqref{Theorem1_3} and arrives at
    \begin{equation} \label{eqn1} \partial \Omega(X) = \left\{ W \left| \ \Omega^{\circ}(W) \leq 1,  \langle W, X \rangle \geq \Omega(X) \right.  \right\}. \end{equation}
    By the definition of $\Omega^{\circ}(X)$, if $\Omega^{\circ}(W) \leq 1$, one naturally has $\langle W, X \rangle \leq \Omega(X)$. Hence, \eqref{eqn1} indicates \eqref{subgradient_Omega}.
\end{proof}

\subsection{Proof of Theorem \ref{GlobalKKT}} \label{GlobalKKT:proof}
    \begin{proof}
        A matrix $\hat{X}$ is globally optimal for \eqref{reformM_new} if and only if
        \begin{equation}\label{Theorem1_1} 0 \in \nabla g(\hat{X}) + \partial \chi_{\mathcal{E}}(\hat{X}) + \lambda \partial \Omega(\hat{X}), \end{equation}
        where the explicit form of $\partial \chi_{\mathcal{E}}(\hat{X})$ and
        $\partial \Omega (\hat{X})$ can be derived from Lemmas \ref{lemma_subgradient_chi} and \ref{lemma_subgradient_Omega}.
        In order to represent condition \eqref{Theorem1_1} in terms of $\hat{U}$ and $\hat{V}$, we tailor the expression of $\partial \Omega(\hat{X})$ in \eqref{subgradient_Omega} to factorization $\hat{X} = \hat{U}\hat{V}^T$,
        \begin{equation}\label{Theorem1_7}
         \partial \Omega(\hat{X}) = \left\{ W  \left| \ \Omega^{\circ}(W) \leq 1, \sum_{j=1}^s \hat{U}_j^TW\hat{V}_j = \Omega(\hat{X}) \right. \right\} .
        \end{equation}
        Plugging \eqref{subgradient_chi} and \eqref{Theorem1_7} in \eqref{Theorem1_1}, we obtain that the global optimality of $\hat{X}$ is equivalent to the existence of a vector $\mu \in \mathbb{R}^d$ such that
        \begin{numcases}{}
            & $ \Omega^{\circ}(\hat{W}) \leq 1$, \label{Theorem1_18} \\
            & $ \sum_{j=1}^s \hat{U}_j^T \hat{W} \hat{V}_j = \Omega(\hat{X}), \quad j = 1, 2, \cdots, s$, \label{Theorem1_19_2}
        \end{numcases}
        where the matrix
        \begin{equation}\label{Theorem1_16} \hat{W} = \lambda^{-1}\left[ \mu {\bf 1}_d^T - \nabla g(\hat{X}) \right] \end{equation}
        is defined by $\mu$.
        It is also desirable to have $\hat{X} = \hat{U}\hat{V}^T$ as an optimal factorization with respect to $\Omega$. Therefore, we will make some adjustments to \eqref{Theorem1_19_2}, so that the new condition and \eqref{Theorem1_18} hold simultaneously if and only if $\hat{X}$ is globally optimal and $\hat{X} = \hat{U}\hat{V}^T$ is an optimal factorization.

        Under the condition $\Omega(\hat{X}) = \sum_{i=1}^s \|\hat{U}_j\|_2\|\hat{V}_j\|_2$, it follows from \eqref{Theorem1_19_2} that \begin{equation} \label{Theorem1_20}
            \sum_{j=1}^s \hat{U}_j^T \hat{W} \hat{V}_j = \sum_{j=1}^s \|\hat{U}_j\|_2 \|\hat{V}_j\|_2.
        \end{equation}
        Condition \eqref{Theorem1_18} also implies $\hat{U}_j^T \hat{W} \hat{V}_j \leq \|\hat{U}_j\|_2 \|\hat{V}_j\|_2$ for $j=1, 2, \cdots, s$. Compared with \eqref{Theorem1_20}, the inequalities need to hold as equalities. Hence, \eqref{Theorem1_19_2} can be reduced to
        \begin{equation}
            \hat{U}_j^T \hat{W} \hat{V}_j = \|\hat{U}_j\|_2 \|\hat{V}_j\|_2, \quad j = 1, 2, \cdots, s. \label{Theorem1_19}
        \end{equation}

        Conversely, if there exists $\mu \in \mathbb{R}^d$ satisfying \eqref{Theorem1_18} and \eqref{Theorem1_19}, then
        \[ \sum_{j=1}^s \|\hat{U}_j\|_2 \|\hat{V}_j\|_2 \overset{\eqref{Theorem1_19}}{=} \sum_{j=1}^s \hat{U}_j^T \hat{W} \hat{V}_j = \langle \hat{W}, \hat{X} \rangle \overset{\eqref{Theorem1_18}}{\leq} \Omega(\hat{X}). \]
        The definition of $\Omega(\hat{X})$ suggests that
        $ \Omega(\hat{X}) \leq \sum_{j=1}^s \|\hat{U}_j\|_2 \|\hat{V}_j\|_2$. Therefore,
        $\Omega(\hat{X}) = \sum_{j=1}^s \|\hat{U}_j\|_2 \|\hat{V}_j\|_2$ and
        $\hat{X} = \hat{U}\hat{V}^T$ is an optimal factorization with respect to
        $\Omega$. Consequently, one can easily derive \eqref{Theorem1_19_2} from \eqref{Theorem1_19}, so $\hat{X}$ is a global solution.

        To sum up, we have shown that \eqref{Theorem1_18} and
        \eqref{Theorem1_19} are sufficient and necessary conditions for
        $\hat{X}$ to be a global solution with an optimal factorization $\hat{X} = \hat{U}\hat{V}^T$.
        In the following paragraphs, we will rewrite these two conditions in a more explicit form.

        According to the definition of $\Omega^{\circ}$, \eqref{Theorem1_18} can be transformed into
        \begin{equation} \label{Theorem1_18_2}
            {\bf u}^T\hat{W}{\bf v} \leq 1, \quad \forall {\bf u},{\bf v} \in \mathbb{R}^d_+ \ s.t. \ \|{\bf u}\|_2 = \|{\bf v}\|_2 = 1.
        \end{equation}
        We next claim that, under \eqref{Theorem1_18}, the following statements \eqref{Theorem1_9} and \eqref{Theorem1_10} are both equivalent to \eqref{Theorem1_19}:
        \begin{equation}\label{Theorem1_9}
            \left[\hat{W}\hat{V}\right]_+ = \hat{U}{\bf diag}\left\{\frac{\|\hat{V}_j\|_2}{\|\hat{U}_j\|_2} \right\}_{j=1}^s,
        \end{equation}
        \begin{equation}\label{Theorem1_10}
            \left[\hat{W}^T\hat{U}\right]_+ = \hat{V}{\bf diag}\left\{\frac{\|\hat{U}_j\|_2}{\|\hat{V}_j\|_2} \right\}_{j=1}^s.
        \end{equation}

        Due to the symmetry of \eqref{Theorem1_9} and \eqref{Theorem1_10}, we only need to prove the equivalence between \eqref{Theorem1_9} and \eqref{Theorem1_19}.
        The condition $\Omega^{\circ}(\hat{W}) \leq 1$ implies that, for an arbitrary $j$, $\frac{\hat{U}_j}{\|\hat{U}_j\|_2}$ optimizes function ${\bf u}^T(\hat{W}\hat{V}_j)$ over the set $\{{\bf u} \in \mathbb{R}_+^{d} \mid \|{\bf u}\|_2 = 1\}$, $i.e.$,
        \begin{equation} \label{eqn2}
            \frac{\hat{U}_j}{\|\hat{U}_j\|_2} = \arg\max_{{\bf u} \in \mathbb{R}_+^{d}, \|{\bf u}\|_2 = 1} u^T(\hat{W}\hat{V}_j) = \frac{[ \hat{W}\hat{V}_j ]_+}{\| [ \hat{W}\hat{V}_j ]_+ \|_2}.
        \end{equation}
        The second equality in \eqref{eqn2} holds because $[\hat{W}\hat{V}_j]_+ \neq {\bf 0}_d$, which is derived by comparing the two sides of \eqref{Theorem1_19} and considering that $\hat{U}_j \neq {\bf 0}_d$, $\hat{V}_j \neq {\bf 0}_d$ and $\hat{U}_j \geq 0$.
        Additionally, we have
        \[ \|[\hat{W}\hat{V}_j]_+\|_2 = \left(\frac{[\hat{W}\hat{V}_j]_+}{\|[\hat{W}\hat{V}_j]_+\|_2}\right)^T\hat{W}\hat{V}_j \overset{\eqref{eqn2}}{=} \left(\frac{\hat{U}_j}{\|\hat{U}_j\|_2} \right)^T\hat{W}\hat{V}_j \overset{\eqref{Theorem1_19}}{=} \|\hat{V}_j\|_2, \]
        which gives \[[\hat{W}\hat{V}_j]_+ = \frac{\|\hat{V}_j\|_2}{\|\hat{U}_j\|_2} \hat{U}_j. \]
        Combining the vectors together, we arrive at \eqref{Theorem1_9}.
        On the other hand, when \eqref{Theorem1_9} holds,
        \[ \begin{aligned}
            & \hat{U}_j^T \hat{W} \hat{V}_j & = & \hat{U}_j^T [\hat{W}\hat{V}_j]_+ \\ & & = & \hat{U}_j^T \left( \frac{\|\hat{V}_j\|_2}{\|\hat{U}_j\|_2} \hat{U}_j \right) = \|\hat{U}_j\|_2\|\hat{V}_j\|_2, \qquad j = 1,\cdots,s.
        \end{aligned} \]
        Therefore, \eqref{Theorem1_19} and \eqref{Theorem1_9} are equivalent.

        Integrating \eqref{Theorem1_18_2}, \eqref{Theorem1_9} and \eqref{Theorem1_10} together, we complete the proof of Theorem \ref{GlobalKKT}.
        We also note that, for an arbitrary $i = 1,2,\cdots,d$, since $(\hat{U}^i)^T{\bf 1}_s = 1$, there exists $j \in \{1,2,\cdots,s\}$ such that $\hat{u}_{ij} > 0$. The vector $\mu \in \mathbb{R}^{d}$ has a closed form,
        \begin{equation}\label{mu}
            \mu_i = \left( \nabla g(\hat{X})\hat{V} \right)_{ij} + \lambda \hat{u}_{ij} \frac{\| \hat{V}_j \|_2}{\| \hat{U}_j \|_2}, \quad \text{for any $j$ such that $\hat{u}_{ij} > 0$}. \end{equation}
        The existence of $\mu$ is verified if and only if for any fixed $i$, the right hand side of \eqref{mu} takes the same value for all $j$ satisfying $\hat{u}_{ij} > 0$.
    \end{proof}

\subsection{Proof of Lemma \ref{LemmaEquivalence}} \label{LemmaEquivalence:proof}
    \begin{proof}
        It is obvious that if $(U,V)$ is feasible for \eqref{fixrM}, then $X=UV^T$ is feasible for \eqref{reformM_new}. By the definition of $\Omega(X)$, for any $U \in \mathcal{U}^{d \times s}$, $V \in \mathcal{V}^{d \times s}$,
        \[f_{\lambda}(X) = g(UV^T) + \lambda \Omega(X) \leq g(UV^T) + \lambda \sum_{j=1}^s \|U_j\|_2\|V_j\|_2 = F_{\lambda}(U,V).\]
        Taking the minimum on both sides, we have
        \begin{equation}\label{LemmaEquivalence_1}
            f_{\lambda}(\hat{X}) \leq F_{\lambda}(\hat{U},\hat{V}),
        \end{equation}
        where $\hat{X}$ and $(\hat{U},\hat{V})$ are respectively the optimal solutions to \eqref{reformM_new} and \eqref{fixrM}.

        We next aim to show that, given a global solution $\hat{X}$ to \eqref{fixrM}, there exists $(\hat{U},\hat{V}) \in \mathcal{U}^{d \times s} \times \mathcal{V}^{d \times s}$ such that $F_{\lambda}(\hat{U},\hat{V}) = f_{\lambda}(\hat{X})$ when $s$ is sufficiently large.
        In fact, the atomic set $\mathcal{A}_+$ in \eqref{definition_atomicset}
        is compact, which implies that its convex hull ${\rm
        conv}(\mathcal{A}_+)$ is also compact. To this end, the infimum in the
        definition of $\Omega(\hat{X})$ can be achieved, $i.e.$, $\hat{X} \in t {\rm conv}(\mathcal{A}_+)$ holds for $t=\Omega(\hat{X})$.
        According to Carath\'eodory's Theorem, there exist atoms $A_1, A_2, \cdots, A_{s_0} \in \mathcal{A}_+$ and parameters $c_1, c_2, \cdots, c_{s_0} > 0$ where $s_0 \leq d^2+1$, such that \[ \hat{X} =\Omega(\hat{X})\sum_{j=1}^{s_0} c_j A_j, \quad \sum_{j=1}^{s_0} c_j = 1. \]
        Here, each atom $A_j$ can be represented by $A_j = {\bf u}_j{\bf v}_j^T$ with ${\bf u}_j, {\bf v}_j \in \mathbb{R}_+^d$, $\|{\bf u}_j\|_2 = \|{\bf v}_j\|_2=1$. When $s \geq s_0$, we take
        \[ \begin{aligned}
         & \hat{U} := \Omega(\hat{X}) \left[ \begin{array}{ccccc} c_1({\bf v}_1^T{\bf 1}_d){\bf u}_1, & c_2({\bf v}_2^T{\bf 1}_d){\bf u}_2, & \cdots, & c_{s_0}({\bf v}_{s_0}^T{\bf 1}_d){\bf u}_{s_0}, & \mathbf{0}_{d \times (s-s_0)} \end{array}\right], \\
         & \hat{V} := \left[ \begin{array}{ccccc} ({\bf v}_1^T{\bf 1}_d)^{-1}{\bf v}_1, & ({\bf v}_2^T{\bf 1}_d)^{-1}{\bf v}_2, & \cdots, & ({\bf v}_{s_0}^T{\bf 1}_d)^{-1}{\bf v}_{s_0}, & V' \end{array} \right], \end{aligned} \]
        where $V'$ is any matrix in $\mathcal{V}^{d \times (s-s_0)}$.
        In this way, $(\hat{U}, \hat{V})$ is feasible for \eqref{fixrM}, $\hat{X} = \hat{U}\hat{V}^T$ and $\Omega(\hat{X}) = \sum_{j=1}^{s}\|\hat{U}_j\|_2\|\hat{V}_j\|_2$, hence
        \begin{equation}\label{LemmaEquivalence_2}
            F_{\lambda}(\hat{U},\hat{V}) = f_{\lambda}(\hat{X}).
        \end{equation}
        It results in the equivalence between \eqref{reformM_new} and \eqref{fixrM}.
    \end{proof}

\subsection{Proof of Theorem \ref{LocalKKT}} \label{LocalKKT:proof}
            \begin{proof}
            In problem \eqref{fixrM}, the feasible set consists of linear
            constraints. Hence, the regularity condition for KKT conditions is satisfied. In other words, the KKT conditions are necessary for local optimality.

            The Lagrangian function of problem \eqref{fixrM} is defined as follows:
            \[ \begin{aligned}
                & \mathcal{L}(U,V,\mu^U,\mu^V,\Theta^U,\Theta^V) = F_{\lambda}(U,V) & - & \left(\mu^U\right)^T \left(U{\bf 1}_s - {\bf 1}_d\right) - \langle \Theta^U, U \rangle \\ & & - & \left(\mu^V\right)^T \left(V^T{\bf 1}_d - {\bf 1}_s\right) - \langle \Theta^V, V \rangle,
            \end{aligned} \]
            where $\mu^U \in \mathbb{R}^{d}$, $\mu^V \in \mathbb{R}^s$, $\Theta^U, \Theta^V \in \mathbb{R}_+^{d \times s}$ are dual variables.
            Given a local solution $(\hat{U},\hat{V})$ to \eqref{fixrM}, the KKT conditions are
            \begin{numcases}{}
                & $ 0 \in \partial_{(U,V)} \mathcal{L} = \partial F_{\lambda}(\hat{U},\hat{V}) - \left[ \begin{array}{c} \mu^U {\bf 1}_s^T \\ {\bf 1}_d\left(\mu^V\right)^T \end{array} \right] - \left[ \begin{array}{c} \Theta^U \\ \Theta^V \end{array} \right]$, \label{Theorem2_1} \\
                & $ \hat{U}{\bf 1}_s = {\bf 1}_d, \quad \langle \Theta^U, \hat{U} \rangle = 0, \hat{U} \geq 0, \Theta^U \geq 0$, \label{Theorem2_2}\\
                & $ \hat{V}^T{\bf 1}_d = {\bf 1}_s, \quad \langle \Theta^V, \hat{V} \rangle = 0, \hat{V} \geq 0, \Theta^V \geq 0$.  \label{Theorem2_3}
            \end{numcases}

            For the convenience of notations, we assume that $\hat{U}_j \neq
            {\bf 0}_d$ for $j = 1,2,\cdots,s_1$ and $\hat{U}_j = 0$ for $j =
            s_1+1,s_1+2,\cdots,s$, then decompose $\hat{U}$ and $\hat{V}$ by
            \[
            \hat{U} = \begin{blockarray}{cc}
    		 {\scriptstyle s_1} & {\scriptstyle s - s_1}\\
    		\begin{block}{[cc]}
             \hat{U}_{\alpha} & {\bf 0}_{d \times (s-s_1)} \\
    		\end{block}
    		\end{blockarray}, \quad
            \hat{V} = \begin{blockarray}{cc}
    		 {\scriptstyle s_1} & {\scriptstyle s - s_1}\\
    		\begin{block}{[cc]}
             \hat{V}_{\alpha} & \hat{V}_{\beta} \\
    		\end{block}
    		\end{blockarray}.
            \]
            Similar as $\hat{U}$ and $\hat{V}$, $\mu^V$, $\Theta^U$ and $\Theta^V$ are also decomposed as
            \[
            \mu^V = \begin{blockarray}{cc}
            \begin{block}{[c]c}
              \mu^V_{\alpha} & {\scriptstyle s_1} \\
                \mu^V_{\beta} & {\scriptstyle s-s_1} \\
            \end{block}
            \end{blockarray}, \quad
            \Theta^U = \begin{blockarray}{cc}
    		 {\scriptstyle s_1} & {\scriptstyle s - s_1}\\
    		\begin{block}{[cc]}
             \Theta^U_{\alpha} & \Theta^U_{\beta} \\
    		\end{block}
    		\end{blockarray}, \quad
            \Theta^V = \begin{blockarray}{cc}
    		 {\scriptstyle s_1} & {\scriptstyle s - s_1}\\
    		\begin{block}{[cc]}
            \Theta^V_{\alpha} & \Theta^V_{\beta} \\
    		\end{block}
    		\end{blockarray}.
            \]

            In this way, any subgrandient $W \in \partial F_{\lambda}(\hat{U},\hat{V})$ can be expressed by
            {\small \begin{equation} \label{Theorem2_4}
                W = \left[
                \begin{array}{cc}
                \nabla g(\hat{X})\hat{V}_{\alpha} + \lambda \hat{U}_{\alpha}{\bf diag}\left\{ \frac{\|\hat{V}_j\|_2}{\|\hat{U}_j\|_2}\right\}_{j=1}^{s_1} &
                \nabla g(\hat{X})\hat{V}_{\beta} + \lambda G {\bf diag}\left\{\|\hat{V}_j\|_2\right\}_{j=s_1+1}^s  \\
                \left(\nabla g(\hat{X})\right)^T \hat{U}_{\alpha} + \lambda \hat{V}_{\alpha}{\bf diag}\left\{ \frac{\|\hat{U}_j\|_2}{\|\hat{V}_j\|_2}\right\}_{j=1}^{s_1}
                & {\bf 0}_{d \times (s-s_1)}
                \end{array}
                \right],
            \end{equation}}
            where $G$ is any $d$-by-($s-s_1$) matrix satisfying $ \|G_j\|_2 \leq 1$ for $ j = 1,2,\ldots,s-s_1 $.
            Plugging \eqref{Theorem2_4} in \eqref{Theorem2_1}, we have
            \begin{numcases}{}
                & $\nabla g(\hat{X})\hat{V}_{\alpha} + \lambda \hat{U}_{\alpha}{\bf diag}\left\{ \frac{\|\hat{V}_j\|_2}{\|\hat{U}_j\|_2}\right\}_{j=1}^{s_1} - \mu^U{\bf 1}_{s_1}^T - \Theta^U_{\alpha} = 0$, \label{Theorem2_5}\\
                & $\left(\nabla g(\hat{X})\right)^T \hat{U}_{\alpha} + \lambda \hat{V}_{\alpha}{\bf diag}\left\{ \frac{\|\hat{U}_j\|_2}{\|\hat{V}_j\|_2}\right\}_{j=1}^{s_1} - {\bf 1}_d \left(\mu^V_{\alpha}\right)^T - \Theta^V_{\alpha} = 0$. \label{Theorem2_6}
            \end{numcases}
            Due to the nonnegativity and complementarity of $\hat{U}_{\alpha}$ and $\Theta^U_{\alpha}$,
            \eqref{Theorem2_5} can be reduced to
            \begin{equation}\label{Theorem2_9}
                \left[\mu^U{\bf 1}_{s_1}^T - \nabla g(\hat{X})\hat{V}_{\alpha}\right]_+ = \lambda \hat{U}_{\alpha}{\bf diag}\left\{ \frac{\|\hat{V}_j\|_2}{\|\hat{U}_j\|_2}\right\}_{j=1}^{s_1}.
            \end{equation}
            Similarly, it follows from \eqref{Theorem2_6} that
            \begin{equation}\label{Theorem2_10}
                \left[{\bf 1}_d \left(\mu^V_{\alpha}\right)^T - \left(\nabla g(\hat{X})\right)^T \hat{U}_{\alpha}\right]_+ = \lambda \hat{V}_{\alpha}{\bf diag}\left\{ \frac{\|\hat{U}_j\|_2}{\|\hat{V}_j\|_2}\right\}_{j=1}^{s_1}.
            \end{equation}

            We next investigate the relationship between $\mu^U$ and $\mu^V_{\alpha}$. Multiplying the $j$-th column in \eqref{Theorem2_9} by $\hat{U}_j^T$ on the left, we have
            \begin{equation}\label{Theorem2_12} \hat{U}_j^T\left(\mu^U - \nabla g(\hat{X})\hat{V}_j \right) = \lambda\|\hat{U}_j\|_2\|\hat{V}_j\|_2, \quad j = 1,2,\cdots,s_1. \end{equation}
            From \eqref{Theorem2_10} we also have,
            \begin{equation}\label{Theorem2_13} \hat{V}_j^T\left(\mu^V_j{\bf 1}_d - \left(\nabla g(\hat{X})\right)^T\hat{U}_j\right) = \lambda \|\hat{U}_j\|_2\|\hat{V}_j\|_2, \quad j = 1,2,\cdots,s_1. \end{equation}
            Comparing \eqref{Theorem2_12} and \eqref{Theorem2_13}, and using the fact that $\hat{V}_j^T{\bf 1}_d =1$, we arrive at
            \begin{equation}\label{Theorem2_14} \hat{U}_{\alpha}^T\mu^U =  \mu^V_{\alpha}. \end{equation}
            Replacing $\mu^V_{\alpha}$ in \eqref{Theorem2_10} by $\hat{U}_{\alpha}^T\mu^U$ and combining \eqref{Theorem2_9} and \eqref{Theorem2_10} together, we complete the proof of Theorem \ref{LocalKKT}. Similar as the discussion in the proof of Theorem \ref{GlobalKKT}, $\mu^U$ can be explicitly evaluated and is unique if existing.
    \end{proof}

\subsection{Proof of Corollary \ref{corollary}} \label{corollary:proof}
    \begin{proof}
            Because $(\hat{U},\hat{V})$ is a local solution to \eqref{fixrM}, the KKT conditions \eqref{Theorem2} imply that
            \[ \mu^T \hat{U}_j - \hat{U}_j^T \nabla g(\hat{X}) V^*_j = \lambda \|\hat{U}_j\|_2 \|\hat{V}_j\|_2, \quad j = 1,2,\cdots, s, \]
            $i.e.$,
            \[ \left\langle \mu {\bf 1}_d^T - \nabla g(\hat{X}), \hat{U}_j\hat{V}_j^T \right\rangle = \lambda \|\hat{U}_j\|_2\|\hat{V}_j\|_2, \quad j = 1,2,\cdots,s. \]
            It follows from the condition $\sum_{j=1}^s \alpha_j \hat{U}_j\hat{V}_j^T = {\bf 0}_{d \times d}$ that
            \[ \sum_{j=1}^s \alpha_j \| \hat{U}_j \|_2 \| \hat{V}_j \|_2 = \lambda^{-1}\left\langle \mu {\bf 1}_d^T - \nabla g(\hat{X}), \sum_{j=1}^s \alpha_j \hat{U}_j\hat{V}_j^T \right\rangle = 0. \]
    \end{proof}

\subsection{Proof of Theorem \ref{NewColumn}} \label{NewColumn:proof}
    \begin{proof}
        It is obvious that when $\kappa > 0$ is sufficiently small, $\bar{U} \in \mathcal{U}^{d \times (s+1)}$, $\bar{V} \in \mathcal{V}^{d \times (s+1)}$, and $(\bar{U}, \bar{V})$ is feasible.  We next analyze the first-order perturbation of $F_{\lambda}$ with respect to $\kappa$. Define the difference $\Delta := {F_{\lambda}}(\bar{U},\bar{V}) - {F_{\lambda}}(\hat{U},\hat{V})$,
        \begin{equation}\label{Theorem3_1}
            \begin{aligned}
                 &\Delta &=& \left[ g\left( \bar{U}\bar{V}^T \right) - g\left( \hat{U}\hat{V}^T \right) \right] + \lambda \left\| \kappa \bar{\bf u} \right\|_2 \left\| \left( \bar{\bf v}^T {\bf 1}_d \right)^{-1} \bar{\bf v} \right\|_2 \\
                 &&& + \lambda \sum_{j=1}^s \left(\left\| {\bf diag}\left\{ {\bf 1}_d - \kappa \bar{\bf u} \right\} \hat{U}_j \right\|_2 \| \hat{V}_j \|_2 - \|\hat{U}_j\|_2  \| \hat{V}_j \|_2 \right).
            \end{aligned}
        \end{equation}
        Because
        \[ \begin{aligned}
            & \bar{U}\bar{V}^T & = & {\bf diag}\left\{ {\bf 1}_d - \kappa \bar{\bf u} \right\}\hat{X} + \left( \kappa \bar{\bf u} \right) \left[ \left( \bar{\bf v}^T {\bf 1}_d \right)^{-1} \bar{\bf v} \right]^T \\
            & &=& \hat{X} + \kappa \left\{ \left( \bar{\bf v}^T {\bf 1}_d \right)^{-1} \bar{\bf u} \bar{\bf v}^T -  {\bf diag}\left\{ \bar{\bf u} \right\}\hat{X} \right\}, \end{aligned} \]
        by the Taylor expansion of $g$ at $\hat{X}$, when $\kappa \rightarrow 0+ $,
        \begin{equation}\label{Theorem3_DiffG}
            \begin{aligned}
                  & g\left( \bar{U}\bar{V}^T \right) - g( \hat{X} ) \\
                = & \kappa \left\langle \nabla g( \hat{X} ), \left( \bar{\bf v}^T {\bf 1}_d \right)^{-1} \bar{\bf u} \bar{\bf v}^T - {\bf diag}\left\{ \bar{\bf u} \right\}\hat{X} \right\rangle + O(\kappa^2) \\
                = & \kappa \left\{\left( \bar{\bf v}^T {\bf 1}_d \right)^{-1} \bar{\bf u}^T \nabla g( \hat{X} ) \bar{\bf v} - \left\langle \nabla g( \hat{X} ), {\bf diag}\left\{ \bar{\bf u} \right\}\hat{X} \right\rangle\right\} + O(\kappa^2).
            \end{aligned}
        \end{equation}
        Under the condition $\hat{U}_j \neq {\bf 0}_d$, $j=1,2,\cdots,s$, the
        Euclidean norm is differentiable at $\hat{U}_j$. Therefore, when $\kappa \rightarrow 0+ $,
        \begin{equation}\label{Theorem3_DiffNorm}
            \begin{aligned}
                & \left\| {\bf diag}\left\{ {\bf 1}_d - \kappa \bar{\bf u} \right\} \hat{U}_j \right\|_2 - \| \hat{U}_j \|_2
                = \left\| \hat{U}_j - \kappa {\bf diag}\left\{ \bar{\bf u} \right\}\hat{U}_j \right\|_2 - \|\hat{U}_j \|_2 \\
                = & \left\langle \frac{\hat{U}_j}{\|\hat{U}_j\|_2}, - \kappa{\bf diag}\left\{ \bar{\bf u} \right\}\hat{U}_j \right\rangle + O(\kappa^2)
                = - \frac{\kappa}{\|\hat{U}_j\|_2}(\hat{U}_j)^T{\bf diag}\left\{ \bar{\bf u} \right\}\hat{U}_j + O(\kappa^2).
            \end{aligned}
        \end{equation}
        Plugging \eqref{Theorem3_DiffG} and \eqref{Theorem3_DiffNorm} into \eqref{Theorem3_1}, we obtain
        \begin{equation}\label{Theorem3_3}
            \begin{aligned}
                  &\Delta &=&
                    \kappa \left\{ \left( \bar{\bf v}^T {\bf 1}_d \right)^{-1} \bar{\bf u}^T\nabla g(\hat{X})\bar{\bf v} - \left\langle \nabla g(\hat{X}), {\bf diag}\left\{ \bar{\bf u} \right\}\hat{X} \right\rangle \right. \\ & && \left. + \lambda \left( \bar{\bf v}^T {\bf 1}_d \right)^{-1} - \lambda \sum_{j=1}^s  \frac{\| \hat{V}_j \|_2}{\|\hat{U}_j\|_2}\hat{U}_j^T{\bf diag}\left\{ \bar{\bf u} \right\} \hat{U}_j \right\} + O(\kappa^2).
            \end{aligned}
        \end{equation}
        $(\hat{U}, \hat{V})$ is a local solution to \eqref{fixrM}, so it follows from the results in Theorem \ref{LocalKKT} that there exists a Lagrangian multiplier $\mu \in \mathbb{R}^d$ such that \eqref{Theorem2} holds. Especially we have
        \[ \left[ \mu{\bf 1}_s^T - \nabla g(\hat{X}) \hat{V} \right]_+ = \lambda \hat{U} {\bf diag}\left\{\frac{\|\hat{V}_j\|_2}{\|\hat{U}_j\|_2} \right\}, \]
        which implies that
        \[ \begin{aligned}
             & \lambda \sum_{j=1}^s \frac{\|\hat{V}_j\|_2}{\|\hat{U}_j\|_2}\hat{U}_j^T{\bf diag}\{\bar{\bf u}\}\hat{U}_j =  \left\langle \lambda \hat{U} {\bf diag}\left\{ \frac{\|\hat{V}_j\|_2}{\|\hat{U}_j\|_2} \right\}_{j=1}^s, {\bf diag}\{\bar{\bf u}\}\hat{U} \right\rangle. \\
             = & \left\langle \mu{\bf 1}_s^T - \nabla g(\hat{X}) \hat{V}, {\bf diag}\{\bar{\bf u}\}\hat{U} \right\rangle = \mu^T\bar{\bf u} - \left\langle \nabla g(\hat{X}), {\bf diag}\{\bar{\bf u}\}\hat{X} \right\rangle.
        \end{aligned} \]
        Then, \eqref{Theorem3_3} can be reduced to
        \[ \begin{aligned}
            & \Delta & = & \kappa \left\{ \left( \bar{\bf v}^T {\bf 1}_d \right)^{-1} \bar{\bf u}^T\nabla g(\hat{X})\bar{\bf v} + \lambda \left( \bar{\bf v}^T {\bf 1}_d \right)^{-1} - \bar{\bf u}^T \mu \right\} + O(\kappa^2) \\
            & & = & \kappa \left( \bar{\bf v}^T {\bf 1}_d \right)^{-1} \left\{ \lambda - \bar{\bf u}^T \left[ \mu {\bf 1}_d^T - \nabla g(\hat{X}) \right] \bar{\bf v} \right\} + O(\kappa^2), \qquad (\kappa \rightarrow 0+).
        \end{aligned} \]
        Under the condition $\bar{\bf u}^T \left[\mu {\bf 1}_d^T - \nabla g(\hat{X}) \right] \bar{\bf v} > \lambda$, $\Delta < 0$ for some sufficiently small $\kappa > 0$.
        \end{proof}

\subsection{Proof of Theorem \ref{GlobalConverge}} \label{GlobalConverge:proof}
        \begin{proof}
        Recall that we rewrite the factorization \eqref{fixrM} into \eqref{PALMfunction}. According to Lemma 3 and Theorem 1 in \cite{Bolte2014}, it suffices to show that,
        \begin{enumerate}
            \item the objective function in \eqref{PALMfunction} is a KL function,
            \item Assumptions 1 and 2 in \cite{Bolte2014} hold for our problem,
            \item $\{ (U^k,V^k) \}_{k \in \mathbb{N}} $ is bounded.
        \end{enumerate}

        Because the feasible set of \eqref{fixrM} is bounded, condition 3 naturally holds. Since $F_{\lambda} \in C^2$ and $\chi_{\mathcal{U}^{d \times s}}$, $\chi_{\mathcal{V}^{d \times s}}$ are lower semi-continuous, Assumption 1 in \cite{Bolte2014} is also satisfied.

        We next verify the KL property of $\tilde{F}_{\lambda}$. $g(UV^T)$ is a polynomial function, so it is semi-algebraic. The Euclidean norm $\|\cdot\|_2$ is a semi-algebraic function, and the finite sums and products of semi-algebraic functions are also semi-algebraic, so $\sum_{j=1}^{s}\|U_j\|_2\|V_j\|_2$ is
        semi-algebraic. $\mathcal{U}^{d \times s}$ and $\mathcal{V}^{d \times
        s}$ consist of linear constraints, so they are semi-algebraic sets. Hence, their characteristic functions $\chi_{\mathcal{U}^{d \times s}}$ and $\chi_{\mathcal{V}^{d \times s}}$ are semi-algebraic. It can be concluded that $\tilde{F}_{\lambda}$ is semi-algebraic. According to {Theorem 3} in \cite{Bolte2014}, it satisfies the KL property.

        As for Assumption 2, the objective function in \eqref{PALMfunction} is nonnegative, so it is bounded below. Assumption 2(i) holds.
        Additionally, we have shown that, when $V$ is fixed, the partial gradient $\nabla_U F_{\lambda}(U,V)$ is globally Lipschitz with moduli $L_1(V)$.
        Due to the boundedness of $V$ and the continuity of $L_1(V)$, $L_1(V)$ has positive lower and upper bounds. Likewise, $\nabla_V F_{\lambda}(U,V)$ also has positive lower and upper bounds for any fixed $U$. For this reason, Assumption 2(ii)\&(iii) in \cite{Bolte2014} are satisfied. Because $F_{\lambda}(U,V) \in C^2$, Assumption 2(iv) holds.

        According to Theorem 1 in \cite{Bolte2014}, the sequence $\{(U^k,V^k)\}_{k \in \mathbb{N}}$ has finite length property and globally converges to a stationary point.
        \end{proof}

\bibliographystyle{siam}
\bibliography{report3}

\end{document}